\documentclass{techReport}
\pdfoutput=1

\usepackage{mathtools}          
\usepackage{stmaryrd}           
\usepackage{xstring}            
\usepackage{enumitem}           

\usepackage[f]{esvect}          

\usepackage{graphicx}           
\usepackage{xcolor}             

\usepackage{yhmath}             

\newtheorem{conjecture}[thrm]{Conjecture}

\newcommand{\liebra}[2]{{[}{#1},{#2}{]}}
\newcommand{\poibra}[2]{{\{}{#1},{#2}{\}}}

\newcommand{\frp}[2]{\frac{\partial {#1}}{\partial {#2}}}

\newcommand{\prodscal}[2]{#1 \cdot #2}

\newcommand{\hampath}{\text{HamPath}}

\definecolor{blue}{rgb}{0.06,0.06,0.6}

\definecolor{blueDrift}{rgb}{0.1 0.1 0.8}
\definecolor{redDrift}{rgb}{0.8 0.1 0.1}
\definecolor{greenDrift}{rgb}{0.2 0.7 0.2}

\graphicspath{{figures/}}

\makeatletter
\newcommand{\reqnomode}{\tagsleft@false\let\veqno\@@eqno}
\newcommand{\leqnomode}{\tagsleft@true\let\veqno\@@leqno}
\makeatother

\makeatletter
\newsavebox\myboxA
\newsavebox\myboxB
\newlength\mylenA

\newcommand*\xoverline[2][0.75]{%
    \sbox{\myboxA}{$\m@th#2$}%
    \setbox\myboxB\null
    \ht\myboxB=\ht\myboxA%
    \dp\myboxB=\dp\myboxA%
    \wd\myboxB=#1\wd\myboxA
    \sbox\myboxB{$\m@th\overline{\copy\myboxB}$}
    \setlength\mylenA{\the\wd\myboxA}
    \addtolength\mylenA{-\the\wd\myboxB}%
    \ifdim\wd\myboxB<\wd\myboxA%
       \rlap{\hskip 0.5\mylenA\usebox\myboxB}{\usebox\myboxA}%
    \else
        \hskip -0.5\mylenA\rlap{\usebox\myboxA}{\hskip 0.5\mylenA\usebox\myboxB}%
    \fi}
\makeatother

\newcommand{\nbSet}[1]{\mathbb{#1}}
\newcommand{\setPositive}{\mathbb{+}}
\newcommand{\setNegative}{\mathbb{-}}
\newcommand{\setStar}{*}
\newcommand{\R}{\nbSet{R}} 

\newcommand{\setDeco}[2]{
    \IfEqCase{#2}{
        {b}{\xoverline{\nbSet{#1}}}
        {s}{\nbSet{#1}^{\setStar}}
        {sb}{\xoverline{\nbSet{#1}}^{\setStar}}
        {pb}{\xoverline{\nbSet{#1}}_{\setPositive}}
        {n}{\nbSet{#1}^{\phantom{\setStar}}_{\setNegative}}
        {p}{\nbSet{#1}^{\phantom{\setStar}}_{\setPositive}}
        {sn}{\nbSet{#1}^{\setStar}_{\setNegative}}
        {sp}{\nbSet{#1}^{\setStar}_{\setPositive}}
    }
}

\newcommand{\Rs}{ \ensuremath{\setDeco{R}{s}} }
\newcommand{\Rp}{ \ensuremath{\setDeco{R}{p}} }
\newcommand{\Rsp}{ \ensuremath{\setDeco{R}{sp}} }

\newcommand{\Ball}{\mathrm{B}}                     
\newcommand{\BallClosed}{\xoverline{\mathrm{B}}}     
\newcommand{\Sphere}{\mathrm{S}}

\newcommand{\MBall}{\mathbb{B}}                     
\newcommand{\MBallClosed}{\xoverline{\mathbb{B}}}     
\newcommand{\MSphere}{\mathbb{S}}
\newcommand{\Cut}{\mathrm{Cut}}

\newcommand{\Split}{\mathrm{Split}}
\newcommand{\T}{\mathbb{T}}
\newcommand{\D}{\mathrm{D}}

\newcommand{\WF}{\mathbb{W}}
\newcommand{\tinj}{t_{\mathrm{inj}}}
\newcommand{\tcut}{t_{\mathrm{cut}}}
\newcommand{\tvor}{t_{\mathrm{vor}}}

\newcommand{\abs}[1]{\lvert#1\rvert}
\newcommand{\norme}[1]{\lVert#1\rVert}
\newcommand{\normeStyle}[1]{\left\lVert#1\right\rVert}

\DeclareMathOperator{\sign}{sign}
\DeclareMathOperator{\im}{Im}
\DeclareMathOperator{\atanh}{atanh}

\newcommand{\intervalle}[4]{\mathopen{#1}#2
                                \mathclose{}\mathpunct{},#3
                                \mathclose{#4}}
\newcommand{\intervalleff}[2]{\intervalle{[}{#1}{#2}{]}}
\newcommand{\intervalleof}[2]{\intervalle{(}{#1}{#2}{]}}
\newcommand{\intervallefo}[2]{\intervalle{[}{#1}{#2}{)}}
\newcommand{\intervalleoo}[2]{\intervalle{(}{#1}{#2}{)}}

\newcommand{\enstq}[2]{\left\{#1\mathrel{}\middle|\mathrel{}#2\right\}}

\newcommand{\umax}{u_\mathrm{max}}  
\newcommand{\uset}{\mathcal{U}}     
\newcommand{\udom}{\mathbf{U}}      
\newcommand{\Aset}{\mathcal{A}}     
\newcommand{\Tudom}{\mathcal{D}}    
\newcommand{\xdom}{M}               
\newcommand{\Fcal}{\mathcal{F}}     

\newcommand{\veps}{\varepsilon}
\newcommand\vphi{\varphi}
\newcommand{\Htrue}{\mathbf{H}}

\newcommand{\tin}{t_\mathrm{in}}
\newcommand{\tout}{t_\mathrm{out}}

\usepackage{circuitikz}
\usetikzlibrary{arrows,shapes,backgrounds,patterns}
\usepackage{pifont}
\usepackage{tikzgraphicx}
\tikzstyle{every picture}+=[remember picture]
\tikzstyle{na} = [baseline=-.5ex]
\makeatletter
\newcommand{\gettikzxy}[3]{%
  \tikz@scan@one@point\pgfutil@firstofone#1\relax
  \edef#2{\the\pgf@x}%
  \edef#3{\the\pgf@y}%
}
\makeatother

\newcommand\gammat{\tilde{\gamma}}
\newcommand\rt{\tilde{r}}
\newcommand\zt{\tilde{z}}
\newcommand\thetat{\tilde{\theta}}
\newcommand\pt{\tilde{p}}

\newcommand\rb{\bar{r}}
\newcommand\tb{\bar{t}}
\newcommand\thetab{\bar{\theta}}
\newcommand\pb{\bar{p}}
\newcommand\mub{\bar{\mu}}
\newcommand\lb{\bar{\lambda}}

\newcommand{\pre}{p_{r,e}}
\newcommand{\re}{r_e}

\newcommand{\diff}{\,{\rm d}}

\newcommand{\fonction}[5]{
    \begin{array}{rcll}
        #1 \colon & #2 & \longrightarrow & #3 \\
                  & #4 & \longmapsto     & #5
    \end{array}
}
\newcommand{\function}[5]{\fonction{#1}{#2}{#3}{#4}{#5}}

\begin{document}
\title{A Zermelo navigation problem with a vortex singularity}\thanks{This research is supported by the French Ministry for Education, Higher Education and Research. 
}
\author{Bernard Bonnard}\address{Inria Sophia Antipolis and Institut de Math\'ematiques de Bourgogne, UMR CNRS 5584, 9 avenue Alain Savary,
21078 Dijon, France: \texttt{Bernard.Bonnard@u-bourgogne.fr}.}
\author{Olivier Cots}\address{Toulouse Univ., INP-ENSEEIHT-IRIT, UMR CNRS 5505, 2 rue Camichel, 31071 Toulouse, France:
\texttt{olivier.cots@irit.fr}.} 
\author{Boris Wembe}\address{Toulouse Univ., IRIT-UPS, UMR CNRS 5505, 118 route de Narbonne, 31062 Toulouse, France:
\texttt{boris.wembe@irit.fr}.}
\date{\today}
\begin{abstract} 
Helhmoltz-Kirchhoff equations of motions of vortices of an incompressible fluid in the plane define a dynamics with singularities
and this leads to a Zermelo navigation problem describing the ship travel in such a field where the control is the heading angle.
Considering one vortex, we define a time minimization problem which can be analyzed with the technics of
geometric optimal control combined with numerical simulations, the geometric frame being the extension of Randers metrics in
the punctured plane, with rotational symmetry. 
Candidates as minimizers are
parameterized thanks to the Pontryagin Maximum Principle as extremal solutions of a Hamiltonian vector field. We analyze the time minimal
solution to transfer the ship between two points where during the transfer the ship can be either in a strong current region in the 
vicinity of the vortex or in a weak current region. The analysis is based on a micro-local classification of the extremals using mainly
the integrability properties of the dynamics due to the rotational symmetry.
The discussion is complex and related to the existence of an isolated extremal (Reeb) circle due to the vortex singularity.
The explicit computation of cut points where the extremal curves cease to be optimal is given
and the spheres are described in the case where at the initial point the current is weak.
\end{abstract}

\begin{resume} 
Les \'equations d'Helhmoltz-Kirchhoff pour le mouvement tourbillonaire d'un fluide incompressible dans le plan d\'efinissent une dynamique
hamiltonienne \`a singularit\'es localis\'ees aux tourbillons. Cela conduit \`a d\'efinir un probl\`eme de Zermelo d\'ecrivant le mouvement
d'un navire o\`u le contr\^ole est l'angle de cap et le crit\`ere \`a minimiser est le temps de transfert entre deux points
du plan. Dans cet article,
on se limite au cas d'un seul tourbillon localis\'e en z\'ero et le probl\`eme est analys\'e avec les techniques du contr\^ole optimal 
g\'eom\'etrique combin\'ees \`a des simulations num\'eriques, le contexte g\'eom\'etrique \'etant l'extension des m\'etriques de Randers
dans le plan poss\'edant une sym\'etrie de r\'evolution. On doit en effet consid\'erer le cas d'un courant faible, mais aussi d'un courant
fort localis\'e au voisinage du tourbillon. Les trajectoires candidates \`a minimiser le temps sont param\'etr\'ees en utilisant le Principe
du Maximum de Pontryaguine comme des extr\'emales solutions d'un syst\`eme hamiltonien dont les projections sur l'espace d'\'etat sont les
g\'eod\'esiques. L'analyse du probl\`eme optimal repose sur la classification micro-locale des solutions extr\'emales utilisant l'int\'egrabilit\'e
de la dynamique. La discussion est complexe et repose sur l'existence d'un cercle g\'eod\'esique dit de Reeb, cons\'equence de la singularit\'e
tourbillonaire. 
La discussion est compl\'et\'ee par l'\'evaluation des points de coupure, les points o\`u les extr\'emales cessent d'\^etre optimales.
Les sph\`eres sont d\'ecrites dans le cas d'un \'etat initial \`a courant faible.
\end{resume}

%
\maketitle


%

\section{Introduction}

Helhmoltz and Kirchhoff originated the model of the displacement of particles in a two dimensional fluid, see \cite{Helmholtz, Kirchhoff:1876} for
the original articles and \cite{Newton:2001,ArnoldKhesin:1998,Saffman:1992} for a modern presentation of Hamiltonian dynamics. In this model, the vorticity of the fluid is concentrated at points $z_i \coloneqq ( x_i, y_i)$, $i = 1,\cdots,N$, with \emph{circulation} parameters $k_i$ and the configuration space is $\R^{2N}$ with
coordinates $(x_1, y_1, \cdots, x_N, y_N)$ endowed with the symplectic form $\omega \coloneqq \sum_{i=1}^{N} k_i\, \diff y_i \wedge \diff x_i$.
The dynamics is given by the Hamiltonian canonical equations
\begin{equation}
k_i\, \dot{x}_i = \frac{\partial H}{\partial y_i}, \quad k_i\, \dot{y}_i = -\frac{\partial H}{\partial x_i},
\label{eq:dynamics_intro}
\end{equation}
$1 \le i \le N$, where the Hamiltonian function $H$ is
\begin{equation}
H \coloneqq -\frac{1}{\pi} \sum_{i < j} k_i \, k_j \, \ln \norme{z_i - z_j},
\label{eq:Hamiltonian_intro}
\end{equation}
where $\norme{z_i - z_j}$ is the Euclidean distance.
In this article, we consider a motionless single vortex given by eq.~\eqref{eq:dynamics_intro}, that we fix at the origin of the reference frame.
This corresponds to set $z_1=(x_1, y_1)$ at $(0, 0)$. We consider a particle as a (point)
vortex with zero circulation, setting $k_2 = 0$ with $N=2$, under the influence of the \emph{current} generated by the vortex and
given by the vector field with a singularity at the origin. The current is defined by \eqref{eq:dynamics_intro} and denoted, omitting indices, 
$F_0(x, y) \coloneqq X_1(x, y) \frac{\partial}{\partial x} + X_2(x, y) \frac{\partial}{\partial y}$, with $z \coloneqq (x, y)$.
These classical notations being not adapted for later considerations, we will denote by $x \coloneqq (x_1, x_2)$ the position of the particle
(instead of $z=(x,y)$) and by $k$ the circulation parameter of the vortex.
In the following we will consider a (Zermelo) time minimization problem.   
%

To define a Zermelo navigation problem, following \cite{Caratheodory, Zermelo:1931} and see \cite{BrysonHo:1975} for the optimal control frame related
to Zermelo's problems, we consider the particle as the \emph{ship} of the navigation problem and the control is defined by $\alpha$ 
the heading angle of the ship axis. Hence the control field is given by 
$u \coloneqq \umax (\cos \alpha, \sin \alpha)$ where $\umax$ is the maximal amplitude and this leads to a control system written as:
\begin{equation}
\frac{\diff x}{\diff t} = F_0(x) + u_1\, F_1(x) + u_2\, F_2(x),
\label{eq:control_system_intro}
\end{equation}
with $F_1 \coloneqq {\partial}/{\partial x_1}$, $F_2 \coloneqq {\partial}/{\partial x_2}$, $x = (x_1, x_2)$ and
$u \coloneqq (u_1, u_2)$ bounded by $\norme{u} \le \umax$.
We then consider the associated time minimal control problem to transfer the ship from an initial configuration $x_0$ to a target $x_f$,
where $x_0$ and $x_f$ are two points of the punctured plane $\R^2 \setminus \{0\}$. By a rescaling we can assume $\umax = 1$ and
denoting by $g$ the Euclidean metric on the plane, 
$\norme{u} \le 1$ bounds the control amplitude by $1$ and we have two cases: the case $\norme{F_0}_g < 1$ of \emph{weak current}
versus the case $\norme{F_0}_g > 1$ of \emph{strong current}. In the weak case, the time minimal problem defines a Randers metric in the 
plane, which is a specific Finsler metric, see \cite{BaoRoblesShen:2004} for this geometric frame. In the neighborhood of the vortex we have $\norme{F_0}_g > 1$, hence, due to the singularity we have a non trivial extension of the classical case.

A neat treatment of the historical Zermelo navigation was made by \cite{Caratheodory,Zermelo:1931} and their study is an important inspiration for
our work. Optimal control, with the Hamiltonian formulation coming from the Pontryagin Maximum
Principle \cite{Pontryagin:1962}, forms the frame that we shall use in our analysis, combined with recent development concerning Hamiltonian dynamics
to deal with $N$ vortices or $N$ bodies dynamics, see \cite{MeyerHallOffin:2009}.
An intense research activity was lead by H.~Poincar\'e on the dynamics of such systems \cite{Poincare:1952} to compute periodic trajectories
avoiding collisions and such techniques lead to the concept of choregraphy developed by \cite{CGMS:2002} for the $N$-body problem
and \cite{CDG:2018} for the $N$-vortex system, showing the relations between both dynamics in the Hamiltonian frame \cite{MeyerHallOffin:2009}.
From the control point of view, there is a lot of development related to space navigation for the $N$-body problem, see \cite{BFT:2006}, valuable in our
study for ship navigation in the $N$-vortex problem. Optimal control problem in this area was developed for this purpose in relation with 
Finsler geometry, see \cite{BonnardSugny:2010} or \cite{Serres:2006} for a general setting in the planar case.
More general results about vortex control may be found in \cite{protas:2008,Vainchtein:2003} for instance.

From the time minimal point of view, using the Maximum Principle, we lift the control dynamics \eqref{eq:control_system_intro}
by defining the pseudo-Hamiltonian 
\[
H(x,p,u) \coloneqq H_0(x,p) + \sum_{i=1}^{2} u_i \, H_{i}(x,p)
\]
where $H_i(x,p) \coloneqq p\cdot F_i(x)$, $i=1,2,3$, are the Hamiltonian lifts of $F_0$, $F_1$ and $F_2$ and the Maximum Principle leads
to analyze the extremal curves solution of the Hamiltonian vector field 
\[
\vv{\Htrue} \coloneqq \frac{\partial\Htrue}{\partial p} \frac{\partial}{\partial x} - 
\frac{\partial\Htrue}{\partial x} \frac{\partial}{\partial p},
\]
defined by the maximized Hamiltonian $\Htrue(x,p) \coloneqq \max_{\norme{u}\le 1} H(x,p,u)$.
%
In the context of geometric optimal control, see \cite{AgrachevSachkov:2004, Jurdjevic:1997} for a general reference, time minimal solutions are obtained by a
\emph{micro-local} analysis of $\vv{\Htrue}$ combined with the computation of the \emph{cut point} along extremal curves, that is the first point
where an extremal curve ceases to be globally optimal. Fixing the initial point $x(0) = x_0$ of the navigation problem, the set of such points is called
the \emph{cut locus}. 
Even in the Riemannian geometry, the determination of the cut locus is a very complicated problem and there are only a few results. A
major recent contribution concerns the case of ellipsoids solving the Jacobi conjecture \cite{ItohKiyohara:2004}. Parallel development were obtained in the frame of
space navigation where geometric analysis is combined with numerical methods, see again \cite{BonnardSugny:2010} for a general reference for such
contributions. In the frame of Zermelo navigation problem with a small current, called Randers metrics, some results were obtained recently
for sphere of revolutions \cite{HKS:2018}.

Our aim is to extend those results for the navigation problem, with a single vortex, which combines different new phenomena in particular
the existence of a singularity localized at the vortex position which leads to the strong current case and the need of extending the Finsler case
\cite{BaoChernShen, BaoRoblesShen:2004}.
The ultimate goal is to analyze the regularity of the value function: $x_f \mapsto V(x_0,x_f,\mu)$ where $V(x_0,x_f,\mu)$ is the minimum time
from $x_0$ to reach the point $x_f$ of the punctured plane, in presence of a vortex with a circulation $k \coloneqq 2\pi\mu$.
The parameter $\mu$ is introduced later for practical convenience. The regularity of $x_f \mapsto V(x_0,x_f,\mu)$ is analyzed, in particular, in relation with
Legendrian and Lagrangian singularities \cite{GoryunovZakalyukin:2006} associated to the Hamiltonian dynamics $\vv{\Htrue}$.
Regularity of the value function in relation with the Hamilton-Jacobi-Bellman equation leads to sufficient global optimality conditions, see the
seminal reference \cite{Boltyanskii:1966}.

The organization of this article is the following. In Section \ref{sec:existence_PMP}, we present the existence theorem to transfer
in minimum time any two points of the punctured plane. We state the Maximum Principle to parameterize the minimizers as extremal curves of a smooth
Hamiltonian vector field. This leads to define a \emph{shooting method} used to compute candidates as minimizers. Extremal curves are classified using
generic assumptions into \emph{hyperbolic} and \emph{abnormal} curves, candidate as time minimizing curves and \emph{elliptic} curves candidate as
time maximizers. \emph{Conjugate points} where an extremal curve ceases to be optimal for the $\xCn{1}$-topology are calculated numerically and leads
to the conclusion of the absence of such points, hence, the optimality problem boils down to compute cut points in the case of an empty
\emph{conjugate locus}. The final Section \ref{sec:microlocal_properties_value_function} is the main contribution of this article with
the existence result. Thanks to the
integrability of the extremal flow due to the rotational symmetry, the micro-local classification of extremals is presented. The two important
phenomena is the existence of abnormal minimizers and of a single extremal circle trajectory called a \emph{Reeb circle}. Using this
classification, the cut points can be computed along any extremal to determine the time minimal value function, combining geometric
analysis and numerical simulations using the \hampath\ code.
We present in details the case where at the initial point the current is weak.
This gives a nontrivial extension of the Finsler situation.

\section{Existence results and Pontryagin Maximum Principle}
\label{sec:existence_PMP}

\subsection{Existence of time minimal solutions}
%
\label{sec:existence_results}


We consider a single vortex centered in the reference frame and thus the control system of our Zermelo navigation problem is given by
\begin{equation*}
\dot{x}_1(t)=-\frac{k}{2\pi}\frac{x_2(t)}{r(t)^2} + u_1(t), \quad \dot{x}_2(t) = \frac{k}{2\pi}\frac{x_1(t)}{r(t)^2} + u_2(t),
\end{equation*}
with $r(t)^2 \coloneqq x_1(t)^2+x_2(t)^2$ the square distance of the ship, that is the particle, to the origin and where $k$ is the circulation of the vortex.
This control system may be written in the following form:
\begin{equation}
\dot{x}(t)= F_0(x(t))+ \sum_{i=1}^{2} u_i(t)F_i(x(t)),
\label{eq:control_system_passive_tracer}
\end{equation}
with $F_0$, $F_1$ and $F_2$ three real analytic (\ie\ $\xCn{\omega}$) vector fields, where the current (or \emph{drift}) is given by
\begin{equation}
F_0(x) \coloneqq \frac{\mu}{x_1^2+x_2^2} \left( -x_2\frac{\partial}{\partial x_1}+ x_1 \frac{\partial}{\partial x_2} \right),
\label{eq:drift}
\end{equation}
with $\mu \coloneqq {k}/{2\pi}$, and where the control fields are
$
F_1 \coloneqq {\partial}/{\partial x_1}$ 
and $F_2 \coloneqq {\partial}/{\partial x_2}.
$
Considering the polar coordinates $(x_1, x_2) \eqqcolon (r\cos\theta, r\sin\theta)$ and an adapted rotating frame for the control,
$v \coloneqq u\, e^{-i\theta}$, the control system \eqref{eq:control_system_passive_tracer} becomes
\begin{equation}
\dot{r}(t) = v_1(t), \quad \dot{\theta}(t) = \frac{\mu}{r(t)^2} + \frac{v_2(t)}{r(t)}.
\label{eq:system_polar}
\end{equation}

We give hereinafter some classical definitions and refer to \cite{BonnardChyba:2003} for more details.
We consider \emph{admissible control laws} in the set
$
\uset \coloneqq \enstq{u \colon \intervallefo{0}{+\infty} \to \udom}{u \text{ measurable}},
$
where the control domain $\udom \coloneqq \BallClosed(0,\umax) \subset \R^2$ denotes the Euclidean closed ball of radius 
$\umax > 0$ centered at the origin. Since the drift introduces a singularity at the origin, we define by
$
\xdom \coloneqq \R^2\setminus \{0\}
$
the state space,
and for any $u \in \uset$ and $x_0 \in \xdom$, we denote by $x_u(\cdot, x_0)$
the unique solution of \eqref{eq:control_system_passive_tracer} associated to the control $u$ such that $x_u(0,x_0) = x_0$.
We introduce for a time $T>0$ and an initial condition $x_0 \in \xdom$, the set $\uset_{T,x_0} \subset \uset$ of control laws
$u \in \uset$ such that the associated trajectory $x_u(\cdot,x_0)$ is well defined over $\intervalleff{0}{T}$, and we denote by
$
\Aset_{T,x_0} \coloneqq \im E_{T,x_0} 
$
the \emph{atteignable set} (or \emph{reachable set}) from $x_0$ in time $T$, where we have introduced the \emph{endpoint mapping}
\[
\function{E_{T,x_0}}{\uset_{T,x_0}}{\xdom}{u}{x_u(T, x_0).}
\]
Then, we denote by $\Aset_{x_0} \coloneqq \cup_{T\ge 0} \Aset_{T,x_0}$ the atteignable set from $x_0$.
We recall that the control system is said to be \emph{controllable from $x_0$} if $\Aset_{x_0} = \xdom$ and
\emph{controllable} if $\Aset_{x_0} = \xdom$ for any $x_0 \in \xdom$.

Now, for a given pair $(x_0, x_f) \in \xdom^2$ and some parameters $\umax \in \Rsp$ and $\mu \in \R$, we define the problem of
steering \eqref{eq:control_system_passive_tracer} in minimum time from the initial condition $x_0$ to the target $x_f$:
\leqnomode
\begin{equation}
\label{eq:minimum_time}
\tag{$P$}
V(x_0, x_f, \mu, \umax) \coloneqq \inf \enstq{T}{(T,u) \in \Tudom_{x_0} ~\text{and}~ x_u(T,x_0) = x_f},
\end{equation}
\reqnomode
where $\Tudom_{x_0} \coloneqq \enstq{(T,u) \in \intervallefo{0}{+\infty} \times \uset}{ u \in \uset_{T,x_0}}$.
We emphasize the fact that the \emph{value function} $V$ depends on the initial condition $x_0$, the target $x_f$ and the parameters
$\umax$ and $\mu$.
The first main result is the following:
\begin{thrm}
\label{thm:existence}
For any $(x_0, x_f, \mu, \umax) \in \xdom^2 \times \Rs \times \Rsp$, the problem \eqref{eq:minimum_time} admits a solution.
\end{thrm}

\begin{rmrk}
Note that when $\mu=0$, the result is clearly true in $\R^2$ but false in $\xdom = \R^2 \setminus \{0\}$.
\end{rmrk}

Up to a time reparameterization $\tau \coloneqq t\,\umax$ and a rescaling of $\mu$, one can fix $\umax=1$ and we have the relation
$
V(x_0, x_f, \mu, \umax) = V(x_0, x_f, {\mu}/{\umax}, 1) / \umax.
$
We thus fix from now $\umax=1$ and write the value function (with a slight abuse of notation)
\begin{equation}
V(x_0, x_f, \mu) \coloneqq V(x_0, x_f, \mu, 1).
\label{eq:value_function}
\end{equation}
We first prove that there exists an admissible trajectory connecting any pair of points in $\xdom$.
\begin{lmm}
The system \eqref{eq:control_system_passive_tracer} is controllable.
\label{lmm:controllability}
\end{lmm}
\begin{proof}
Let consider a pair $(x_0, x_f) \in \xdom^2$. We introduce $r_0 \coloneqq \norme{x_0}$ and $r_f \coloneqq \norme{x_f}$.
From $x_0$, we can apply a constant control $v(t) = (\pm 1,0)$ (depending on whether $r_f$ is smaller or greater than $r_0$)
until the distance $r_f$ is reached and then apply a constant control $v(t) = (0, \sign(\mu))$ until the target $x_f$ is reached.
\end{proof}

\begin{rmrk}
The controllability gives us that the value function is finite while the existence of solutions implies that the value function is lower semi-continuous.
\end{rmrk}

The existence of time-optimal solutions relies on the classical Filippov's theorem \cite[Theorem 9.2.i]{Cesari} 
and the main idea is to prove that the
problem \eqref{eq:minimum_time} is equivalent to the same problem with the restriction that the trajectories remain in a compact set.
To prove this, we need a couple of lemmas.
Let us introduce some notations for the first lemma:
for a trajectory-control pair $(x,u)$, we associate the pair $(q,v)$ with $q \coloneqq (r,\theta)$ the polar coordinates and $v=u\, e^{-i\theta}$.
We denote by $q_v(\cdot, q_0)$ the solution of \eqref{eq:system_polar} with control $v$ such that $q_v(0,q_0) = q_0$.
We define for $(\veps, R, \mu) \in \Rp\times\Rsp\times \Rs$ and $\theta_0 \in \R$, two optimization problems:

\begin{enumerate}[label=(\alph*)]
\item   The minimum time to make a complete round at a distance $R$ to the vortex:
\[
T_\theta(R,\theta_0, \mu) \coloneqq \inf \enstq{T}{(T,u) \in \Tudom_{x_0} ~\text{and}~ q_v(T,(R,\theta_0))=(R,\theta_0+s\,2\pi)},
\]
where $s \coloneqq \sign(\mu)$, $x_0 \coloneqq (R\cos \theta_0, R\sin \theta_0)$ and where the control $u$ is related to $v$ by $u = v e^{i\theta}$.
\item   The minimum time to reach the circle of radius $\veps$ from a distance $R$ to the vortex:
\[
T_r(\veps,R,\theta_0, \mu) \coloneqq \inf \enstq{T}{(T,u) \in \Tudom_{x_0} ~\text{and}~ r_v(T,(R,\theta_0))=\veps}.
\]
\end{enumerate}
Since it is clear, due to the rotational symmetry of the problem, that $T_\theta(R,\cdot,\mu)$ and $T_r(\veps,R,\cdot,\mu)$ are invariant
with respect to $\theta_0$, one can fix $\theta_0 = 0$ and set $T_\theta(R,\mu)\coloneqq T_\theta(R,0,\mu)$ and $T_r(\veps,R,\mu)\coloneqq T_r(\veps,R,0,\mu)$.
Besides, from the proof of lemma \ref{lem:T}, one can notice that $T_r$ does not depend on $\mu$, hence, one can define $T_r(\veps, R) \coloneqq T_r(\veps, R, 0)$.
Under these considerations, we have the following comparison between $T_\theta$ and $T_r$:
\begin{lmm}
\label{lem:T}
For any $(\veps, R, \mu)$ s.t. $\mu \ne 0$, $0 < R < R_\mu$ and $0 \le \veps < \veps_{\mu,R}$,
with
\[
R_\mu \coloneqq \frac{\abs{\mu}}{2\pi-1} \quad\text{and}\quad \veps_{\mu,R} \coloneqq R \left(1-\frac{2\pi R}{\abs{\mu}+R}\right),
\]
we have $0 \le \veps < \veps_{\mu,R} < R < R_\mu$ and
$
T_\theta(R, \mu) < T_r(\veps, R),
$
that is the minimum time to make a complete round at a distance $R$ to the vortex is strictly smaller than the minimum time to reach the
circle of radius $\veps < R$.
\end{lmm}
\begin{proof}
It is clear from \eqref{eq:system_polar} that $T_\theta(R, \mu)$ is given by the control
$v(t) = (0,s)$.
This gives by a simple calculation
$
T_\theta(R, \mu) = {2\pi R^2}/{(\abs{\mu}+R)}.
$
It is also clear that
$T_r (\veps, R, \mu)$ is given by $v(t)=(-1,0)$, whence
$
T_r (\veps, R, \mu) =  R - \veps = T_r(\veps, R)
$
and indeed $T_r$ does not depend on $\mu$. Fixing $\veps = 0$, we have
$
T_\theta(R, \mu) = T_r (0, R) \Leftrightarrow R = {\abs{\mu}}/{(2\pi-1)} \eqqcolon R_\mu.
$
Besides, we have
\[
T_\theta(R, \mu) < T_r (\veps, R) \Leftrightarrow \veps < R
\left( 1 - \frac{2\pi R}{\abs{\mu}+R} \right) \eqqcolon \veps_{\mu,R}
\]
but also we have
$
0 < \veps_{\mu,R} \Leftrightarrow R < R_\mu,
$
whence the conclusion.
\end{proof}


Next an admissible trajectory $x$ associated to a pair $(T,u) \in \Tudom_{x_0}$ is such that $x(T) = x_u(T,x_0)=x_f$.
Let us fix $(x_0, x_f, \mu) \in \xdom^2 \times \Rs$ and introduce $r_0 \coloneqq \norme{x_0}$ and $r_f \coloneqq \norme{x_f}$.
Then:
\begin{lmm}
\label{lem:admissible_traj}
%
There exists $\veps > 0$ such that any optimal trajectory is contained in $\xdom \setminus \BallClosed(0,\veps)$.
\end{lmm}

\begin{proof}
Let consider $(\veps, R)$ s.t. $0 < R < \min\{R_\mu, r_0, r_f\}$ and $0 < \veps < \veps_{\mu,R}$.
Let us recall that $\veps < \veps_{\mu,R} < R$ since $R < R_\mu$ and consider an admissible trajectory $x$ intersecting $\BallClosed(0, \veps)$
and associated to a pair denoted $(T,u)$.
Then, there exists two times $0 < t_1 \le t_2 < T$ s.t. $x(\intervalleff{t_1}{t_2}) \subset \BallClosed(0, \veps)$.
Since $0 < \veps < R < \min\{r_0,r_f\}$,
there exists also
$
0 < \tin < t_1 \le t_2 < \tout < T
$
s.t. $x(\intervalleff{\tin}{\tout}) \subset \BallClosed(0, R)$
and s.t. $x(\tin)$ and $x(\tout)$ belong to $\partial \BallClosed(0, R) = \Sphere(0,R)$, the sphere of radius $R$ centered at the origin.
In all generality, one can assume that
$
\forall\, t\in \intervallefo{0}{\tin}\cup\intervalleof{\tout}{T}$, $x(t)\notin\BallClosed(0,R).
$
Let consider the circular arc from $x(\tin)$ to $x(\tout)$ obtained with a control $v = (0, s)$, $s \coloneqq \sign(\mu)$,
realized in a time denoted $\tau>0$. It is clear that
$
\tau \le T_\theta(R, \mu)
$
since $T_\theta(R, \mu)$ is the time to make a circular arc of angle
$2\pi$. It is also clear from lemma \ref{lem:T} and from the definition of $T_r$ that
$
\tau \le T_\theta(R, \mu) < T_r(\veps, R) \le \tout-\tin.
$
Let us replace the part $x(\intervalleff{\tin}{\tout})$
by the circular arc. Then, the new trajectory associated to the pair denoted $(T',u')$ is still admissible
and is strictly better than $x$ since $T' = T - (\tout-\tin) + \tau < T$.
Moreover, this new trajectory is by construction contained in $\xdom \setminus \BallClosed(0,\veps)$. Whence the conclusion.
\end{proof}

\begin{figure}[ht!]
\def\scalingExistence{0.97}
\centering
\begin{tikzpicture}[scale=\scalingExistence]

\pgfmathsetmacro{\Rvortex}{0.08}
\pgfmathsetmacro{\Rorbite}{0.4}
\pgfmathsetmacro{\R}{2}
\pgfmathsetmacro{\Rpoint}{0.07} node[yshift=1];

\draw[color=red] (0,0) circle (\Rorbite);
\draw[color=cyan] (0,0) circle (\R);

\coordinate (x0) at ( 4,0.3);
\coordinate (xf) at (-4,1);
\coordinate (xi) at ( 2,0);
\coordinate (xo) at (-1.9021,0.6180);

\coordinate(P0) at ( 3.3, 0.4);
\coordinate(P1) at ( 2.7, 0.3);
\coordinate(P2) at ( 0.8,-0.4);
\coordinate(P3) at (-0.8,-0.3);
\coordinate(P4) at (-2.5, 1.0);
\coordinate(P5) at (-3.2, 1.2);



\draw[black, thick] (x0) .. controls (P0) and (P1) .. (xi) node[midway, sloped] {\scriptsize \ding{228}};
\draw[black, thick] (xi) .. controls (P2) and (P3) .. (xo) node[near end, sloped] {\scriptsize \ding{228}};
\draw[black, thick] (xo) .. controls (P4) and (P5) .. (xf) node[midway, sloped] {\scriptsize \ding{228}};

\draw (0,0)  to (0.4,0);
\draw (0,0)  to (0,2);
\draw (0.2,0.15) node[black]{$\varepsilon$} ;
\draw (0.0,1.0) node[black, left]{R} ;

\shade[ball color=red] (0,0) circle (\Rvortex);
\filldraw [blue] (x0) circle [radius=1.5pt] node[below] {$x_f$};
\filldraw [gray] (xi) circle [radius=1.0pt];
\filldraw [gray] (xo) circle [radius=1.0pt];
\filldraw [blue] (xf) circle [radius=1.5pt] node[below] {$x_0$};

\end{tikzpicture}
%
%
%
%
\begin{tikzpicture}[scale=\scalingExistence]

\pgfmathsetmacro{\Rvortex}{0.08}
\pgfmathsetmacro{\Rorbite}{0.4}
\pgfmathsetmacro{\R}{2}
\pgfmathsetmacro{\Rpoint}{0.07} node[yshift=1];

\draw[color=red] (0,0) circle (\Rorbite);
\draw[color=cyan] (0,0) circle (\R);

\coordinate (x0) at ( 4,0.3);
\coordinate (xf) at (-4,1);
\coordinate (xi) at ( 2,0);
\coordinate (xo) at (-1.9021,0.6180);

\coordinate(P0) at ( 3.3, 0.4);
\coordinate(P1) at ( 2.7, 0.3);
\coordinate(P2) at ( 0.8,-0.4);
\coordinate(P3) at (-0.8,-0.3);
\coordinate(P4) at (-2.5, 1.0);
\coordinate(P5) at (-3.2, 1.2);


\draw[black, thick]         (x0) .. controls (P0) and (P1) .. (xi) node[midway, sloped] {\scriptsize \ding{228}};
\draw[black, thin, dashed]  (xi) .. controls (P2) and (P3) .. (xo);
\draw[black, thick]         (xo) .. controls (P4) and (P5) .. (xf) node[midway, sloped] {\scriptsize \ding{228}};
\draw[black, thick]         (xo) arc (162:360:2) node[midway, sloped] {\scriptsize \ding{228}};

\draw (0,0)  to (0.4,0);
\draw (0,0)  to (0,2);
\draw (0.2,0.15) node[black]{$\varepsilon$} ;
\draw (0.0,1.0) node[black, left]{R} ;

\shade[ball color=red] (0,0) circle (\Rvortex);
\filldraw [blue] (x0) circle [radius=1.5pt] node[below] {$x_f$};
\filldraw [gray] (xi) circle [radius=1.0pt];
\filldraw [gray] (xo) circle [radius=1.0pt];
\filldraw [blue] (xf) circle [radius=1.5pt] node[below] {$x_0$};

\end{tikzpicture}

\caption{Illustration of the construction of a strictly better admissible trajectory. The vortex is represented by a red ball,
while the trajectories are the solid black lines.
One can see on the left, a trajectory crossing the ball of radius $\veps$.
This trajectory is replaced on the right subgraph by a strictly better admissible trajectory.}

\end{figure}
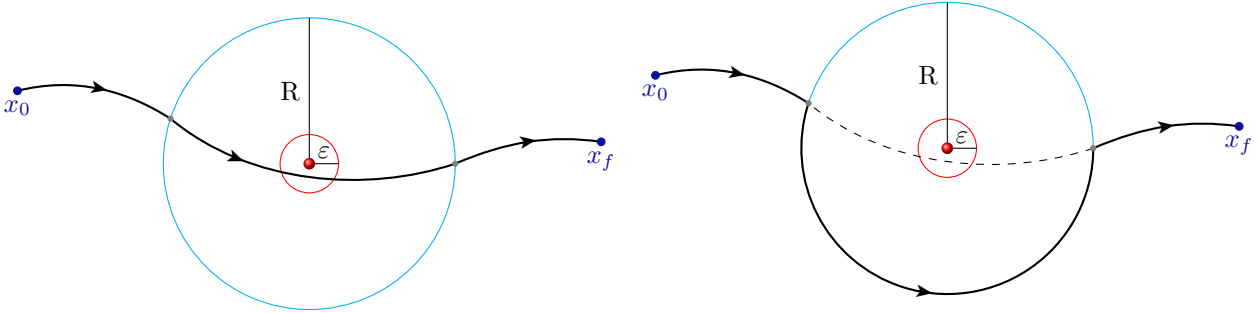

We are now in position to prove Theorem \ref{thm:existence}.
\begin{proof}[Proof of Theorem \ref{thm:existence}]
By lemma \ref{lmm:controllability}, there exists an admissible trajectory $x$. Let $T^*$ denote the first time s.t. $x(T^*) = x_f$.
Let us introduce $R_1 \coloneqq \veps$ from lemma \ref{lem:admissible_traj} and $R_2 \coloneqq r_0 + T^*$, with $r_0 \coloneqq \norme{x_0}$.
By lemma \ref{lem:admissible_traj}, the problem \eqref{eq:minimum_time} is equivalent to the same problem with the additional constraint
$R_1 \le r(t)$. Since $\dot{r}(t) = v_1(t)$ and $v_1(t) \le 1$, then for any $t \in \intervalleff{0}{T^*}$ we have $r(t) \le r_0 + T^*$.
The problem \eqref{eq:minimum_time} is thus equivalent to the same problem with the additional constraints: $R_1 \le r(t) \le R_2$.
The trajectories of the equivalent problem are contained in the compact set $\BallClosed(0,R_2) \setminus \Ball(0,R_1)$.
The result follows from the Filippov's existence theorem.
\end{proof}

\subsection{Classification of the extremal curves }
\label{sec:classification}

In this section, we recall concepts and results from \cite{BonnardSugny:2010}.
Let $x_0 \in \xdom$ and $(T,u) \in \Tudom_{x_0}$ be an optimal solution of problem \eqref{eq:minimum_time} with $x\coloneqq x_u(\cdot,x_0)$ the associated
optimal trajectory. According to the Pontryagin Maximum Principle \cite{Pontryagin:1962}, there exists an absolutely
continuous function $p \colon \intervalleff{0}{T} \to \R^2$ satisfying the adjoint equation almost everywhere over $\intervalleff{0}{T}$:
\begin{equation}
\dot{p}(t) = -\nabla_x H(x(t),p(t),u(t)),
\label{eq:adjointe}
\end{equation}
where\footnote{The standard inner product is written $a\cdot b$ or $\langle a, b \rangle$.}
$H(x,p,u) \coloneqq p \cdot (F_0(x) + u_1\, F_1(x) + u_2\, F_2(x))$ is the pseudo-Hamiltonian associated to
the problem \eqref{eq:minimum_time}.
Besides, there exists $p^0 \le 0$ such that:
\begin{equation}
\text{the pair $(p(\cdot), p^0)$ never vanishes}
\label{eq:pcout}
\end{equation}
and such that the optimal control satisfies the maximization condition almost everywhere over $\intervalleff{0}{T}$:
\begin{equation}
H(x(t), p(t), u(t)) = \max_{w \in \udom} H(x(t), p(t), w) = -p^0.
\label{eq:max}
\end{equation}

\begin{dfntn}
An \emph{extremal} is a $4$-uplet $(x(\cdot),p(\cdot),p^0,u(\cdot))$ satisfying \eqref{eq:control_system_passive_tracer} 
and \eqref{eq:adjointe}--\eqref{eq:max}. It is said \emph{abnormal} whenever $p^0=0$ and \emph{normal} whenever $p^0\ne 0$.
It is called \emph{strict} if $p(\cdot)$ is unique up to a factor.
An extremal $(x(\cdot),p(\cdot),p^0,u(\cdot))$ is called a \emph{BC-extremal} if $x(0)=x_0$ and if there is a time $T\ge 0$ s.t. $x(T)=x_f$.
\end{dfntn}

Let us introduce the Hamiltonian lifts $H_i(x,p) \coloneqq p \cdot F_i(x)$, $i=0,1,2$, 
the function $\Phi \coloneqq (H_1, H_2)$ and the \emph{switching function} $\vphi$ defined for $t \in \intervalleff{0}{T}$ by
$
\vphi(t) \coloneqq \Phi(z(t)) = p(t)$, $z(\cdot) \coloneqq (x(\cdot),p(\cdot)).
$
The maximization condition \eqref{eq:max} implies for a.e. $t \in \intervalleff{0}{T}$:
\[
u(t) = \frac{\vphi(t)}{\norme{\vphi(t)}} = \frac{p(t)}{\norme{p(t)}},
\]
whenever $\vphi(t) \ne 0$. Introducing the \emph{switching surface}
$
\Sigma \coloneqq \enstq{z \in \xdom \times \R^2}{\Phi(z) = 0} = \xdom \times \{0\}
$
and denoting $z\coloneqq (x,p) \in \xdom \times \R^2$, one can define outside $\Sigma$ the Hamiltonian:
\begin{equation}
\Htrue(z) \coloneqq H(z, \Phi(z) / \norme{\Phi(z)}) = H_0(z) + \norme{\Phi(z)}
= H_0(z) + \norme{p}.
\label{eq:True_Hamiltonian}
\end{equation}

\begin{dfntn}
An extremal $(x(\cdot),p(\cdot),p^0,u(\cdot))$ contained outside the switching surface $\Sigma$ is called of \emph{order zero}.
\end{dfntn}

Let us recall that a \emph{switching time} $0 < t < T$ is a time s.t. $\vphi(t)=0$ and s.t. for any $\veps > 0$ (small enough)
there exists a time $\tau \in (t-\veps, t+\veps) \subset \intervalleff{0}{T}$ s.t. $\vphi(\tau)\ne 0$. 
We can show that the extremals are only of order zero, and thus are smooth:
\begin{prpstn}
\label{prop:order_zero}
All the extremals are of order zero, that is there are no switching times.
\end{prpstn}
\begin{proof}
Let $(x(\cdot),p(\cdot),p^0,u(\cdot))$ be an extremal.
If there exists a time $t$ s.t. $\vphi(t)=0$, then $p(t)=0$ and we have $H(x(t),p(t),u(t)) = 0 = -p^0$ which is impossible
by \eqref{eq:pcout}.
\end{proof}

We have the standard following result:
\begin{prpstn}
The extremals of order zero are smooth responses to smooth controls on the boundary of $\norme{u} \le 1$.
They are singularities of the endpoint mapping $E_{T, x_0}$ for the $L^\infty$-topology when $u$ is restricted to the unit 
sphere $\Sphere^1$.
\end{prpstn}

For any Hamiltonian $\Htrue(z)$, resp. pseudo-Hamiltonian $H(z,u)$, we denote by $\vv{\Htrue}(z) \coloneqq (\nabla_p\Htrue(z), -\nabla_x \Htrue(z))$,
resp.  $\vv{H}(z,u) \coloneqq (\nabla_p H(z,u), -\nabla_x H(z,u))$, its associated \emph{Hamiltonian vector field}, resp. \emph{pseudo-Hamiltonian vector field}.
With these notations, we have the following classical but still remarkable fact:
\begin{prpstn}
\label{prop:True_Hamiltonian}
Let $(x(\cdot),p(\cdot),p^0,u(\cdot))$ be an extremal. Denoting $z\coloneqq(x,p)$, then, we have over $\intervalleff{0}{T}$:
\begin{equation}
\dot{z}(t) = \vv{H}(z(t),u(t)) = \vv{\Htrue}(z(t)) = \vv{H_0}(z(t)) + \left(\frac{p(t)}{\norme{p(t)}}, 0 \right),
\label{eq:True_Hamiltonian_System}
\end{equation}
that is the extremals are given by the flow of the Hamiltonian vector field associated to the maximized  Hamiltonian $\Htrue$.
\end{prpstn}
\begin{proof}
Since the extremal is of order zero, the control $t \mapsto u(t)$ is smooth and the adjoint equation \eqref{eq:adjointe} is satisfied all over
$\intervalleff{0}{T}$.
Besides, denoting (with a slight abuse of notation) $u(z) \coloneqq \Phi(z)/\norme{\Phi(z)}$, we have:
\[  
\begin{aligned}
\Htrue'(z)  &= \frac{\partial H}{\partial z}(z,u(z)) + \frac{\partial H}{\partial u}(z,u(z)) \cdot u'(z) \\
&= \frac{\partial H}{\partial z}(z,u(z)) +
\underbrace{\Phi(z)^T \left(\frac{I_2}{\norme{\Phi(z)}}-\frac{\Phi(z)\Phi(z)^T}{\norme{\Phi(z)}^3}\right)}_{=0} \cdot\, \Phi'(z)
= \frac{\partial H}{\partial z}(z,u(z)) = H'_0(z) + \left(0, \frac{p}{\norme{p}} \right).
\end{aligned}
\vspace{-2em}
\]
\end{proof}
This proposition shows the importance of the true Hamiltonian $\Htrue$ which encodes all the information we need and gives a more geometrical
point of view: we will thus refer to trajectories as \emph{geodesics}. 
%
Besides, from the maximum principle, optimal extremals are contained in the level set $\{ \Htrue \ge 0 \}$. 
Let $z(\cdot, x_0, p_0) \coloneqq (x(\cdot, x_0, p_0), p(\cdot, x_0, p_0))$ be a reference extremal curve solution of $\dot{z} = \vv{\Htrue}(z)$ with initial
condition $z(0, x_0, p_0) = (x_0, p_0)$ and defined over the time interval $\intervalleff{0}{T}$.
\begin{lmm}
One has $x(t, x_0, \lambda p_0) = x(t, x_0, p_0)$ and $p(t, x_0, \lambda p_0) = \lambda p(t, x_0, p_0)$.
\end{lmm}
Thanks to this lemma, and since $p$ never vanishes, we can fix by homogeneity $\norme{p_0} = 1$.
We thus introduce the following definition that gives us a way to parameterize the extremals of order zero.
\begin{dfntn}
We define the \emph{exponential mapping} by
\begin{equation}
\exp_{x_0}(t, p_0) \coloneqq \Pi \circ e^{t \vv{\Htrue}}(x_0, p_0),
\label{eq:exponential_mapping}
\end{equation}
where $e^{t \vv{\Htrue}}(x_0, p_0)$ is the solution at time $t$ of $\dot{z}(s) = \vv{\Htrue}(z(s))$, $z(0) = (x_0, p_0)$.
It is defined for small enough nonnegative time $t$ and we can assume that $p_0$ belongs to $\Sphere^1$.
\end{dfntn}

\begin{dfntn}
Let $z(\cdot) \coloneqq (x(\cdot), p(\cdot))$ be a reference extremal of order zero, defined on $\intervalleff{0}{T}$. Let
$\Htrue$ be the Hamiltonian defined by \eqref{eq:True_Hamiltonian}. The associated geodesic $x(\cdot)$
is called \emph{exceptional} if $\Htrue = 0$, \emph{hyperbolic} if $\Htrue > 0$ and \emph{elliptic} if $\Htrue < 0$,
along the reference extremal $z(\cdot)$.
\label{def:classification_order_zero}
\end{dfntn}

\begin{rmrk}
The previous definition is related to the more classical definition \ref{def:classification_single_input}. Even if the elliptic
geodesics are not optimal according to the PMP, they still play a role in the analysis of the optimal
synthesis, in particular in the computation of the cut locus when the current (or drift) is strong, see
Section \ref{sec:microlocal_properties_value_function}.
\end{rmrk}

In cartesian coordinates, the Hamiltonian writes 
\[
\Htrue(x_1,x_2,p_1,p_2) = \frac{\mu}{x_1^2+x_2^2} \left( -p_1\, x_2 + p_2\, x_1 \right) + \sqrt{p_1^2+p_2^2},
\]
and the extremals are solution of the following Hamiltonian system:
\begin{equation*}
\begin{aligned}
\dot{x}_1 &=          - \mu \frac{x_2}{r^2} + \frac{p_1}{\norme{p}}, \quad
\dot{x}_2  = \mu \frac{x_1}{r^2} + \frac{p_2}{\norme{p}}, \\
\dot{p}_1 &=          - \frac{\mu}{r^4} \left( 2\, x_1\, x_2\, p_1 - (x_1^2 - x_2^2)\, p_2 \right), \quad
\dot{p}_2 = \frac{\mu}{r^4} \left( (x_1^2 - x_2^2)\, p_1 - 2\, x_1\, x_2\, p_2 \right). \\
\end{aligned}
\end{equation*}
Introducing the \emph{Mathieu transformation}
\begin{equation}
\begin{pmatrix}
p_r \\ p_\theta
\end{pmatrix}
=
\begin{pmatrix}
\cos\theta & \sin\theta \\ -r \sin\theta & r \cos\theta
\end{pmatrix}
\begin{pmatrix}
p_1 \\ p_2
\end{pmatrix}
\label{eq:Mathieu_transform}
\end{equation}
then, in polar coordinates, the Hamiltonian is given by (we still denote by $p$ the covector in polar coordinates)
\[
\Htrue(r, \theta, p_r, p_\theta) = p_\theta \frac{\mu}{r^2} + \norme{p}_{r},
\]
where $\norme{p}_r \coloneqq \sqrt{p_r^2 + p_\theta^2/r^2}$. It is clear from $\Htrue$ that $\theta$ is a \emph{cyclic variable} 
and thus the problem has a symmetry of revolution and by Noether theorem, the adjoint variable $p_\theta$ is a \emph{first integral}.
This relation $p_\theta = \text{constant}$ corresponds to the \emph{Clairaut} relation on surfaces of revolution.
Hence, we can fix $\theta(0) = 0$ and consider $p_\theta$ has a parameter of the associated Hamiltonian system in polar coordinates:
\begin{equation}
\dot{r}         = \frac{p_r}{\norme{p}_r}, \quad
\dot{\theta}    = \frac{1}{r^2} \left( \mu + \frac{p_\theta}{\norme{p}_r} \right), \quad
\dot{p}_r       = \frac{p_\theta}{r^3} \left( 2 \mu + \frac{p_\theta}{\norme{p}_r} \right), \quad
\dot{p}_\theta  = 0.
\label{eq:hamiltonian_system_polar}
\end{equation}


\subsection{$\xCn{1}$-second order necessary conditions in the regular case}
\label{sec:CN2}

Since the extremals are of order zero, one can restrict $u(t)$ to the $1$-sphere $\Sphere^1$. Writing $u(t) \eqqcolon (\cos\alpha(t), \sin\alpha(t))$, we have with some abuse of notations
$H=H_0+u_1 H_1+u_2 H_2= H_0 + \cos\alpha\, H_1+ \sin\alpha\, H_2$, with $\alpha$ the new control.
Differentiating twice with respect to $\alpha$, we have
\[
\frac{\partial H}{\partial \alpha}      = -\sin\alpha\, H_1 + \cos\alpha\, H_2, \quad
\frac{\partial^2 H}{\partial \alpha^2}  = -(\cos\alpha\, H_1+ \sin\alpha\, H_2),
\]
and since $u=(\cos\alpha, \sin\alpha)=\Phi/\norme{\Phi}$ and $\Phi=(H_1,H_2)$ never vanishes along any extremal, we have
\[
\frac{\partial^2 H}{\partial \alpha^2}  = -\sqrt{H_1^2+H_2^2} < 0
\]
along any extremal.
Hence, the \emph{strict Legendre-Clebsch} condition is satisfied and we are in the \emph{regular} case \cite{BCT:2007}, but with a \emph{free} final time $T$.
\begin{dfntn}
Let $z(\cdot)$ be a reference extremal curve solution of $\dot{z} = \vv{\Htrue}(z)$ given by \eqref{eq:True_Hamiltonian_System}.
The variational equation
\begin{equation}
\dot{\wideparen{\delta z}}(t) = \vv{\Htrue}'(z(t)) \cdot \delta z(t), 
\label{eq:Jacobi}
\end{equation}
is called a \emph{Jacobi equation}. A \emph{Jacobi field} is a non trivial solution $J$ of \eqref{eq:Jacobi}. It is said to be \emph{vertical}
at time $t$ if $\delta x(t) \coloneqq \Pi'(z(t)) \cdot J(t) = 0$, where $\Pi \colon (x,p) \mapsto x$ is the standard projection.
\end{dfntn}

Let $z(\cdot, x_0, p_0) \coloneqq (x(\cdot, x_0, p_0), p(\cdot, x_0, p_0))$ with $p_0 \in \Sphere^1$, be a reference extremal curve solution
of $\dot{z} = \vv{\Htrue}(z)$ with initial condition $z(0, x_0, p_0) = (x_0, p_0)$ and defined over the time interval $\intervalleff{0}{T}$.
%
Following \cite{BCT:2007}, we make the following generic assumptions on the reference extremal in order to derive second order optimality conditions:
\begin{itemize}
\item[\textbf{(A1)}] The trajectory $x(\cdot, x_0, p_0)$ is a one-to-one immersion on $\intervalleff{0}{T}$.
\item[\textbf{(A2)}] The reference extremal is normal and strict.
\end{itemize}

\begin{dfntn}
Let $z = (x, p)$ be the reference extremal defined hereinabove. Under assumptions {(A1)} and {(A2)},
the time $0 < t_c \le T$ is called \emph{conjugate} 
if the exponential mapping is not an immersion at $(t_c, p_0)$.
The associated point $\exp_{x_0}(t_c, p_0) = x(t_c, x_0, p_0)$ is said to be \emph{conjugate to} $x_0$. We denote by $t_{1c}$ the first conjugate time.
\label{def:conjugate_time}
\end{dfntn}

The following result is fundamental, see \cite{BonnardKupka:1993}.
\begin{thrm}
Under assumptions {(A1)} and {(A2)}, the extremities being fixed, the reference geodesic $x(\cdot)$ is locally 
time minimizing (resp. maximizing) for the $L^\infty$-topology on the set of controls up to the first conjugate time
in the hyperbolic (resp. elliptic) case.
\label{thm:optimalite_Linfty}
\end{thrm}

\paragraph*{\textbf{Algorithm to compute conjugate times}}
Writing the reference trajectory  $x(t) \coloneqq x(t, x_0, p_0)$ and considering a Jacobi field $J(\cdot) \coloneqq (\delta x(\cdot), \delta p(\cdot))$
along the reference extremal, which is vertical at the initial time (\ie\ $\delta x(0) = 0$) and normalized by $p_0 \cdot \delta p(0) = 0$ 
(since $p_0$ is restricted to $\Sphere^1$), then $t_c$ is a conjugate time if and only if $t_c$ is a solution of the following equation:
\begin{equation}
t \mapsto \det(\delta x(t), \dot{x}(t)) = 0.
\label{eq:conjugate_time_eq}
\end{equation}
See \cite{BCT:2007, Co2012} for more details about algorithms to compute conjugate times in a more general setting and \cite{CCG:2010} 
for details about the numerical implementation of these algorithms into the \hampath\ software.

\subsection{The Zermelo-Carath\'eodory-Goh point of view}

From the historical point of view in the Zermelo navigation problem, Zermelo and Carath\'eodory use for the parameterization
of the geodesics the derivative of the heading angle $\alpha$ instead of the angle $\alpha$ itself, see \cite{BrysonHo:1975}.
This corresponds precisely to the so-called \emph{Goh transformation} for the analysis of singular trajectories in 
optimal control, see for instance the reference \cite{BonnardChyba:2003}. This is presented next, in relation with the problem,
to derive sufficient $\xCn{0}$-optimality conditions under generic assumptions, see \cite{BonnardKupka:1993}.

\subsubsection{Goh transformation}





Retricting to extremals of order zero, the \emph{Goh transformation} amounts to set (we use the same notation $u$ for the new control
but no confusion is possible):
\[
\dot{\alpha} = u,
\]
that is to take $\dot{\alpha}$ as the control of the ship.
Note that such a transformation transforms $L^\infty$-optimality conditions on the set of controls into $\xCn{0}$-optimality
conditions on the set of trajectories. For the geodesics computations, this amounts only to a reparameterization of extremal
curves of order zero.
Considering the vortex problem in cartesian coordinates $(x_1, x_2)$, we introduce $x \coloneqq (x_1, x_2, x_3)$,
$x_3 \coloneqq \alpha$, and the control system becomes
\[
\dot{x} = F(x) + u\, G(x),
\]
with
\[
F(x) \coloneqq
\begin{pmatrix}
F_0(x_1,x_2) + \cos x_3\, F_1(x_1,x_2) + \sin x_3\, F_2(x_1,x_2) \\
0
\end{pmatrix}
\quad
\text{and}
\quad
G = \frac{\partial}{\partial x_3}.
\]
The associated pseudo-Hamiltonian reads
\[
H(x,p,u) \coloneqq \prodscal{p}{(F(x) + u\, G(x))}, \quad p \coloneqq (p_1, p_2, p_3),
\]
and relaxing the bound on the new control, the maximization condition implies $p \cdot G = 0$ along any extremal. These extremals are called
\emph{singular} and the associated control is also called singular. Let us recall how we compute the singular extremals.
First we need to introduce the concepts of Lie and Poisson brackets.
The \emph{Lie bracket} of two $\xCn{\omega}$ vector fields $X$, $Y$ on an open subset $V \subset \R^n$ is computed with the convention:
\[
\liebra{X}{Y}(x) \coloneqq \frp{X}{x}(x)\,Y(x) - \frp{Y}{x}(x)\,X(x),
\]
and denoting $H_X$, $H_Y$ the Hamiltonian lifts: $H_X(z) \coloneqq \prodscal{p}{X(x)}$, $H_Y(z) \coloneqq \prodscal{p}{Y(x)}$,
with $z\coloneqq(x,p) \in V \times \R^n$, the \emph{Poisson bracket} reads:
\[
\poibra{H_X}{H_Y} \coloneqq H'_X \cdot \vv{H}_Y = \prodscal{p}{\liebra{X}{Y}(x)},
\]
where 
\[
\vv{H}_X = \frp{H}{p} \frp{}{x} - \frp{H}{x} \frp{}{p}.
\]
Differentiating twice $\prodscal{p(\cdot)}{G(x(\cdot))}$ with respect to the time $t$, one gets:
\begin{lmm}
Singular extremals $(z(\cdot),u(\cdot))$ are solution of the following equations:
\begin{align*}
H_G(z(t)) = \poibra{H_G}{H_F}(z(t)) = 0, \\[0.0em]
\poibra{\poibra{H_G}{H_F}}{H_F}(z(t)) + u(t)\, \poibra{\poibra{H_G}{H_F}}{H_G}(z(t)) = 0.
\end{align*}
If $\poibra{\poibra{H_G}{H_F}}{H_G} \ne 0$ along the extremal, then the singular control
is called of \emph{minimal order} and it is given by the dynamic feedback:
\begin{equation*}
u_s(z(t)) \coloneqq -\frac{ \poibra{\poibra{H_G}{H_F}}{H_F}(z(t))}{\poibra{\poibra{H_G}{H_F}}{H_G}(z(t))}.
\end{equation*}
\end{lmm}
Plugging the control $u_s$ in feedback form into the pseudo-Hamiltonian leads to define a true Hamiltonian denoted
$\Htrue_s(z) \coloneqq H(z, u_s(z))$, and one has:
\begin{lmm}
Singular extremals of minimal order are the solutions of
$
\dot{z}(t) = \vv{\Htrue}_s(z(t)),
$
with the constraints $H_G(z(t)) = \poibra{H_G}{H_F}(z(t)) = 0$.
\end{lmm}

\subsubsection{The case of dimension 3 applied to the Zermelo problem}

Consider the following affine control system: $\dot{x} = F(x) + u\, G(x)$, with $u \in \R$ and $x \in \R^3$,
where $F$, $G$, are $\xCn{\omega}$ vector fields.
Let $z(\cdot) \coloneqq (x(\cdot), p(\cdot))$ be a reference singular extremal curve on $\intervalleff{0}{T}$.
We assume the following:
\begin{itemize}
\item[\textbf{(B1)}] The reference geodesic $t \mapsto x(t)$ is a one-to-one immersion on $\intervalleff{0}{T}$.
\item[\textbf{(B2)}] $F$ and $G$ are linearly independent along $x(\cdot)$.
\item[\textbf{(B3)}] $G$, $\liebra{G}{F}$, $\liebra{\liebra{G}{F}}{G}$ are linearly independent along $x(\cdot)$.
\end{itemize}

From {(B3)}, $p$ is unique up to a factor and the geodesic is strict and moreover $u_s$ can be computed
as a true feedback:
\[
u_s(x) = -\frac{D'(x)}{D(x)},
\]
where we denote:
\[
\begin{aligned}
D   & \coloneqq \det( G, \liebra{G}{F}, \liebra{\liebra{G}{F}}{G}  ), \\
D'  & \coloneqq \det( G, \liebra{G}{F}, \liebra{\liebra{G}{F}}{F}  ).
\end{aligned}
\]
Moreover, let us introduce $D'' \coloneqq \det( G, \liebra{G}{F}, F  )$. In the vortex problem, one has:
\[
\dot{\alpha} = -\frac{D'(x)}{D(x)}.
\]
In our problem with the Goh extension, one orients $p(\cdot)$ using the convention of the maximum principle:
$\prodscal{p(t)}{F(x(t))} \ge 0$ on $\intervalleff{0}{T}$ and we introduce the following definition consistent with definition \ref{def:classification_order_zero}:
\begin{dfntn}
Under assumptions {(B1)}, {(B2)} and {(B3)}, a geodesic is called:
\begin{itemize}
\item \emph{hyperbolic} if $D D'' > 0$,
\item \emph{elliptic} if $D D'' < 0$,
\item \emph{abnormal} (or \emph{exceptional}) if $D'' = 0$.
\end{itemize}
\label{def:classification_single_input}
\end{dfntn}
Note that the condition $D D'' \ge 0$ amounts to the \emph{generalized Legendre-Clebsch} condition
\[
\frac{\partial}{\partial u}\frac{\diff^2}{\diff t^2} \frac{\partial H}{\partial u}(z(t))\ge 0
\]
and according to the \emph{higher-order maximum principle} \cite{Krener:1977},
this condition is a necessary (small) time minimization condition. See \cite{BonnardKupka:1993} for the general frame
relating the optimal control problems using the Goh transformation and applicable to our study and for the following result.
\begin{thrm}
Under assumptions {(B1)}, {(B2)} and {(B3)}, a reference hyperbolic (resp. elliptic) geodesic $x(\cdot)$
defined on $\intervalleff{0}{T}$
is time minimizing (resp. maximizing) on $\intervalleff{0}{T}$ with respect to all trajectories contained in a 
$\xCn{0}$-neighborhood of $x(\cdot)$ if $T< t_{1c}$ where $t_{1c}$ is the first conjugate time along $x(\cdot)$ as defined
by \ref{def:conjugate_time} for the projection of $x(\cdot)$ on the $(x_1,x_2)$ plane. In the exceptional case, the reference
geodesic is $\xCn{0}$-time minimizing and maximizing.
\label{thm:optimality_singular}
\end{thrm}


\subsection{Influence of the circulation}

\subsubsection{Influence of the circulation on the drift}

Denoting the drift \eqref{eq:drift} $F_0(x,\mu)$ to emphasize the role of $\mu$, one introduces for $(x,\mu) \in \xdom \times \Rs$ the set
\begin{equation}
\Fcal(x,\mu) \coloneqq \enstq{F_0(x,\mu) + \sum_{i=1}^2 u_i\,F_i(x)}{ u\coloneqq (u_1,u_2) \in \udom}.
\label{eq:drift_domain}
\end{equation}
Then, we have (noticing that if $u \in \udom$, then $-u \in \udom$):
\vspace{-1.0em}
\begin{equation*}
\begin{aligned}
0 \in \Fcal(x,\mu) \Leftrightarrow \exists\, u\coloneqq (u_1,u_2) \in \udom ~\text{s.t.}~ F_0(x,\mu) = \sum_{i=1}^2 u_i\,F_i(x) 
\Leftrightarrow \norme{F_0(x,\mu)}_g = \norme{F_0(x,\mu)} \le 1 \Leftrightarrow |\mu| \le \norme{x} = r.
\end{aligned}
\end{equation*}
This leads to introduce the following definition.
\begin{dfntn}
The drift $F_0(x,\mu)$ is said to be
\emph{weak at the point $x$} if $\norme{F_0(x,\mu)} < 1$,
\emph{strong at $x$}   if $\norme{F_0(x,\mu)} > 1$ and
\emph{moderate at $x$} if $\norme{F_0(x,\mu)} = 1$.
\end{dfntn}

\begin{rmrk}
Note that if the drift could have been weak at any point $x \in \xdom$, then we would have been in the Finslerian case \cite{BaoChernShen} with
a metric of Randers type. However, this is not possible for $\mu \ne 0$ (the case $\mu=0$ is trivially Euclidean on $\R^2$),
since in this case, we have for any $x \in \xdom \cap \BallClosed(0, \abs{\mu}) \ne \emptyset$ that the drift $F_0(x, \mu)$ is not weak.
\end{rmrk}

\begin{figure}[ht!]
\centering
\def\sizefig{0.25}
\includegraphics[width=\sizefig\textwidth]{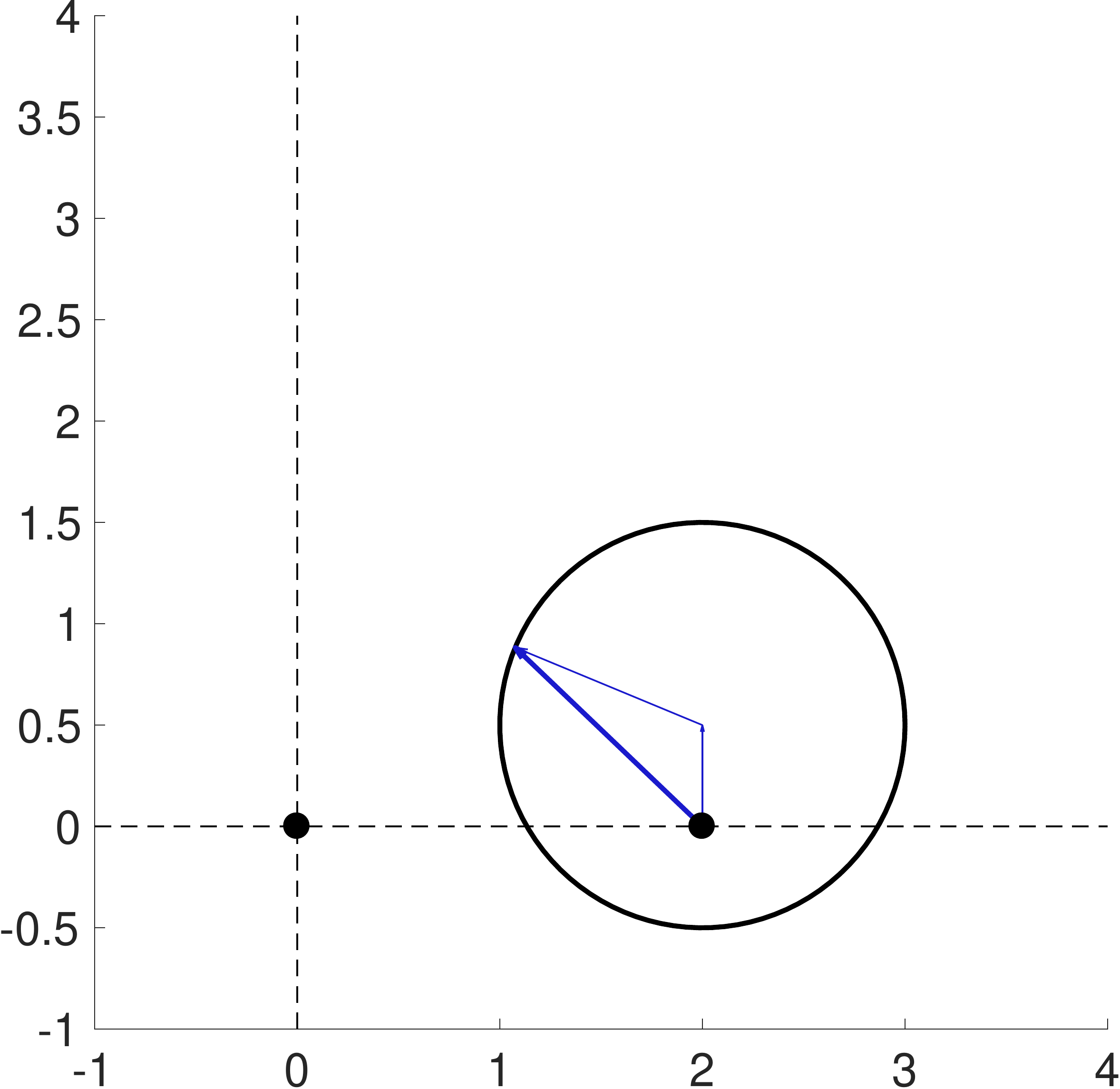}
\hspace{1,0em}
\includegraphics[width=\sizefig\textwidth]{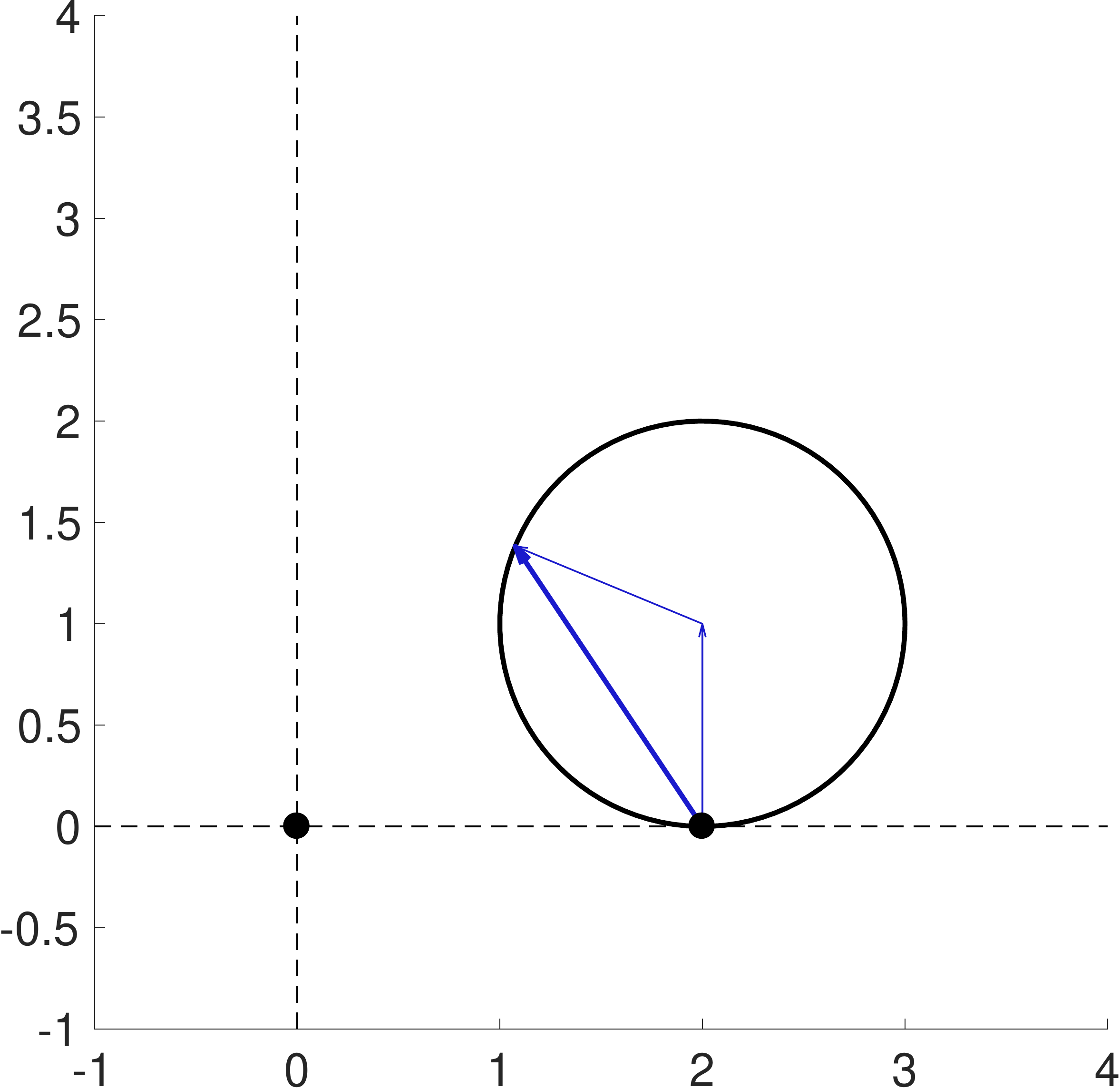}
\hspace{1,0em}
\includegraphics[width=\sizefig\textwidth]{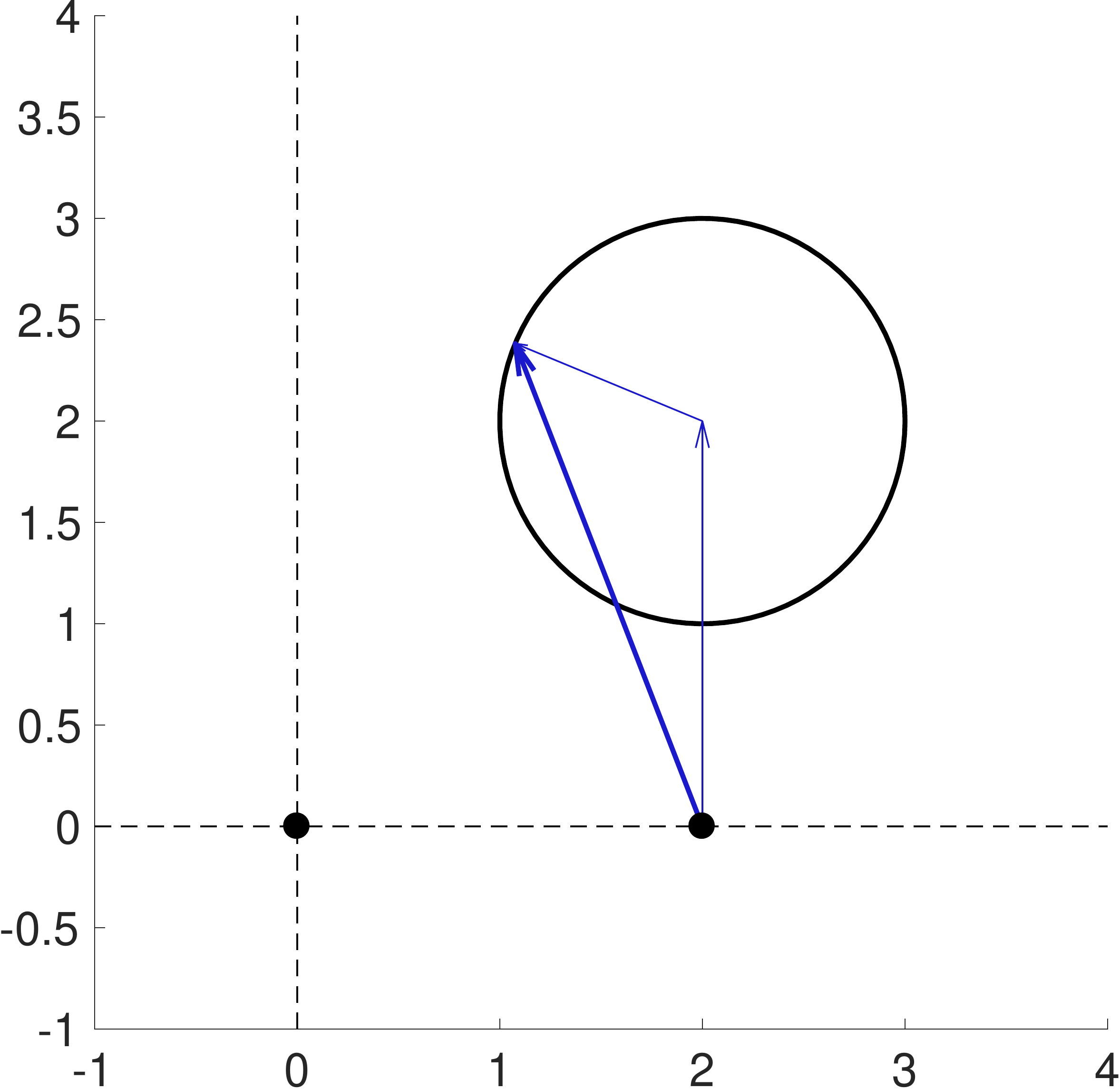}
\caption{The vortex is placed at the origin and marked by a black dot, as the initial point $x_0 \coloneqq (2,0)$ in cartesian
coordinates. The {black circle} corresponds to the set of initial directions $x_0 + \Fcal(x_0,\mu)$ and the thick {{blue vector}}
is the initial direction associated to the initial control $u(0) \coloneqq (\cos \alpha, \sin \alpha)$, with $\alpha = 7\pi/8$.
It is decomposed as the sum of the drift (oriented vertically) and the control field.
On the left, we have a weak drift ($\mu=0.5\, r_0$) at $x_0$, in the middle we have a moderate drift ($\mu=r_0$)
and on the right a strong drift ($\mu=2\, r_0$).
}
\label{fig:drift}
\end{figure}

\begin{rmrk}
One can also notice that at the initial time, the strength of the drift depends on the ratio $\abs{\mu}/r_0$.
See Figure~\ref{fig:drift} for an illustration of the different possible strengths of the drift.
\end{rmrk}

\subsubsection{Influence of the circulation on the abnormal extremals}
\label{sec:abnormal_extremals}

\begin{prpstn}
Let $(x(\cdot),p(\cdot),p^0,u(\cdot))$ be an abnormal extremal, that is $p^0=0$.
Then, the drift is strong or moderate all along the geodesic.
\end{prpstn}
\begin{proof}
Since the abnormal extremal is of order $0$, then all along the extremal we have
$
\Htrue(x(t),p(t)) = p(t) \cdot F_0(x(t), \mu) + \norme{p(t)} = 0.
$
So, the Cauchy-Schwarz inequality gives
$
\norme{p(t)} = |p(t) \cdot F_0(x(t), \mu)| \le \norme{p(t)} \norme{F_0(x,\mu)}
$
and since $\norme{p(t)} \ne 0$, the result follows.
\end{proof}

\begin{rmrk}
According to the previous proposition, the abnormal geodesics (that is the projection of the abnormal extremals on the state manifold)
are contained in the ball $\BallClosed(0, \abs{\mu})$.
\end{rmrk}

\begin{figure}[ht!]
\centering
\def\sizefig{0.32}
\includegraphics[width=\sizefig\textwidth]{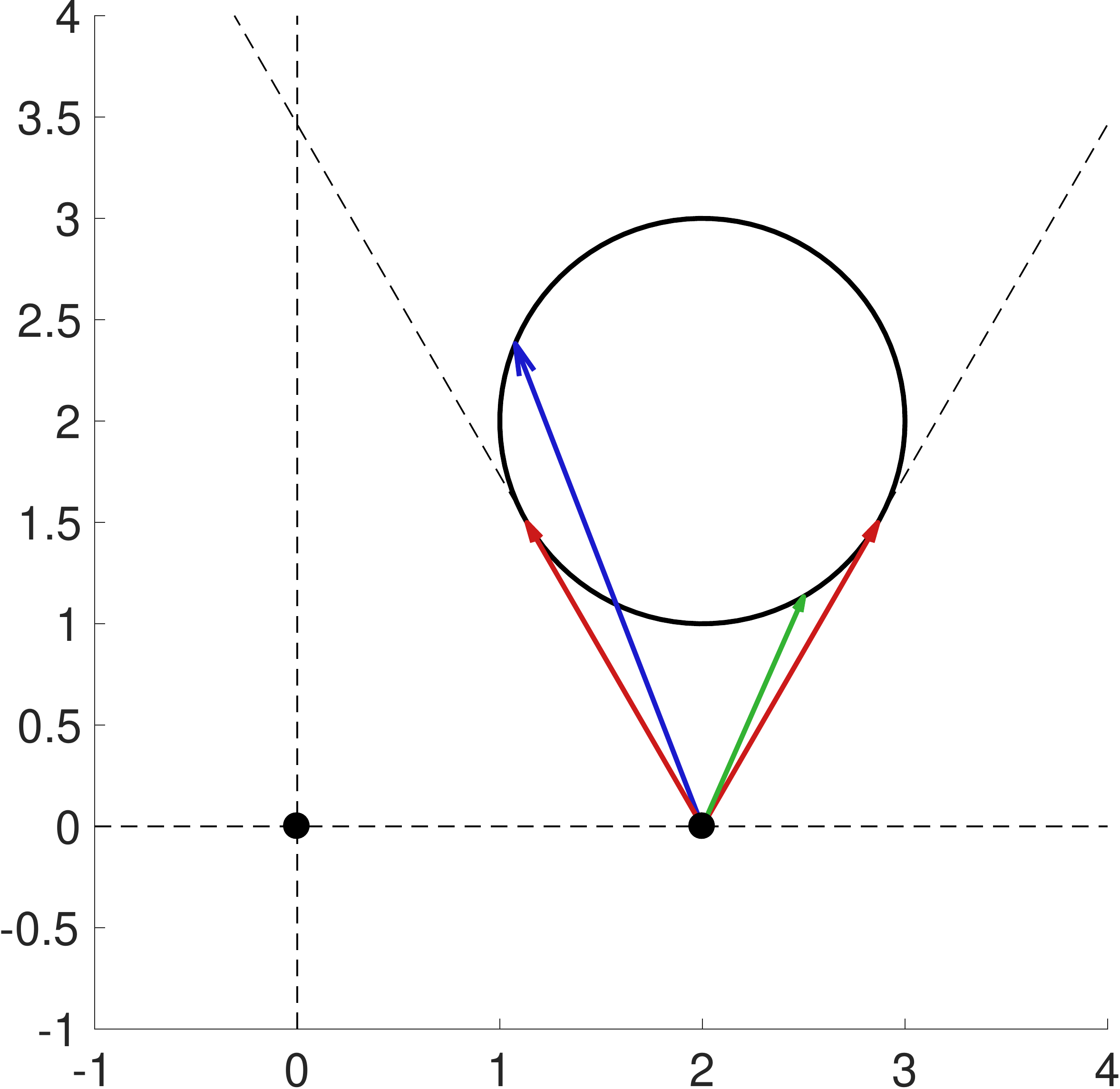}
\caption{The vortex is placed at the origin and marked by a black dot, as the initial point $x_0 \coloneqq (2,0)$
in cartesian coordinates. The drift is strong at the initial point since $\mu = 2\, r_0$.
The {black circle} represents the set of initial directions $x_0 + \Fcal(x_0,\mu)$ and the {{blue vector}}
is the initial hyperbolic direction associated to the initial control $u(0) \coloneqq (\cos \alpha, \sin \alpha)$, with $\alpha = 7\pi/8$.
The direction in {{green}} is elliptic while the two {{red}} directions
are the abnormal directions located at the boundary of the cone of admissible directions.}
\label{fig:abnormal_directions}
\end{figure}

According to the PMP, $p^0=0$ for the abnormal extremals while $p^0 < 0$ for the normal extremals. In the normal case, by homogeneity, one can
fix $p^0 = -1$ and the initial adjoint vector $p_0 \coloneqq p(0)$ of normal extremals lives in the one dimensional space $\enstq{p \in \R^2}{\Htrue(x_0, p) = 1}$.
This parameterization is very classical. Another possibility is to set $\norme{(p_0, p^0)} = 1$ since by the PMP, the pair $(p(\cdot),p^0)$
does not vanish. Finally, in our case, since all the extremals are of order zero, that is since $p$ does not vanish, we can also fix $\norme{p_0}=1$.
We consider this third possibility but in polar coordinates, that is, denoting $p_0 \coloneqq (p_{r}(0), p_\theta)$ (recalling that $p_\theta$ is constant)
we parameterize the initial adjoint vector by:
\[
p_0 \in \enstq{p \in \R^2}{\norme{p}_{r_0}=1}.
\]
We thus introduce $\alpha \in \intervallefo{0}{2\pi}$ such that
$
p_{r}(0) = \cos\alpha$, $p_\theta = r_0 \sin \alpha,
$
which gives the initial control
$
v(0) = (p_{r}(0), {p_\theta}/{r_0}) = (\cos\alpha, \sin\alpha).
$
According to the Mathieu transformation \eqref{eq:Mathieu_transform}, one has in cartesian coordinates that
$
p_x(0) = \cos\theta_0 \cos\alpha - \sin\theta_0 \sin\alpha$ and
$p_y(0) = \sin\theta_0 \cos\alpha + \cos\theta_0 \sin\alpha,
$
so in the particular case $\theta_0=0$, we have
$
u(0) = (p_x(0), p_y(0)) = (\cos\alpha, \sin\alpha).
$
This parameterization has the advantage to cover the normal and the abnormal extremals. According to the PMP, we have at the initial time:
\[
\Htrue(q_0, p_0) = p_\theta \frac{\mu}{r_0^2} + \norme{p_0}_{r_0} = p_\theta \frac{\mu}{r_0^2} + 1 = -p^0 \ge 0,
\]
with $q_0 = (r_0, \theta_0)$.
Introducing (with a slight abuse of notation)
$
\Htrue(\alpha) \coloneqq {\mu} \sin\alpha/ {r_0} + 1,
$
then, the abnormal extremals are characterized by
$
\Htrue(\alpha) = 0 \Leftrightarrow \sin\alpha = -{r_0}/{\mu}.
$
We have three cases:
\begin{itemize}
\item If the drift is weak at the initial point, then this equation has no solution which explains why there is no abnormal extremals.
In this case, $\Htrue(\alpha) > 0$ for any $\alpha$ and thus there are only hyperbolic geodesics.
\item If the drift is moderate at the initial point, that is if $\abs{\mu}=r_0$, then the single abnormal extremal is parameterized by
$\alpha = {\pi}/{2}$ if $\mu < 0$ and by $\alpha = {3\pi}/{2}$ if $\mu > 0$, 
\item In the last case when the drift is strong, then for a given $\mu$, the equation $\Htrue(\alpha) = 0$ has two distinct solutions
$\alpha_1^a < \alpha_2^a$ in $\intervallefo{0}{2\pi}$. If $\mu < 0$, then $\alpha_1^a$ and $\alpha_2^a$ are contained in $\intervalleoo{0}{\pi}$
while if $\mu > 0$, then $\alpha_1^a$ and $\alpha_2^a$ are contained in $\intervalleoo{\pi}{2\pi}$. We have in addition the following symmetry:
\[
\alpha_2^a =  \pi - \alpha_1^a \quad\text{if}\quad \mu < 0 \quad \text{and} \quad 
\alpha_2^a = 3\pi - \alpha_1^a \quad\text{if}\quad \mu > 0.
\]
The normal extremals solution of the PMP are parameterized by the set $\enstq{\alpha \in \intervallefo{0}{2\pi}}{\Htrue(\alpha) > 0} = 
\intervallefo{0}{\alpha_1^a} \cup \intervalleoo{\alpha_2^a}{2\pi}$ while for $\alpha \in \intervalleoo{\alpha_1^a}{\alpha_2^a}$
we have $\Htrue(\alpha) < 0$. 
One can see on Figure~\ref{fig:abnormal_directions}, the two abnormal directions with two hyperbolic and elliptic
directions (that is resp. associated to hyperbolic and elliptic geodesics).
The two abnormal directions define the boundary of the cone of admissible directions and reveal a lack of accessibility
in the neighborhood of $x_0$.
\end{itemize}

\subsection{Numerical results}

\subsubsection{Resolution of the shooting equation}



We introduce the \emph{shooting mapping}
\begin{equation}
S(T,p_0) \coloneqq \exp_{x_0}(T, p_0) - x_f,
\label{eq:shooting_function}
\end{equation}
where $x_f$ is the target and $\exp$ is the exponential mapping defined by \eqref{eq:exponential_mapping}.
The shooting mapping is defined for 
\[
(T,p_0) \in \enstq{(T,p_0) \in \Rp \times \Sphere^1}{ T < t_{p_0}}
\]
where $t_{p_0} \in \Rsp \cup \{+\infty\}$ is the maximal time such that $\exp_{x_0}(\cdot, p_0)$ is well defined 
over $\intervallefo{0}{t_{p_0}}$.
Let $(T, p_0)$ be a solution of $S=0$ such that the associated extremal is normal.
The shooting mapping is differentiable at $(T, p_0)$ and if $T$ is not a conjugate time, then
its Jacobian is of full rank at $(T, p_0)$,
which is a necessary condition to compute numerically the BC-extremals by means of Newton-like algorithms.
We present in the following some examples of hyperbolic geodesics fixing the initial condition to $x_0 \coloneqq (2,0)$
and solving the shooting equations $S=0$ thanks to the \hampath\ code \cite{CCG:2010}, for different final conditions and for 
different strengths of the drift.

\medskip
\paragraph{\textbf{\hampath\ code}}
A Newton-like algorithm is used to solve the shooting equation $S(T,p_0)=0$. Providing $\Htrue$ and $S$ to \hampath, the code generates
automatically the Jacobian of the shooting function. To make the implementation of $S$ easier, \hampath\ supplies the exponential mapping.
Automatic Differentiation is used to produce $\vv{\Htrue}$ and is combined with Runge-Kutta integrators to assemble the exponential mapping
and the variational equations \eqref{eq:Jacobi} used to compute conjugate times.
See \cite{CCG:2010, Co2012} for more details about the code.

\medskip
\paragraph{\textbf{Example 1}}
For this first example we want to steer the particle from $x_0$ to $x_f \coloneqq (-2,0)$ with $\mu\coloneqq 2\, \norme{x_0}$ (strong drift).
In this case, we obtain a final time $T \approx 1.641$ and the shooting equation $S=0$, is solved with a very good accuracy of order $1e^{-12}$
(it is the same for the others examples but it won't be mentioned anymore).
The associated hyperbolic geodesic is portrayed on Figure~\ref{fig_ex_1_2_strong_drift}.
The point vortex is represented by a black dot as the initial condition. The initial velocity $\dot{x}(0)$ is given with
the boundary (the black circle) of $x_0 + \mathcal{F}(x_0, \mu)$, cf. eq.~\eqref{eq:drift_domain}. One can see that the drift is strong since
$x_0 \not\in x_0 + \mathcal{F}(x_0, \mu)$.

\paragraph{\textbf{Example 2}}
To emphazise the influence of the final condition, let us take again $\mu\coloneqq 2\,\norme{x_0}$ and set $x_f\coloneqq (2.5,0)$.
We can note from Figure~\ref{fig_ex_1_2_strong_drift} that the solution turns around the point vortex and profits from the circulation.
In this case we obtain a final time $T \approx 2.821$.

\def\sizefig{0.4}
\begin{figure}[ht!]
\centering
\hspace{-2em}\includegraphics[width=\sizefig\textwidth]{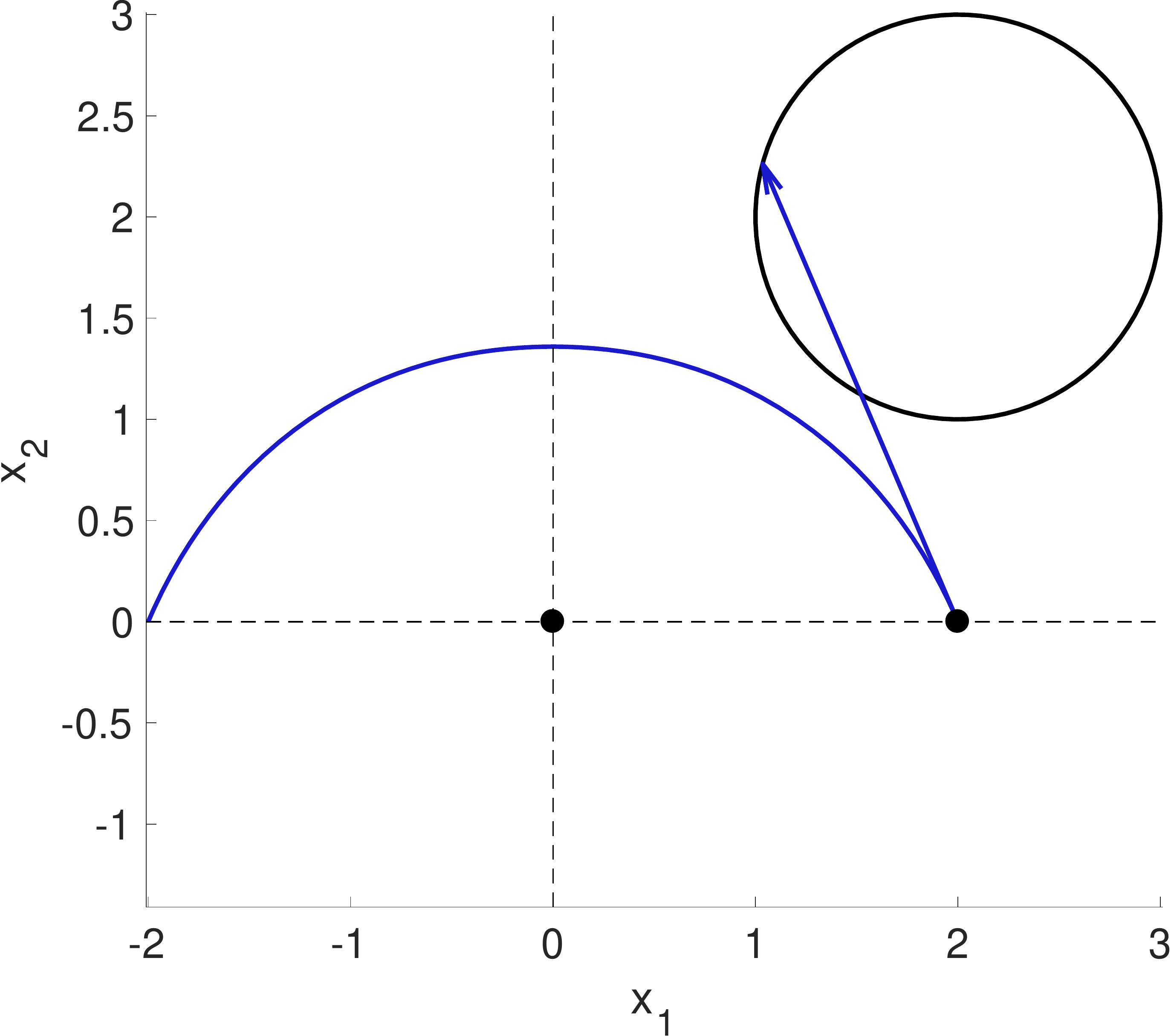}
\hspace{2em}\includegraphics[width=\sizefig\textwidth]{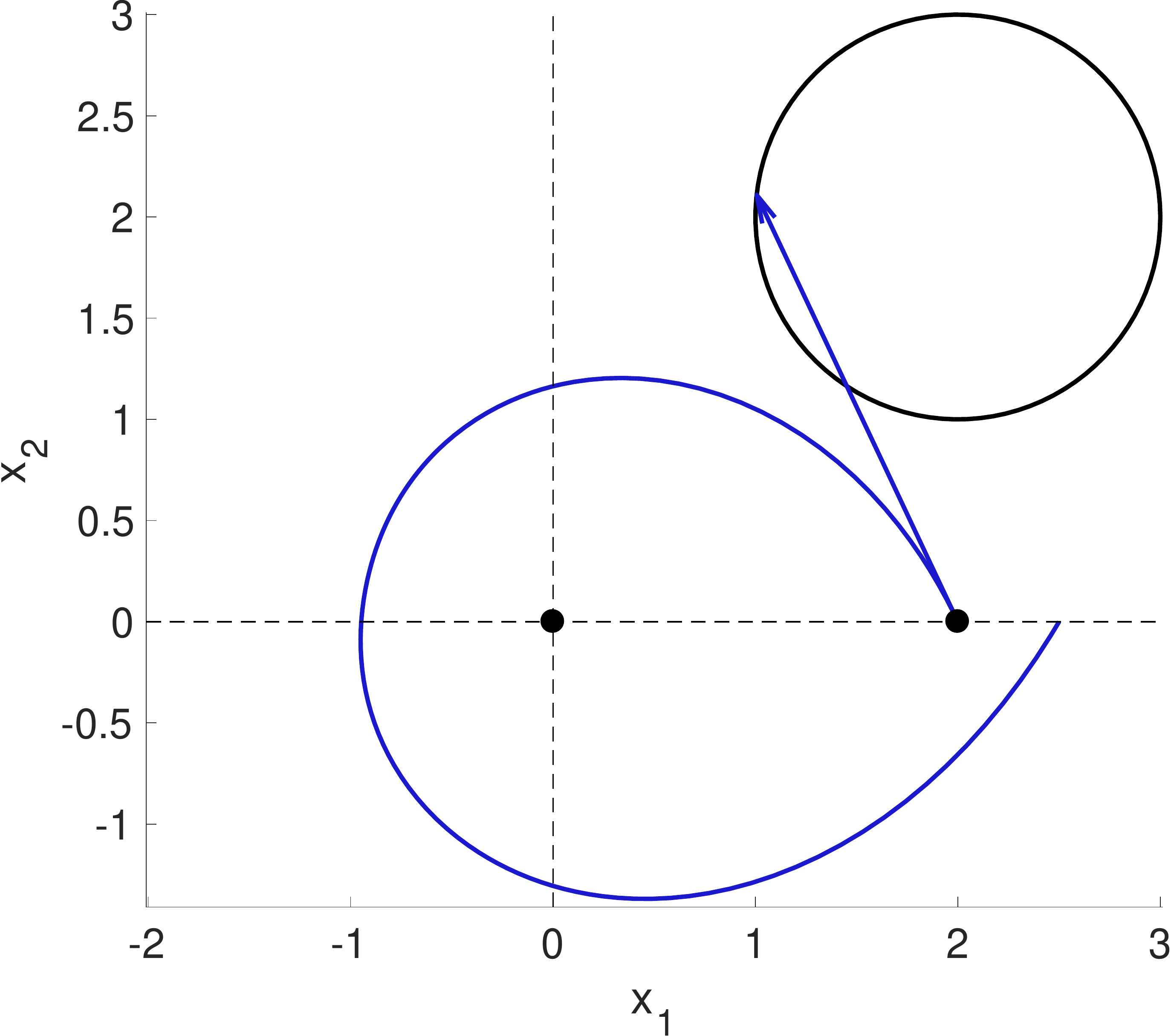}
\vspace{-0.5em}
\caption{\textbf{Example 1 and 2.} Geodesic with a strong drift at the initial point: $\mu=2\,\norme{x_0}$.
(Left: example~1) $x_0=(2,0)$, $x_f=(-2,0)$ and the final time is $T\approx 1.641$.
(Right: example~2) $x_0=(2,0)$, $x_f=(2.5,0)$ and the final time is $T\approx 2.821$.}
\label{fig_ex_1_2_strong_drift}
\end{figure}

\begin{figure}[ht!]
\centering
\hspace{-2em}\includegraphics[width=\sizefig\textwidth]{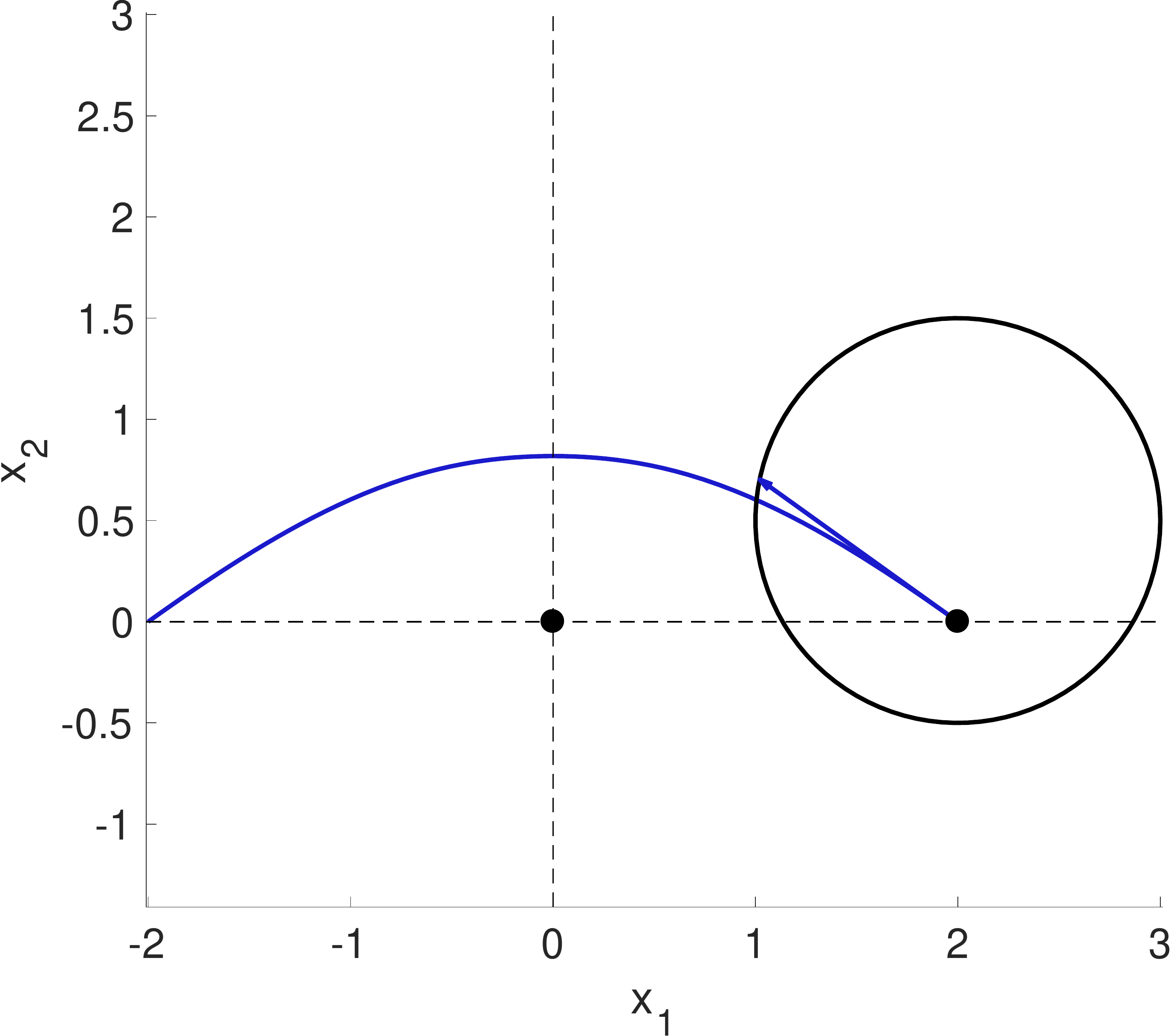}
\hspace{2em}\includegraphics[width=\sizefig\textwidth]{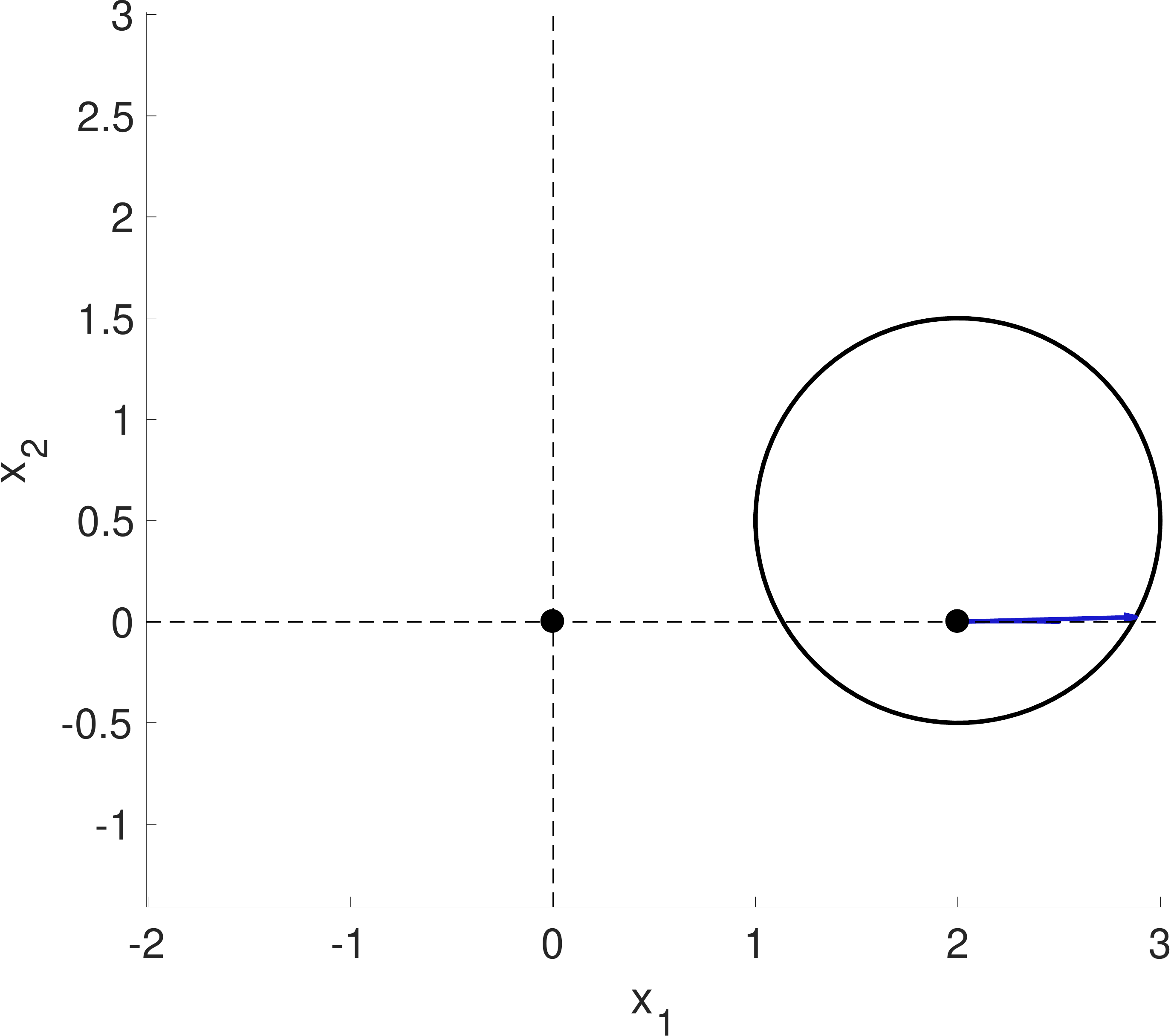}
\vspace{-0.5em}
\caption{\textbf{Example 3 and 4.} Geodesic with a weak drift at the initial point: $\mu=0.5\,\norme{x_0}$.
(Left: example~3) $x_0=(2,0)$, $x_f=(-2,0)$ and the final time is $T\approx 2.826$.
(Right: example~4) $x_0=(2,0)$, $x_f=(2.5,0)$ and the final time is $T\approx 0.56$.
Compare this to Figure~\ref{fig_ex_1_2_strong_drift}.}
\label{fig_ex_3_4_weak_drift}
\end{figure}

\paragraph{\textbf{Examples 3-4}}
Here we want to observe what happens for a weak drift. We set $\mu\coloneqq 0.5\norme{x_0}$ and present two cases with
$x_f\coloneqq (-2,0)$ (cf. left subgraph of Figure~\ref{fig_ex_3_4_weak_drift}) and $x_f\coloneqq (2.5,0)$ 
(cf. right subgraph of Figure~\ref{fig_ex_3_4_weak_drift}).
When $x_f = (-2,0)$, the final condition is the same as in the example 1 but since the drift is weaker, the final time is longer. 
This is because the particle takes advantage of the vortex circulation.
On the other hand, for $x_f=(2.5,0)$ (same final condition as example 2) and considering a weak drift, then the particle does not turn around the vortex.

\subsubsection{Numerical results on the absence of conjugate points}
\label{sec:NumercialResultsConj}

\begin{figure}[ht!]
\centering
\includegraphics[width=0.45\textwidth]{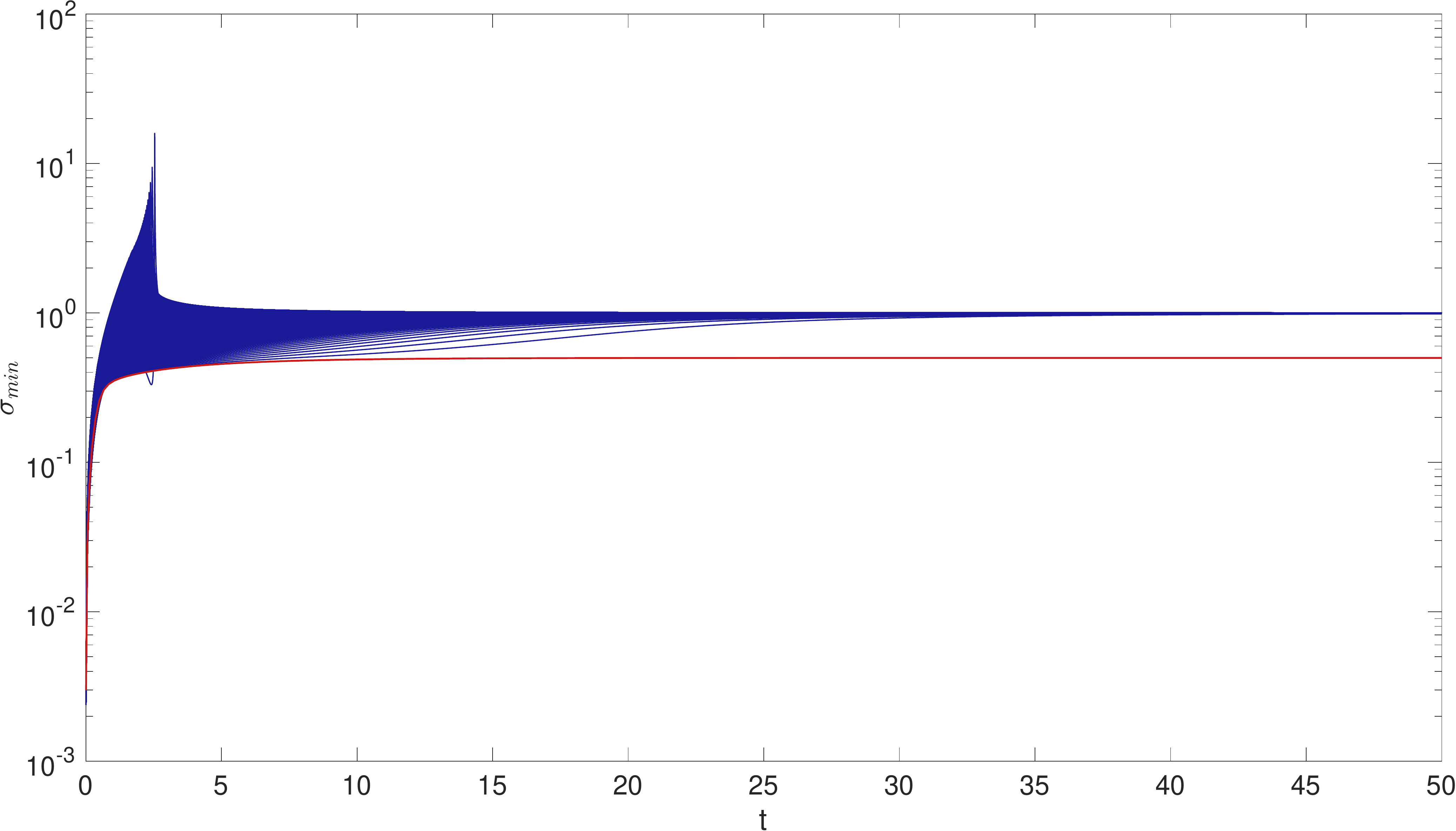}
\includegraphics[width=0.52\textwidth]{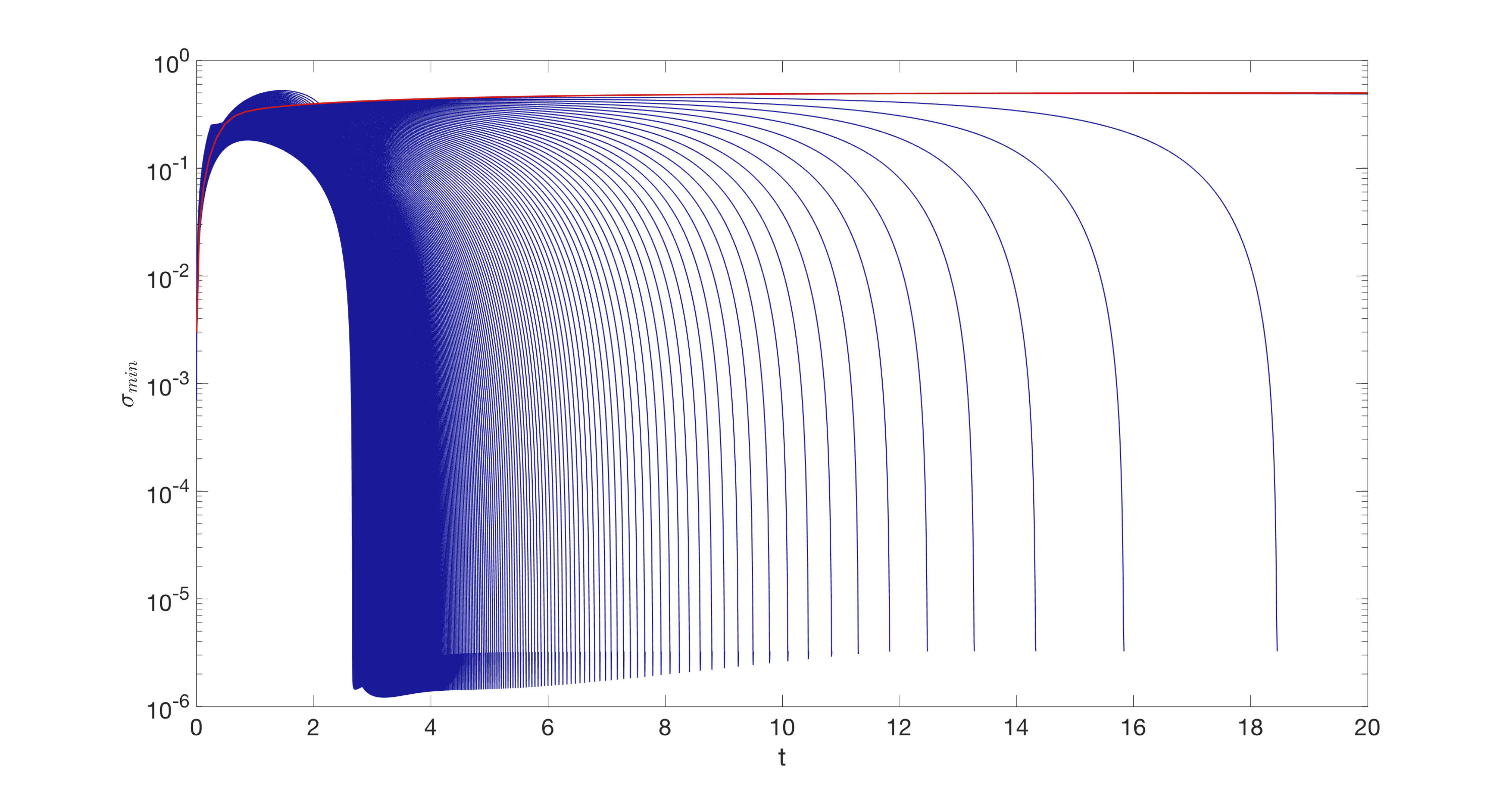}
\caption{Smallest singular value with respect to time. Setting: $\theta_0 = 0$, $\mu = 2.0$, $r_0 = 4\mu/3$ (weak drift at the initial point).
The set of parameters $\intervallefo{0}{2\pi}$ is uniformly discretized in $N\coloneqq 1000$ sub-intervals and we define 
$0\eqqcolon \alpha_1 < \cdots < \alpha_N < 2\pi$ the associated parameters.
Each blue curve is the graph of the smallest singular value
$\sigma_{\mathrm{min}}(t)$ for one $\alpha \in \Theta(r_0,\mu) \cap \cup_{i=1}^N \{\alpha_i\}$ and for 
$t \in \intervalleff{0}{50}$ on the left subgraph, and for one $\alpha \in \Lambda(r_0,\mu) \cap \cup_{i=1}^N \{\alpha_i\}$ on the right subgraph. In this case, we stop the numerical integration 
when $r(t) \le 10^{-3}$ which explains why the minimal value of $\sigma_{\mathrm{min}}$ is between $10^{-5}$ and $10^{-6}$.
The red curve on each plot corresponds to the separatrix.}
\label{fig:conjugate_test}
\end{figure}

According to Section \ref{sec:micro_local}, there are two types of geodesics reaching
the boundary of the domain: the bounded geodesics that reach the vortex and the unbounded geodesics. There exists also a unique geodesic separating these
two cases, which is called the \emph{separatrix}.
Fixing $\theta_0 = 0$ and using the parameterization
$
u(0) = (\cos \alpha, \sin \alpha)$, $\alpha \in \intervallefo{0}{2\pi},
$
from Section \ref{sec:abnormal_extremals},
we define for a pair $(r_0, \mu) \in \Rsp \times \Rs$, the set $\Lambda(r_0,\mu)$ of parameters $\alpha$ 
such that the associated geodesics converge
to the vortex and $\Theta(r_0,\mu)$ the set of parameters such that the associated geodesics go to infinity (in norm). The sets $\Lambda(r_0,\mu)$ and
$\Theta(r_0,\mu)$ depend on the strength of the drift and the sign of $\mu$, and they are given in Section \ref{sec:micro_local}.
One can find in Figure~\ref{fig:conjugate_test} the smallest singular value,
denoted $\sigma_{\mathrm{min}}(\cdot)$, of
$
\det(\dot{x}(\cdot), \delta x(\cdot))
$
over the time, see eq.~\eqref{eq:conjugate_time_eq}.
For a fixed $\alpha \in \intervallefo{0}{2\pi}$, we denote by $t_\alpha \in \Rsp \cup \{+\infty\}$ the maximal time such that the associated geodesic
is well defined over $\intervallefo{0}{t_\alpha}$.
If there exists a time $t_c \in \intervalleoo{0}{t_\alpha}$ such that $\sigma_{\mathrm{min}}(t_c) = 0$, then the time $t_c$ is a conjugate time.
If not, the geodesic has no conjugate time.
One can see from the left subgraph of Figure~\ref{fig:conjugate_test}, that for any weak drift and for any $\alpha \in \Theta(r_0,\mu)$, that
there are no conjugate times and $\sigma_{\mathrm{min}}(t) \to 1$ when $t \to t_\alpha = +\infty$.\footnote{It is clear that $t_\alpha = +\infty$ 
for any $\alpha \in \Theta(r_0,\mu)$, since $\abs{\dot{r}}\le 1$.}
From the right subgraph of Figure~\ref{fig:conjugate_test}, it is clear that for any weak drift and for any 
$\alpha \in \Lambda(r_0,\mu)$, that there are no conjugate times and that $\sigma_{\mathrm{min}}(t) \to 0$ when $t \to t_\alpha$.
On the two figures, the red curve corresponds to the separatrix. It is clear that for the separatrix 
$\sigma_{\mathrm{min}}(t) \to 0.5$ when $t \to t_\alpha = +\infty$.
Since one has similar numerical results in the strong drift case, we make the following conjecture:
\begin{conjecture}
\label{conj:no_conj}
For any $(x_0, \mu) \in M \times \Rs$, the conjugate locus is empty.
\end{conjecture}

\section{Micro-local analysis and properties of the value function}
\label{sec:microlocal_properties_value_function}

\subsection{Poincar\'e compactification on $\Sphere^3$ of the extremal dynamics and integrability results}
\label{sec:compactification}

The objective of this section is to provide the geometric frame to analyze the extremals of order zero. Reparameterizing, the flow defines a polynomial vector field wich can be compactified using Poincar\'e method to analyze the behaviors of extremal curves converging either to the origin or to the infinity, see Section \ref{sec:classification}. Thanks to the rotational symmetry, this foliation can be integrated which is crucial to define in the next section the concept of Reeb circle.



Introducing $c \coloneqq -p^0$, then from the condition $\Htrue = p_\theta {\mu}/{r^2} + \norme{p}_{r} = -p^0$, one gets 
$r^2 \norme{p}_r = c r^2 - \mu p_\theta$. Plugging this into \eqref{eq:hamiltonian_system_polar}, one obtains:
\[
r^3(cr^2-\mu\, p_\theta)\, \dot{r}      = r^5 p_r, \quad
r^3(cr^2-\mu\, p_\theta)\, \dot{\theta} = \lambda_3 r^3-\lambda_4 r, \quad
r^3(cr^2-\mu\, p_\theta)\, \dot{p}_r    = \lambda_1 r^2- \lambda_2,
\]
where $\lambda_1 \coloneqq (2 \mu c+ p_\theta) p_\theta$, $\lambda_2 \coloneqq 2 (\mu p_\theta)^2$,
$\lambda_3 \coloneqq \mu c+ p_\theta$ and $\lambda_4 \coloneqq \mu^2 p_\theta$.
The aim is to get a polynomial system to integrate, so we use the time reparameterization:
\[
\frac{\diff t}{r^3(cr^2-\mu p_\theta)} = \diff s.
\]
Denoting by $'$ the derivative with respect to $s$, the system is written:
\[
r'      = r^5 p_r, \quad
\theta' = \lambda_3 r^3- \lambda_4 r, \quad
p_r'    = \lambda_1 r^2- \lambda_2.
\]
We use Poincar\'e compactification where the sphere is identified to $x=1$ 
%
%
%
%
%
%
and this leads to the following dynamics:
\begin{equation}
\label{eq:compactifier}
r'      = r^5 p_r, \quad
\theta' = \lambda_3 r^3 x^3- \lambda_4 r x^5, \quad
p_r'    = \lambda_1 r^2 x^4- \lambda_2 x^6, \quad
x'      = 0,
\end{equation}
which can be projected onto the $3$-sphere $r^2+\theta^2+p_\theta^2+x^2=1$, up to a time reparameterization.
Integrating the system~\eqref{eq:compactifier}, one has $x=c_1$ and we get: 
\[
\frac{\diff r}{\diff p_r} = \frac{r^5 p_r}{\lambda_1 r^2 x^4-\lambda_2 x^6}
\]
which can integrated by separating the variables:
\[
\frac{\diff r}{r^5}\left( \lambda_1 r^2 x^4 - \lambda_2 x^6 \right) = p_r \diff p_r.
\]
Hence, we obtain:
\[
\frac{\lambda_2 c_1^6}{4r^4} - \frac{\lambda_1 c_1^4}{2r^2} = \frac{1}{2}p_r^2 + K_1.
\]
Using quadratures from \eqref{eq:compactifier}, one can obtain $s = f(r)$, $\theta = g(r)$. 
To illustrate the integrability properties, we only give the formula in the abnormal case for the $r$ component:
\begin{prpstn}
For $\alpha \in \{\alpha_1^a, \alpha_2^a\}$ (see Section \ref{sec:abnormal_extremals} for the definition of $\alpha_1^a$ and
$\alpha_2^a$), one has:
\[
r(t) = \left(s(t)^2-\frac{2\cos\alpha }{r_0} s(t) +\frac{1}{r_0^2}\right)^{-\frac{1}{2}},
~\text{where}~
s(t) \coloneqq \frac{\sin\alpha}{r_0}\tan\left(\frac{\sin\alpha}{r_0}t+
\arctan\left(-\frac{\cos\alpha}{\sin\alpha}\right)\right)+\frac{\cos\alpha}{r_0}.
\]
\end{prpstn}

\subsection{Micro-local analysis of the extremal solutions}
\label{sec:micro_local}

Since $\theta$ is a cyclic variable, and so $p_\theta$ is a parameter, then we can focus the analysis on the subsystem
\begin{equation}
\dot{r}         = \frac{p_r}{\norme{p}_r}, \quad 
\dot{p}_r       = \frac{p_\theta}{r^3} \left( 2 \mu + \frac{p_\theta}{\norme{p}_r} \right),
\label{eq:subsystem_r_pr}
\end{equation}
to determine the different types of extremals. First of all, it is clear that the unique equilibrium point satisfying $(p_r, p_\theta) \ne (0,0)$
and $r \ge 0$ is given by
\[
(\re, \pre) \coloneqq (2\abs{\mu}, 0),
\]
and $p_\theta$ must satisfy $\sign(p_\theta) = - \sign(\mu)$.
Now, introducing
\[
a(p_\theta, \mu) \coloneqq \mu^2\, p_\theta^2
\quad
\text{and}
\quad
b(p_\theta, r_0, \mu) \coloneqq - p_\theta \left( p_\theta + 2\, \mu + 2\, p_\theta\, \frac{\mu^2}{r_0^2} \right),
\]
one can define
\[
\vphi(r, p_\theta, r_0, \mu) \coloneqq a(p_\theta, \mu) \left( \frac{1}{r^4} - \frac{1}{r_0^4} \right)
+ b(p_\theta, r_0, \mu) \left( \frac{1}{r^2} - \frac{1}{r_0^2} \right)  + 1 - \frac{p_\theta^2}{r_0^2}
\]
and one has the following result:
\begin{lmm}
Along any extremal parameterized by $\norme{(p_r(0), p_\theta)}_{r(0)} = 1$ holds
$
p_r(t)^2 = \vphi(r(t), p_\theta, r(0), \mu)
$.
\end{lmm}
\begin{proof}
From \eqref{eq:subsystem_r_pr}, one has along the extremal (omitting the time variable):
\begin{equation}
p_r \dot{p}_r = \frac{p_\theta}{r^3} \left(p_\theta + 2 \mu\, \norme{p}_{r} \right) \dot{r}. 
\label{eq:pr_dot_pr}
\end{equation}
Using the parameterization $\norme{p_0}_{r_0} = 1$ (with $p_0\coloneqq (p_r(0), p_\theta)$ and $r_0 \coloneqq r(0)$), one has
\begin{equation}
\norme{p(t)}_{r(t)} = - \left( p^0 + p_\theta \frac{\mu}{r(t)^2} \right) = p_\theta\, \mu \left( \frac{1}{r_0^2} - \frac{1}{r(t)^2} \right) + 1,
\label{eq:norme_p_de_t}
\end{equation}
since along the extremal the Hamiltonian is constant and equal to $-p^0$. Putting \eqref{eq:norme_p_de_t} in \eqref{eq:pr_dot_pr}, one has
\[
p_r \dot{p}_r = - \left(\frac{2a}{r^5} + \frac{b}{r^3}\right) \dot{r}.
\]
Integrating this equation, we have
\[
\frac{1}{2} p_r^2 =  \frac{1}{2} \left( \frac{a}{r^4} + \frac{b}{r^2} + c \right),
\quad
c \coloneqq p_r(0)^2 - \frac{a}{r_0^4} - \frac{b}{r_0^2} ,
\]
and since
$
p_r(0)^2 = 1 - {p_\theta^2}/{r_0^2},
$
we get the conclusion.
\end{proof}

Let us fix  $p_\theta$, $r_0$ and $\mu$ and introduce the polynomial function of degree two
$
P(X) \coloneqq a X^2 + b X + c,
$
with
\[
\Delta \coloneqq b^2 - 4 a  c = p_\theta^3 \left( p_\theta \left( 1 + 4 \frac{\mu^2}{r_0^2} \right) + 4 \mu \right),
\]
its discrimant. Then, we have $\vphi(r, p_\theta, r_0, \mu) = P(1/r^2)$ and the sign of the discrimant plays a crucial role in the analysis.

\begin{rmrk}
If $p_\theta = 0$ then $P(X) = 1$ and if $\mu = 0$ then $P(X) = -p_\theta^2 X + c$.
The polynomial $P$ is of degree two if and only if $p_\theta \ne 0$ and $\mu \ne 0$.
The case $\mu = 0$ is the simple Euclidean case that we do not develop.
If $p_\theta = 0$ then $\dot{p}_r = 0$, $\dot{r} = \pm 1$ and this case is clear: these are the fastest geodesics converging either to the vortex
($\dot{r}=-1$) or going to infinity ($\dot{r}=+1$).
\end{rmrk}

Let us introduce now
\[
p_\theta^* \coloneqq -\frac{4\, \mu}{1 + 4 \frac{\mu^2}{r_0^2}},
\]
such that, assuming $p_\theta \ne 0$, $\mu > 0$ and $r_0 > 0$, we have:
\[
\Delta(p_\theta, r_0, \mu) 
\left\{
\begin{aligned}
< 0 & ~\text{ if }~ p_\theta \in \intervalleoo{p_\theta^*}{0}, \\
= 0 & ~\text{ if }~ p_\theta = p_\theta^*, \\
> 0 & ~\text{ if }~ p_\theta \in \R \setminus \intervalleff{p_\theta^*}{0}. 
\end{aligned}
\right.
\]

\begin{rmrk}
If $\mu < 0$, then $p_\theta^* > 0$ and we have $\Delta(p_\theta, r_0, \mu) < 0$ if $p_\theta \in \intervalleoo{0}{p_\theta^*}$,
$\Delta(p_\theta^*, r_0, \mu) = 0$ and $\Delta(p_\theta, r_0, \mu) > 0$ if $p_\theta \in \R \setminus \intervalleff{0}{p_\theta^*}$.
\end{rmrk}

\begin{rmrk}
Note that with our parameterization, we are only interested in parameters $p_\theta$ satisfying $\abs{p_\theta} \le r_0$.
We can check that for any $r_0 > 0$ and $\mu \ne 0$, we have $\abs{p_\theta^*} \le r_0$. The particular case $\abs{p_\theta^*} = r_0$ 
is given by $r_0 = 2\abs{\mu}$.
Thus, setting $p_\theta = r_0 \sin \alpha$, we can find two parameters $0 < \alpha_1^* \le \alpha_2^* < 2\pi$ such that 
$p_\theta^* \eqqcolon r_0 \sin \alpha_1^* = r_0 \sin \alpha_2^*$. This gives us
\[
\sin\alpha_1^* = \sin\alpha_2^* = -\frac{4\,\lambda}{1+4\lambda^2}, \quad \lambda \coloneqq \frac{\mu}{r_0},
\]
and we can see again that the ratio $\mu/r_0$ is of particular interest. Note that $\alpha_1^* = \alpha_2^*$ if $r_0 = 2\abs{\mu}$.
\label{rmk:ratio_important}
\end{rmrk}

Let us analyze the three different cases depending on the sign of the discrimant, considering that $p_\theta \ne 0$
and fixing the parameters $\mu > 0$ and $r_0 > 0$.
\begin{itemize}
\medskip
\item \textbf{Case $\Delta < 0$.} In this case, the polynomial $P$ has no roots and so $p_r(t)^2 = P(1/r(t)^2)$ never vanishes.
See Figure~\ref{fig:delta_negatif} for an illustrative example.
Thus, $p_r$ is of constant sign such as $\dot{r}$ and so $r$ is monotone, even strictly monotone since $p_r$ is nonzero.
This case is given by $p_\theta \in \intervalleoo{p_\theta^*}{0}$, thus, in this case $\abs{p_\theta} \le \abs{p_\theta^*} \le r_0$
and setting $(p_r(0), p_\theta) = (\cos\alpha, r_0\sin\alpha)$, then $r(t) \to 0$ (\ie\ to the vortex)
when $t \to t_\alpha$ if $\cos\alpha<0$ and $r(t) \to +\infty$ when $t \to t_\alpha$ if $\cos\alpha>0$, 
where $t_\alpha \in \Rsp \cup \{+\infty\}$
is the maximal positive time such that the associated geodesic is defined over $\intervallefo{0}{t_\alpha}$.
Denoting
\[
\alpha^* \coloneqq \arcsin\left(\frac{p_\theta^*}{r_0}\right) \in \intervallefo{-\frac{\pi}{2}}{0}
\]
and introducing (cf. remark \ref{rmk:ratio_important})
$
\alpha_1^* \coloneqq \pi - \alpha^*$ and $\alpha_2^* = 2\pi + \alpha^* = 3\pi - \alpha_1^*,
$
then, for $\alpha \in \intervalleoo{\pi}{\alpha_1^*}$, we have $r(t) \to 0$,
while for $\alpha \in \intervalleoo{\alpha_2^*}{2\pi}$, we have $r(t) \to +\infty$, when $t \to t_\alpha$.

\begin{figure}[ht!]
\centering
\includegraphics[width=0.4\textwidth]{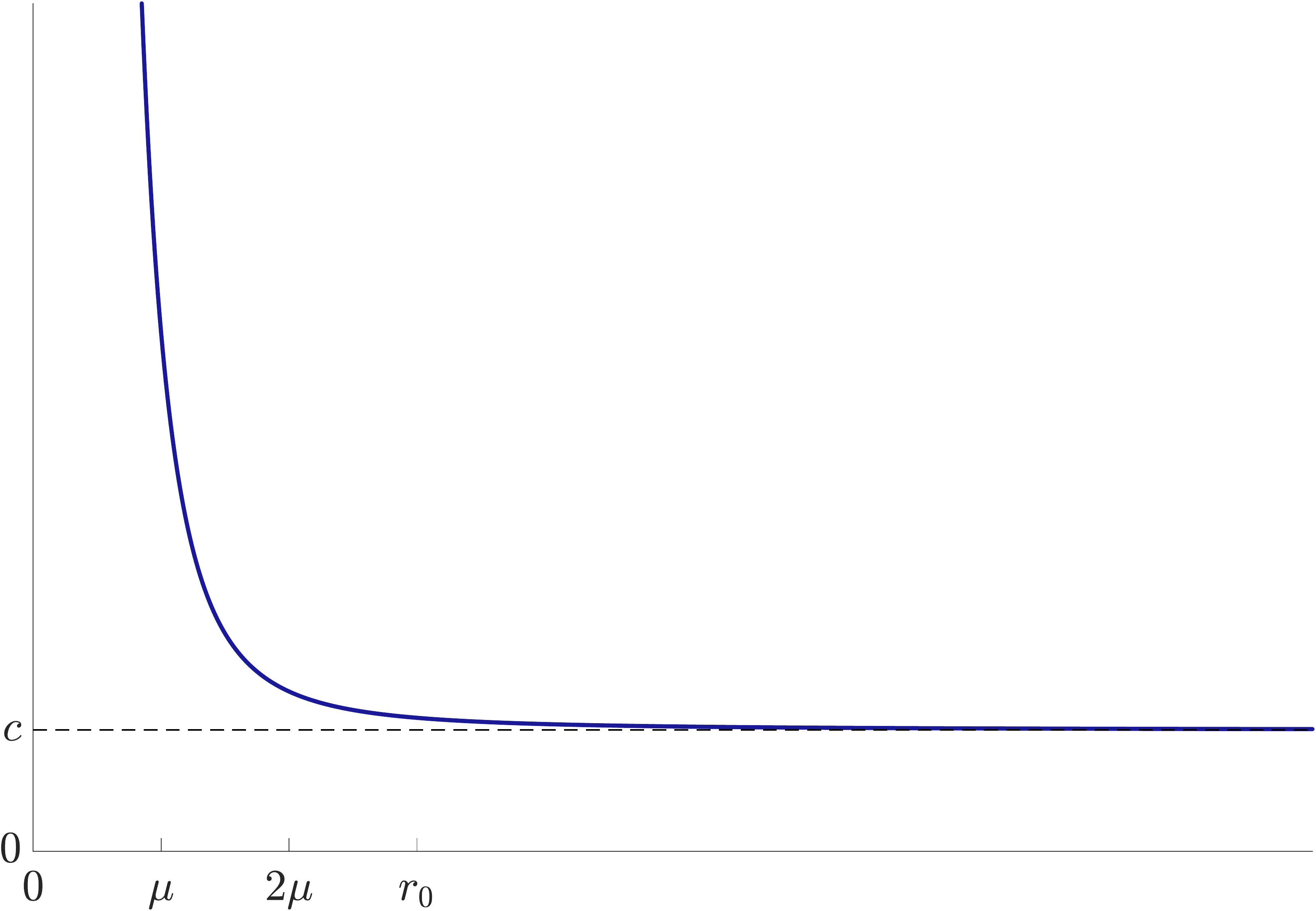}
\caption{\textbf{Case $\Delta < 0$}. Graph of $r \mapsto P(1/r^2)$.
Setting: $\mu = 2$, $r_0 = 3\mu$, $p_\theta = p_\theta^* / 2 \in \intervalleoo{p_\theta^*}{0}$.
}
\label{fig:delta_negatif}
\end{figure}

\medskip
\item \textbf{Case $\Delta = 0$.} Since we assume $p_\theta \ne 0$ then $\Delta(p_\theta, r_0, \mu) = 0$ if and only if
$p_\theta = p_\theta^*$. In this case $P(1/r^2) = 0$ if and only if $r = r^* \coloneqq 2 \mu = r_e$.
See Figure~\ref{fig:delta_zero} for an illustrative example.
We have four possibilities:
\begin{itemize}
\item If $r_0 > r^*$ and $p_r(0) > 0$, that is $p_r(0) = \cos\alpha_2^*$, then $r(t) \to +\infty$ when $t \to t_{\alpha_2^*}$.
\item If $r_0 < r^*$ and $p_r(0) < 0$, that is $p_r(0) = \cos\alpha_1^*$, then $r(t) \to 0$ when $t \to t_{\alpha_1^*}$.
\item If $r_0 > r^*$ and $p_r(0) = \cos\alpha_1^* < 0$, or if $r_0 < r^*$ and $p_r(0) = \cos\alpha_2^* > 0$, then $r(t) \to r^*$.
\item If $r_0 = r^*$ and $p_r(0) = \cos\alpha_1^* = \cos\alpha_2^* = 0$, then $r(t) = r^*$ for any time $t$ and thus, in this
case, the geodesic describes a circle.
\end{itemize}

\begin{rmrk}
We can notice that $(r_e, p_{r,e}) = (2\mu, 0)$ is an unstable hyperbolic fixed point and its
associated linearized system has $\pm 1/2\mu$ as eigenvalues with associated eigenvectors $(\pm \mu(4\mu^2+r_0^2)/r_0^2, 1)$.
\end{rmrk}

\begin{rmrk}
Let us consider the $\theta$ dynamics: 
\[
\dot{\theta} = \frac{1}{r^2} \left( \mu + \frac{p_\theta}{\norme{p}_r} \right).
\]
Fixing $p_\theta = -4\mu r^2/(r^2+4\mu^2)$ and using the parameterization $\norme{p}_r = 1$ (this amounts to consider
$r$ as the initial distance to the vortex), then we have $\dot{\theta} = 0 \Leftrightarrow \mu + p_\theta = 0 
\Leftrightarrow r = 2\abs{\mu}/\sqrt{3}$. See also eq.~\eqref{eq:separating_geo_dyn}.
\label{rmk:dot_theta_zero}
\end{rmrk}

\begin{figure}[ht!]
\centering
\includegraphics[width=0.4\textwidth]{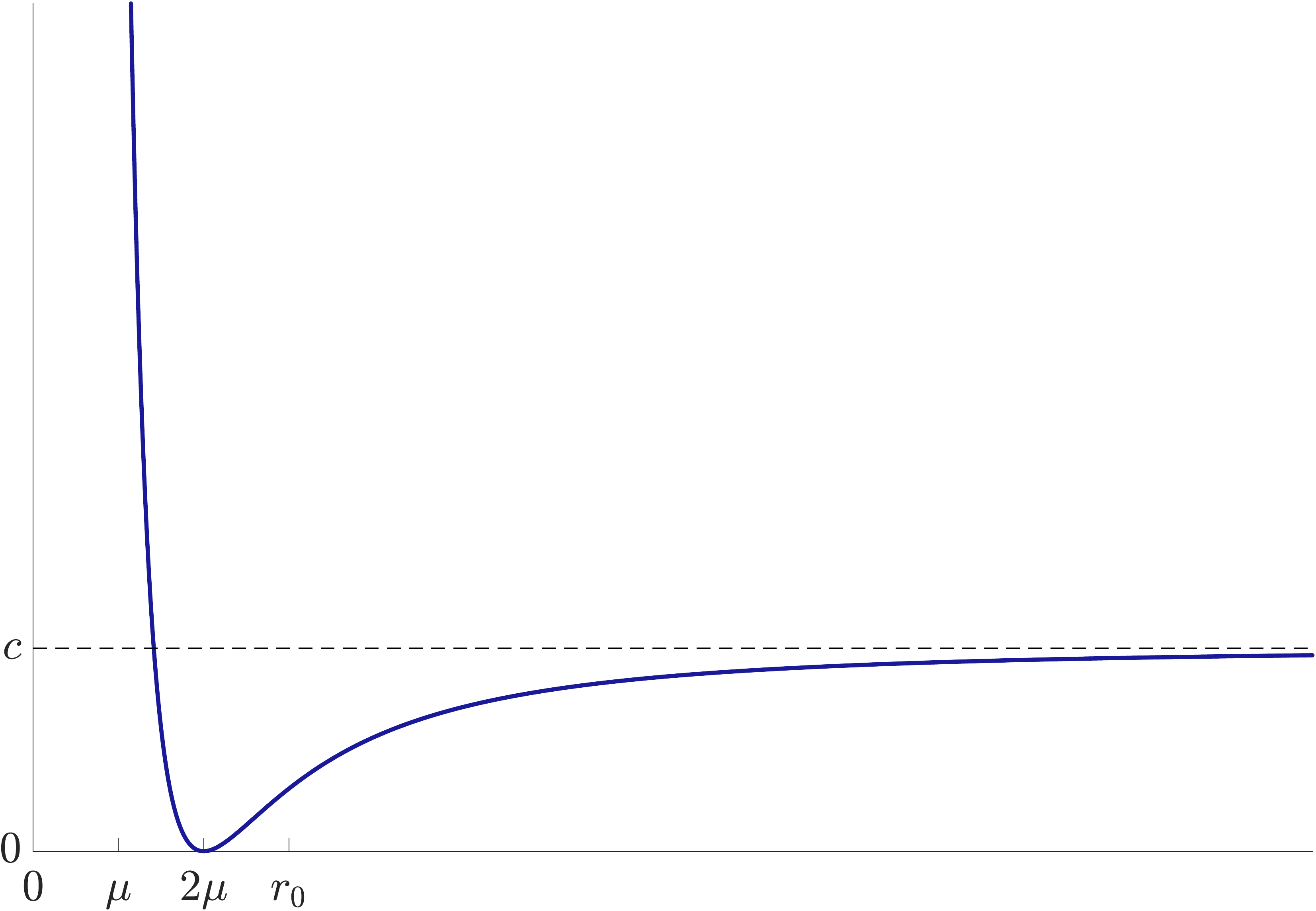}
\caption{\textbf{Case $\Delta = 0$}. Graph of $r \mapsto P(1/r^2)$.
Setting: $\mu = 2$, $r_0 = 3\mu$, $p_\theta = p_\theta^*$.
}
\label{fig:delta_zero}
\end{figure}

\medskip
\item \textbf{Case $\Delta > 0$.} In this case, $P$ has two roots $X_- < X_+$ given by:
\[
X_\pm \coloneqq \frac{p_\theta^2\, r_0^2 + 2\, \mu^2 p_\theta^2 + 2\, \mu\, p_\theta\, r_0^2 \pm
r_0 \sqrt{p_\theta^3 \left( 4\, p_\theta\, \mu^2 + 4\, \mu\, r_0^2 + p_\theta\, r_0^2 \right)}
}{2\, \mu^2\, p_\theta^2\, r_0^2}.
\]
By a tedious calculation, one can prove that $X_- \ge 0$
and $X_- = 0$ if and only if $p_\theta = - r_0^2 / \mu < p_\theta^* < 0$. Let us introduce 
\[
r_1 \coloneqq \frac{1}{\sqrt{X_+}} \quad \text{and} \quad r_2 \coloneqq \frac{1}{\sqrt{X_-}},
\]
with the convention $1/0 = +\infty$. Then, for any $r \in \intervalleoo{r_1}{r_2}$, we have $P(1/r^2) < 0$ and thus it clear that $r_0 \not\in 
\intervalleoo{r_1}{r_2}$, since by definition $P(1/r_0^2) = p_r(0)^2 \ge 0$.

\medskip
\begin{enumerate}
\item \textbf{Subcase $p_\theta > 0$.} Since $p_\theta$ is positive we have $0 < X_- < X_+$ and thus $0 < r_1 < r_2 < +\infty$.
We can easily prove that $p_\theta > 0 \Rightarrow r_1 < r_0$ since:
\[
        r_1 < r_0   \Leftrightarrow 1 < r_0^2 \, X_+ 
                    \Leftrightarrow 0 < p_\theta^2\, r_0^2 + 2\, \mu\, p_\theta\, r_0^2 + 
                    r_0 \sqrt{p_\theta^3 \left( 4\, p_\theta\, \mu^2 + 4\, \mu\, r_0^2 + p_\theta\, r_0^2 \right)}.
\]
This can be seen on Figure~\ref{fig:delta_positif_1}.
Since $r_0 \not\in \intervalleoo{r_1}{r_2}$, then $p_\theta > 0 \Rightarrow r_2 \le r_0$ and thus $r(t)$ will go to $+\infty$ in any case.
More precisely, setting $(p_r(0), p_\theta) = (\cos\alpha, r_0\sin\alpha)$, then we have:

\medskip
\begin{itemize}
    \item If $\alpha \in \intervalleoo{0}{{\pi}/{2}}$, then $\dot{r}(0) > 0$ and thus $r_2 < r_0$. In this case we have $\dot{r}(t) > 0$
        for any $t \ge 0$ and thus $r$ is strictly increasing and $r(t) \to +\infty$ when $t \to t_\alpha$.

    \medskip
    \item If $\alpha = {\pi}/{2}$, then $\dot{r}(0) = 0$, that is $r_2 = r_0$, but $r(t)$ is stricty increasing for $t > 0$ and
        we still have $r(t) \to +\infty$ when $t \to t_\alpha$.

    \medskip
    \item If $\alpha \in \intervalleoo{{\pi}/{2}}{\pi}$, then $\dot{r}(0) < 0$ and again $r_2 < r_0$. In this case, $r(t)$ is decreasing
        over $\intervalleff{0}{\tb_\alpha}$, where $\tb_\alpha > 0$ is the time such that $P(1/r(\tb_\alpha)^2) = 0$, that is such
        that $r(\tb_\alpha) = r_2$. Since $(r_2, 0)$ is not an equilibrium point (a necessary condition is $\sign(p_\theta) = -\sign(\mu)$), and since along the geodesic holds $p_r(t)^2 = P(1/r(t)^2)$, then necessarily, $r(t)$
        is increasing over $\intervallefo{\tb_\alpha}{t_\alpha}$ (the sign of $p_r$ changes at $t = \tb_\alpha$ and no more changing may occur after),
        even strictly increasing for $t > \tb_\alpha$ and we still have $r(t) \to +\infty$ when $t \to t_\alpha$.
\end{itemize}

\begin{figure}[ht!]
    \centering
    \includegraphics[width=0.4\textwidth]{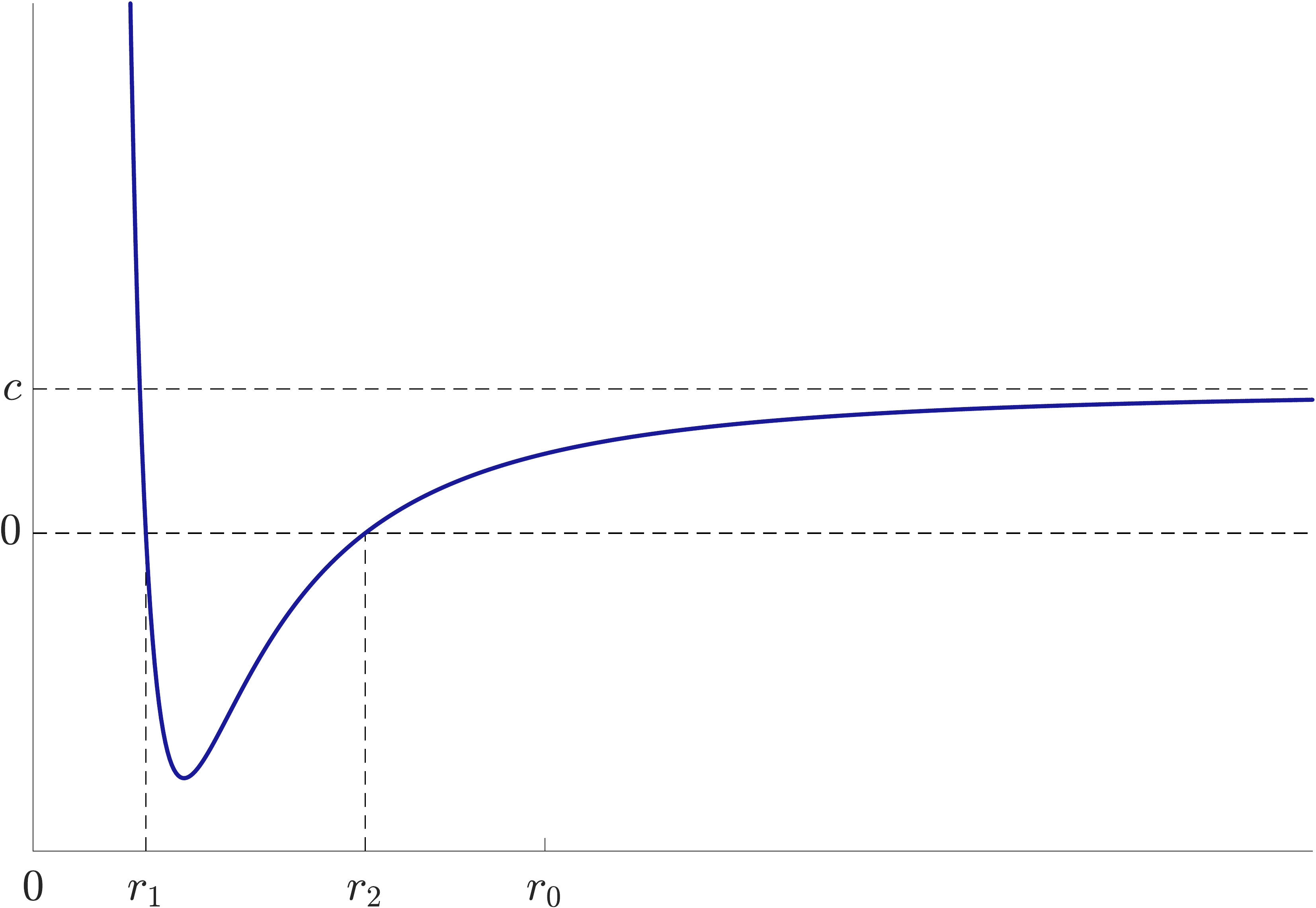}
    \caption{\textbf{Case $\Delta > 0$ and $p_\theta > 0$}. Graph of $r \mapsto P(1/r^2)$.
    Setting: $\mu = 2$, $r_0 = 3\mu$, $p_\theta = r_0/2$.
    }
    \label{fig:delta_positif_1}
\end{figure}

\medskip
\item \textbf{Subcase $p_\theta < p_\theta^*$.} Here we have $-r_0 \le p_\theta < p_\theta^* < 0$.
Setting $(p_r(0), p_\theta) = (\cos\alpha, r_0\sin\alpha)$, this case corresponds to $\alpha \in \intervalleoo{\alpha_1^*}{\alpha_2^*}$.
Let us recall that if $r_0 = r^* = 2\mu$, then ${p_\theta^*} = -r_0$ (\ie\ $\alpha_1^* = \alpha_2^* = 3\pi/2$),
and so if $r_0 = r^*$, then this last case is empty. We thus have either $r_0 < r^*$ or $r_0 > r^*$.
 Let us prove now that $r_1 < r^*$. To do this, we need the first following result:
\[
    p_\theta < p_\theta^* \Leftrightarrow p_\theta\, r_0^2 + 4\, \mu^2\, p_\theta + 4\, \mu\, r_0^2 < 0.
\]
Besides,
\[
    r_1 < r^* = 2\, \mu \Leftrightarrow 1 < 4\, \mu^2 \, X_+
    \Leftrightarrow  0 < p_\theta \left( p_\theta\, r_0^2 + 4\, \mu^2\, p_\theta + 4\, \mu\, r_0^2 \right) + 2 \sqrt{\Delta},
\]
which is true since $p_\theta < 0$ and $p_\theta < p_\theta^*$. We are now in position to conclude:

\medskip
\begin{itemize}
    \item If $r_0 > r^*$ then necessarily $r_0 \ge r_2$ (see Figure~\ref{fig:delta_positif_2a})
        since $r_0 \not\in \intervalleoo{r_1}{r_2}$ and thus $r(t) \to +\infty$ when $t \to t_\alpha$
        for any $\alpha \in \intervalleoo{\alpha_1^*}{\alpha_2^*}$. More precisely, $r(t)$ is strictly increasing if 
        $\alpha \in \intervalleoo{3\pi/2}{\alpha_2^*}$, it is increasing if $\alpha = 3\pi/2$ ($\dot{r}(0)=0$) and has one oscillation if 
        $\alpha \in \intervalleoo{\alpha_1^*}{3\pi/2}$.

    \begin{figure}[ht!]
        \centering
        \includegraphics[width=0.4\textwidth]{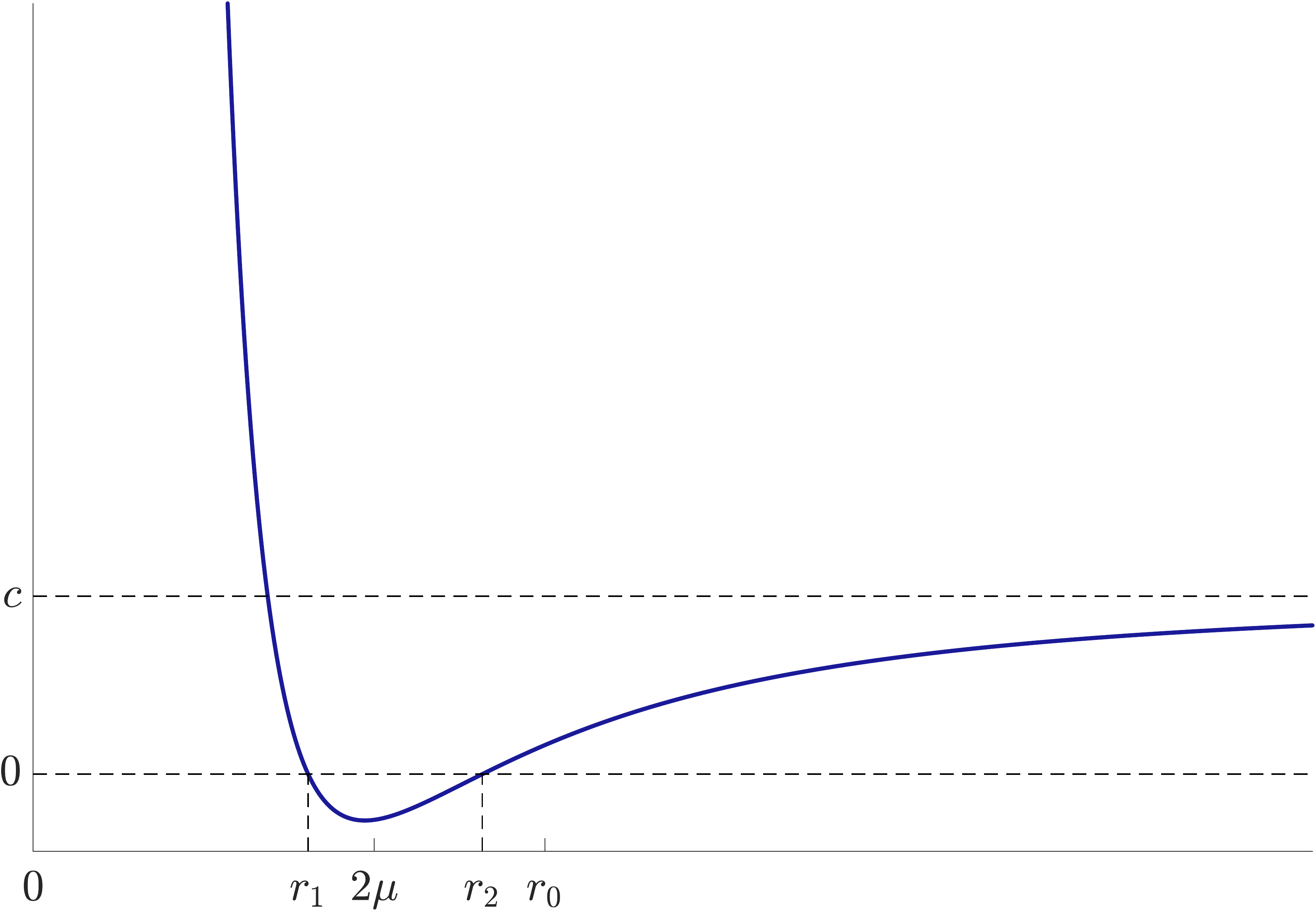}
        \caption{\textbf{Case $\Delta > 0$ and $p_\theta < p_\theta^*$}: $r \mapsto P(1/r^2)$.
        Setting: $\mu = 2$, $r_0 = 3\mu$, $p_\theta = (p_\theta^*-r_0)/2$.
        }
        \label{fig:delta_positif_2a}
    \end{figure}

    \medskip
\item If $r_0 < r^*$, using the fact that $\abs{p_\theta} \le r_0 < r^* = 2\mu$, then one can prove by a tedious calculation
    that $r_0 \le r_1$  (see Figure~\ref{fig:delta_positif_2b}).
    Then, we have $r(t) \to 0$ when $t \to t_\alpha$ for any $\alpha \in \intervalleoo{\alpha_1^*}{\alpha_2^*}$.
    More precisely, $r(t)$ is strictly decreasing if $\alpha \in \intervalleoo{\alpha_1^*}{3\pi/2}$, it is decreasing if 
    $\alpha = 3\pi/2$ ($\dot{r}(0)=0$) and has one oscillation if $\alpha \in \intervalleoo{3\pi/2}{\alpha_2^*}$.

    \begin{figure}[ht!]
        \centering
        \includegraphics[width=0.4\textwidth]{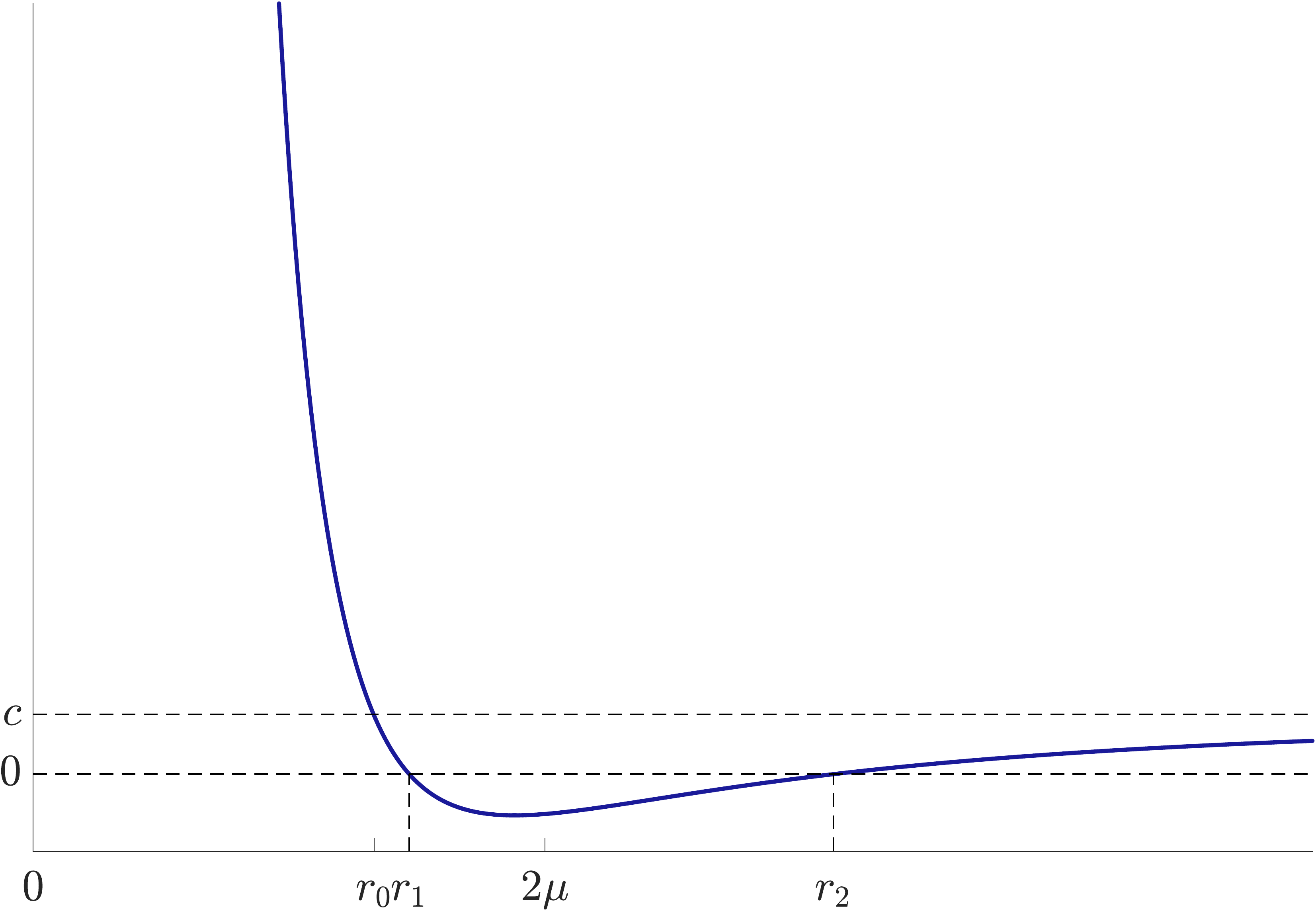}
        \caption{\textbf{Case $\Delta > 0$ and $p_\theta < p_\theta^*$}: $r \mapsto P(1/r^2)$.
        Setting: $\mu = 2$, $r_0 = 4\mu/3$, $p_\theta = (p_\theta^*-r_0)/2$.
        }
        \label{fig:delta_positif_2b}
    \end{figure}

\end{itemize}

\begin{rmrk}
    According to Section \ref{sec:abnormal_extremals}, the abnormal extremals are given by $p_\theta^a \coloneqq -r_0^2/\mu$,
    that is by $\alpha = \alpha_1^a$ and $\alpha = \alpha_2^a$ (defined in Section \ref{sec:abnormal_extremals}). Since
    (recalling that $\mu > 0$)
    \[
        p_\theta^a - p_\theta^* = -\frac{r_0^4}{\mu(r_0^2 + 4 \mu^2)} < 0,
    \]
    it is clear that $p_\theta^a < p_\theta^*$, which gives $\alpha_1^* < \alpha_1^a \le \alpha_2^a < \alpha_2^*$,
    and thus the abnormal case is contained in this last case of $\Delta > 0$. Besides, the abnormal
    extremals exist only if the drift is strong or moderate, that is if $r_0 \le \mu < r^*$. To end this remark, let us mention that
    if the drift is strong, then for any $\alpha \in \intervalleoo{\alpha_1^a}{\alpha_2^a} \ne \emptyset$ the associated extremal is elliptic 
    and from the previous analysis, one can say that all the elliptic 
    geodesics converge to the vortex, that is $r(t) \to 0$ when
    $t \to t_\alpha$.
\end{rmrk}

\end{enumerate}

We want now to group together the geodesics converging to the vortex and separate them from the ones going to infinity. We thus introduce
for $r_0 > 0$ and $\mu \ne 0$, the following sets contained in $\intervallefo{0}{2\pi}$:
\[
\begin{aligned}
\Lambda(r_0, \mu) &\coloneqq
\left\{
\begin{aligned}
\intervallefo{\pi}{\alpha_1^*} & \text{ if $r_0 \ge 2\mu$ and $\mu > 0$} \\
\intervallefo{\pi}{\alpha_2^*} & \text{ if $r_0 < 2\mu$ and $\mu > 0$} \\
\intervalleof{\alpha_2^*}{\pi} & \text{ if $r_0 \ge 2\abs{\mu}$ and $\mu < 0$} \\
\intervalleof{\alpha_1^*}{\pi} & \text{ if $r_0 < 2\abs{\mu}$ and $\mu < 0$}
\end{aligned}
\right. \\[0.5em]
\Theta(r_0, \mu) &\coloneqq
\left\{
\begin{aligned}
\intervallefo{0}{\pi} \cup \intervalleoo{\alpha_1^*}{2\pi} & \text{ if $r_0 > 2\mu$ and $\mu > 0$} \\
\intervallefo{0}{\pi} \cup \intervalleoo{\alpha_2^*}{2\pi} & \text{ if $r_0 \le 2\mu$ and $\mu > 0$} \\
\intervallefo{0}{\alpha_2^*} \cup \intervalleoo{\pi}{2\pi}  & \text{ if $r_0 > 2\abs{\mu}$ and $\mu < 0$} \\
\intervallefo{0}{\alpha_1^*} \cup \intervalleoo{\pi}{2\pi} & \text{ if $r_0 \le 2\abs{\mu}$ and $\mu < 0$}
\end{aligned}
\right. \\[0.5em]
\Psi(r_0, \mu) &\coloneqq
\left\{
\begin{aligned}
\{\alpha_1^*\} & \text{ if ($r_0 > 2\mu$ and $\mu > 0$) or if ($r_0 \le 2\abs{\mu}$ and $\mu < 0$)} \\
\{\alpha_2^*\} & \text{ if ($r_0 \le 2\mu$ and $\mu > 0$) or if ($r_0 > 2\abs{\mu}$ and $\mu < 0$)}
\end{aligned}
\right.
\end{aligned}
\]

\begin{dfntn}
The geodesics parameterized by $\alpha \in \Psi(r_0, \mu)$ are called \emph{separating geodesics}.
\end{dfntn}

We can now summarize the analysis:
\begin{thrm}
Let us fix $x_0 \in \xdom$ and $\mu \ne 0$. Then, to any $\alpha \in \intervallefo{0}{2\pi}$ is associated a unique
zero order extremal $z(\cdot) \coloneqq (x(\cdot), p(\cdot))$ solution of the Hamiltonian system $\dot{z} = \vv{\Htrue}(z)$
given by  eq.~\eqref{eq:True_Hamiltonian_System}, satisfying $x(0) = x_0$,
and parameterized in polar coordinates by $\norme{p(0)}_{r_0} = 1$ (with $r_0 \coloneqq \norme{x_0}$).
Besides, we have:
\[
\begin{aligned}
    \forall\, \alpha \in \Lambda(r_0, \mu)  &~:~ \lim_{t \to t_\alpha} r(t) = 0, \\ 
    \forall\, \alpha \in \Theta(r_0, \mu)   &~:~ \lim_{t \to t_\alpha} r(t) = +\infty \quad \text{with} \quad t_\alpha = +\infty,\\
    \forall\, \alpha \in \Psi(r_0, \mu)        &~:~ \lim_{t \to t_\alpha} r(t) = r^* = 2\abs{\mu} \quad \text{with} \quad t_\alpha = +\infty. \\
\end{aligned}
\]
\label{thm:types_extremals}
\end{thrm}

\end{itemize}

\subsection{Reeb foliations}

\subsubsection{$3$-dimensional compact Reeb component}

A concrete example of foliation
of codimension 1 on the $3$-sphere $\Sphere^3$ is given by the Reeb foliation \cite{Gobdillon:1991}.
To construct such a foliation, the idea is to glue together two $3$-dimensional foliations on two solid tori along their boundary.
Let us recall only the construction of such a foliation on the solid torus $\D^2 \times \Sphere^1$, where 
$\D^2 \coloneqq \enstq{(x,y) \in \R^2}{ x^2+y^2 \le 1} = \BallClosed(0, 1)$, following the presentation of \cite{Hector:1994}.
%
Let $f\colon \D^2 \to \R$ be a smooth mapping such that the following conditions are satisfied:
\begin{enumerate}
\item $f$ is of the form: $f(x,y) \eqqcolon \vphi(x^2+y^2)$, that is in polar coordinates, $f(r, \theta) = \vphi(r^2)$,
\item $f(\partial \D^2) = \{0\}$ and $f(x,y) > 0$ for $(x,y) \not\in \partial \D^2 = \Sphere^2$,
\item $f$ has no critical points on $\partial \D^2$.
\end{enumerate}
Let consider the smooth mapping $F\colon \D^2 \times \R \to \R$ defined by $F(x,y, t) \coloneqq f(x,y)\, e^t$, where $(x,y) \in \D^2$
and $t\in \R$. $F$ is a submersion from $\D^2 \times \R$ to $\R$ and defines a foliation of $\D^2 \times \R$ whose leaves 
are the level sets of $F$: $F^{-1}(a)$, $a \ge 0$.
This foliation is called the \emph{$3$-dimensional non compact Reeb component} and it is denoted by $\mathcal{R}^3$ in \cite{Hector:1994}.
Since $F(x,y, t+1) = e\, F(x,y,t)$, then, $\mathcal{R}^3$ defines a foliation
on the quotient manifold $\D^2 \times \Sphere^1 = \D^2 \times \R \,/\, (x,y,t) \sim (x,y,t+1)$.
This foliation being called the \emph{$3$-dimensional compact Reeb component} and denoted $\mathcal{R}^3_c$.
This foliation admits a unique compact leaf: the boundary of $\D^2 \times \R$ diffeomorphic to the torus $\T^2 = \Sphere^1 \times \Sphere^1$.
All the others leaves are diffeomorphic to $\R^2$ and are also given by the quotient of the graphs of functions of the form (in polar
coordinates): $(r,\theta) \in \Ball(0, 1) \mapsto -\ln(\vphi(r^2)) + b$, $b \in \R$.

\begin{rmrk}
A simple example which is detailed in \cite{Hector:1994} is given by the mapping $f(r,\theta) = 1 - r^2$.
Note that it is possible to construct a Reeb component on the solid torus with a function $F$ which is a submersion only
on $\Ball(0, 1) \times \R$. Such an example is given by the function
$f(r,\theta) = \exp(-\exp(1/(1 - r^2)))$. 
This possibility will be useful in the following section
to construct a $2$-dimensional Reeb component from the Vortex application.
\end{rmrk}

Let us exhibit now a function $f$ in the Vortex application that can be used to construct a $3$-dimensional compact Reeb component
and thus a Reeb foliation on $\Sphere^3$.
To this end, we focus on the separating geodesics, that is the geodesics associated to parameters $\alpha \in \Psi(r_0, \mu)$.
These geodesics are parameterized by $p_\theta = p_\theta^* = -4\mu r_0^2/(r_0^2+4\mu^2)$ and $\norme{p(0)}_{r_0}=1$,
where $p(0) = (p_r(0), p_\theta)$ is the initial covector.
Let us fix the circulation $\mu > 0$ and the initial distance $r_0 > 0$.
Let $\alpha \in \Psi(r_0, \mu)$ denotes the single element contained in this set.
Let us denote by $(r(\cdot), \theta(\cdot))$ the associated separating geodesic and by $\gamma$ its orbit, that
is $\gamma \coloneqq \enstq{(r(t), \theta(t))}{t \in \intervalleoo{t^\alpha}{t_\alpha}}$, where $t^\alpha<0$ is the
minimal negative time such that the geodesic is well defined by backward integration. The time $t_\alpha>0$ being the
maximal positive time such that the geodesic is well defined by forward integration. Let us consider now an arbitrary point 
$(r_1, \theta_1) \in \gamma$ and define $\sigma$ as the orbit associated to the unique separating geodesic when we consider
$(r_1, \theta_1)$ as the new initial state. Then, it is clear that $\gamma = \sigma$ since again the set $\Psi(r_1, \mu)$ has only one element.
Thus, considering that at any time $t$, $r(t)$ is the initial distance to the vortex, leads to reparameterize the separating geodesics
by setting $p_\theta = -4\mu r(t)^2/(r(t)^2+4\mu^2)$ and $\norme{p(t)}_{r(t)}=1$, which gives 
\[
p_r(t)^2 = 1 - \frac{p_\theta^2}{r(t)^2} = \frac{(r(t)^2-4\mu^2)^2}{(r(t)^2+4\mu^2)^2}.
\]
In the case $0< r_0 < 2\mu$, the parameter $\alpha \in \Psi(r_0, \mu)$ is given by $\alpha_2^*$ and thus $p_r(t) > 0$ for any time $t$
(see Section \ref{sec:micro_local}). 
In this case, the dynamics of the separating geodesics reduces to the $2$-dimensional differential equations:
\begin{equation}
\dot{r}      = \frac{p_r}{\norme{p}_r} = p_r = \frac{4\mu^2-r^2}{4\mu^2+r^2}, \quad
\dot{\theta} = \frac{1}{r^2} \left( \mu + \frac{p_\theta}{\norme{p}_r} \right) = 
\frac{1}{r^2} \left( \mu + p_\theta \right) = \frac{\mu}{r^2} \left( \frac{4\mu^2-3r^2}{4\mu^2+r^2} \right).
\label{eq:separating_geo_dyn}
\end{equation}
Writing 
\[
\diff t = \frac{4\mu^2+r^2}{4\mu^2-r^2} \diff r,
\]
then by integration we have
\[
t(r) = 4\mu \atanh\left( \frac{r}{2\mu} \right) - r + c 
= 2\mu \ln \left( \frac{2\mu+r}{2\mu-r} \right) - r + c
= \ln \left( \left(\frac{2\mu+r}{2\mu-r}\right)^{2\mu} e^{-r} K \right),
\]
with a constant $c \coloneqq \ln K$. Introducing $2\mu\rho^2 \coloneqq r$, $\rho \in \intervalleff{0}{1}$, then one can define 
the function
\[
f(\rho, \theta) \coloneqq \vphi(\rho^2) \coloneqq \left(\frac{1-\rho^2}{1+\rho^2}\right)^{2\mu} e^{2\mu \rho^2}
\]
which satisfies the conditions (1)-(2)-(3) to construct a $3$-dimensional compact Reeb component from 
$F(\rho, \theta, t) \coloneqq f(\rho, \theta)\, e^t$.

\subsubsection{$2$-dimensional compact Reeb component}

In the Vortex application, one can find a foliation in the $2$-dimensional state space given by the separating geodesics.
To present this foliation, we need to introduce a generalized $2$-dimensional Reeb component following the presentation of 
\cite{Hector:1994}.
We thus first recall what is a $2$-dimensional Reeb component and then we introduce the generalization.
Let $f\colon \intervalleff{-1}{1} \to \R$ be a smooth mapping such that the following conditions are satisfied:
\begin{enumerate}
\setcounter{enumi}{3}
\item $f$ is of the form: $f(x) \eqqcolon \vphi(x^2)$,
\item $f(\{-1,1\}) = \{0\}$ and $f(\intervalleoo{-1}{1}) > 0$,
\item $f$ has no critical points on $\{-1,1\}$.
\end{enumerate}
Let consider the smooth mapping $F\colon \intervalleff{-1}{1} \times \R \to \R$ defined by $F(x, t) \coloneqq f(x)\, e^t$,
where $x \in \intervalleff{-1}{1}$ and $t\in \R$. $F$ is a submersion from $\intervalleff{-1}{1} \times \R$ to $\R$ and defines
a foliation of $\intervalleff{-1}{1} \times \R$ whose leaves are the level sets of $F$: $F^{-1}(a)$, $a \ge 0$.
This foliation is called the \emph{$2$-dimensional non compact Reeb component} and it is denoted by $\mathcal{R}^2$ in \cite{Hector:1994}.
Since $F(x, t+1) = e\, F(x, t)$, then, $\mathcal{R}^2$ defines a foliation on the annulus 
$\intervalleff{-1}{1} \times \Sphere^1 = \intervalleff{-1}{1} \times \R \,/\, (x, t) \sim (x, t+1)$.
This foliation being called the \emph{$2$-dimensional compact Reeb component} and denoted $\mathcal{R}^2_c$.
This foliation admits two compact leaves: the two boundaries of the annulus, diffeomorphic to $\Sphere^1$.
All the others leaves are diffeomorphic to $\R$ and each of their two extremities wind up around one of
the two compact leaves.
These non compact leaves are also given by the quotient of the graphs of functions of the form: 
$x \in \intervalleoo{-1}{1} \mapsto -\ln(\vphi(x^2)) + b$, $b \in \R$.



\begin{rmrk}
A simple example which is detailed in \cite{Hector:1994} is given by the function $f(x) = 1 - x^2$. The maximum
of this function is given by solving $f'(x)=2x=0$, that is for $x=0$. 
Note that it is possible to construct a Reeb component on the annulus with a function $F$ which
is a submersion only on $\intervalleoo{-1}{1} \times \R$ considering for instance
$f(x) = \exp(-\exp(1/(1-x^2)))$.
\end{rmrk}

\begin{rmrk}
Note that we can construct a $3$-dimensional Reeb component from the $2$-dimensional Reeb component
if $f$ satisfies the condition (4)-(5)-(6). See \cite{Hector:1994} for details.
\end{rmrk}

Let us slightly generalize the construction of a Reeb component on the annulus. To do so, we consider
a function $f \colon \intervalleff{a}{b} \to \R$ smooth on $\intervalleoo{a}{b}$, $a<b$ in $\R$,
and such that:
\begin{enumerate}
\setcounter{enumi}{6}
\item $f(\{a,b\}) = \{0\}$, $f(\intervalleoo{a}{b})>0$.
\end{enumerate}
Then, $F(x,t) \coloneqq f(x)\, e^t$ is a submersion from $\intervalleoo{a}{b} \times \R$ to $\R$
and defines a foliation on $\intervalleff{a}{b} \times \R$ which is called the 
\emph{generalized $2$-dimensional non compact Reeb component} and denoted $\mathcal{R}^{2,G}$. 
As in the two previous cases, $\mathcal{R}^{2,G}$ defines a foliation on the annulus 
$\intervalleff{a}{b} \times \Sphere^1$ which is called the \emph{generalized $2$-dimensional compact Reeb
component} and denoted $\mathcal{R}_c^{2,G}$. This foliation, like $\mathcal{R}_c^{2}$, has
two compact leaves diffeomorphic to $\Sphere^1$ and all the others leaves are diffeomorphic to $\R$
and given by the quotient of the graphs of functions of the form:
$x \in \intervalleoo{a}{b} \mapsto -\ln(f(x)) + c$, $c \in \R$.

\begin{rmrk}
The most important difference with the classical case is that we do not impose that $f(x)$ is
of the form $\vphi(x^2)$. Even in the case $a=-b$, if $f(x)$ is not of the form $\vphi(x^2)$,
then we do not have necessarily $f(-x)=f(x)$ and if we do not have this symmetry, then it is
not possible to construct a $3$-dimensional Reeb component from a $2$-dimensional one using
the construction detailed in \cite{Hector:1994}.
\end{rmrk}

\begin{figure}[ht!]
\def\sizefigreeb{0.3}
\centering
\includegraphics[width=\sizefigreeb\textwidth]{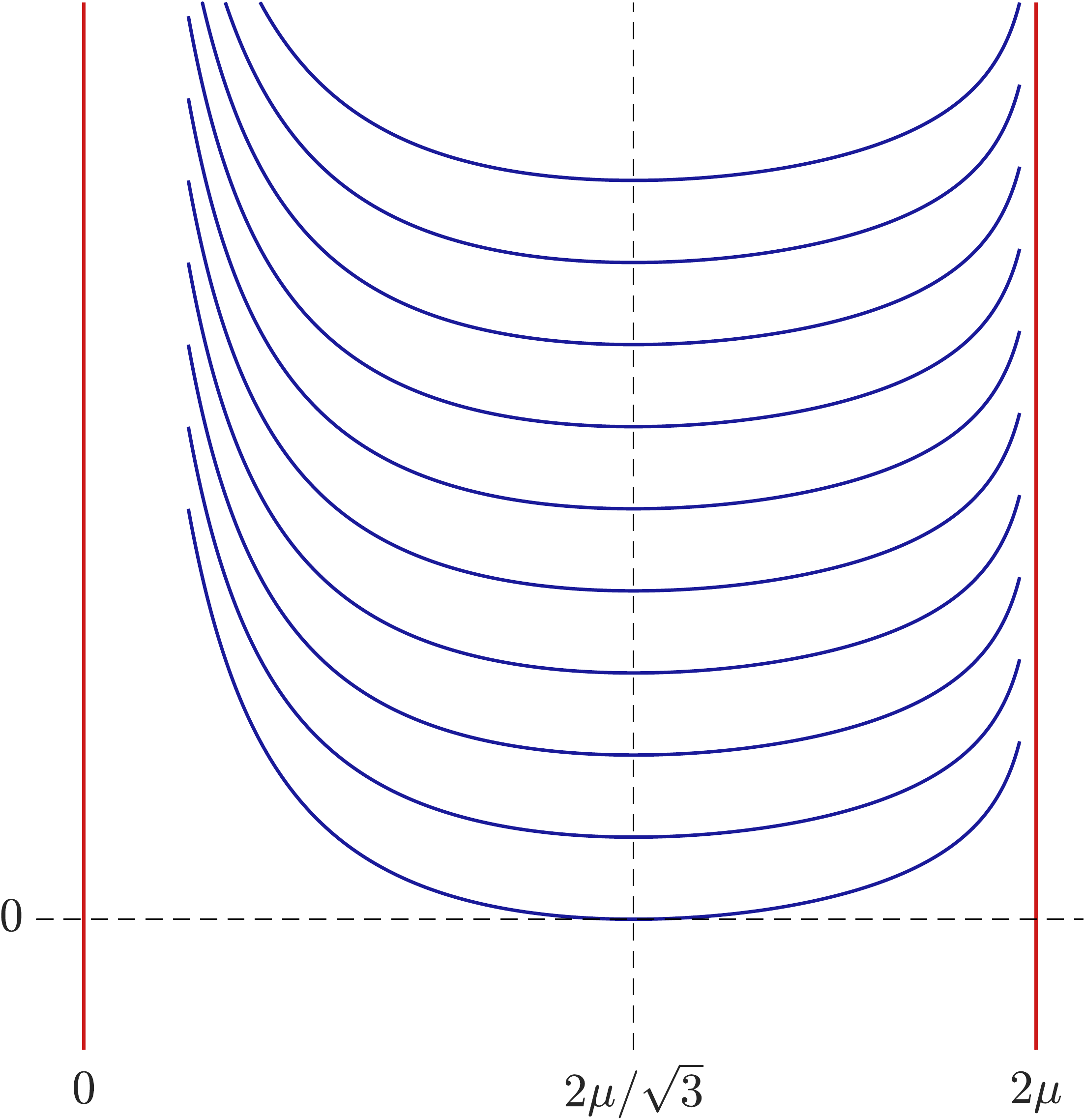}
\hspace{3.0em}
\includegraphics[width=\sizefigreeb\textwidth]{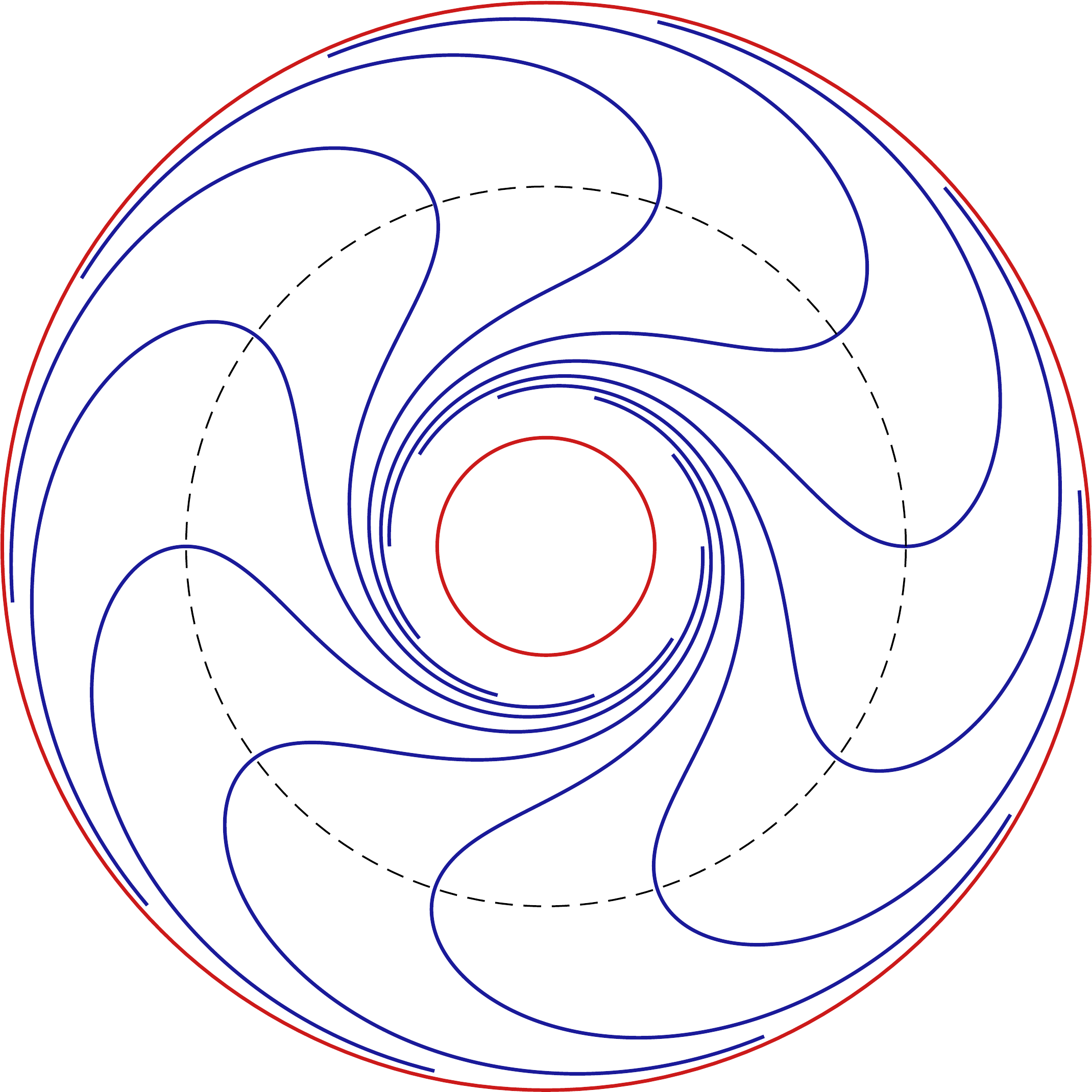}
\caption{(Left) Graphs of $r \in \intervalleoo{0}{2\mu} \mapsto -\ln(f(r)) + c$, $c \in \R$,
with $f$ defined by eq.~\eqref{eq:f_Reeb_Vortex}.
The minimum is attained at $r=2\mu/\sqrt{3}$ where $f'(r)=0$.
(Right) Generalized $2$-dimensional compact Reeb component of the annulus given by the function
$F(r,\theta) = f(r)\, e^\theta$. The intermediate dashed circle of radius $2\mu/\sqrt{3}$ is orthogonal to all
the non compact leaves in blue. The two compact leaves are the two red circles.}
\label{fig:Reeb_vortex_1}
\end{figure}

Let us exhibit now a function $f$ from the Vortex application that satisfies the condition (7).
We fix $\mu > 0$.
Coming back to eq.~\eqref{eq:separating_geo_dyn}, one can write
\[
\diff \theta = \frac{\mu}{r^2} \left( \frac{4\mu^2-3r^2}{4\mu^2-r^2} \right) \diff r.
\]
Integrating we have:
\[
\theta(r) = \atanh\left( \frac{r}{2\mu} \right) + \frac{\mu}{r} + c 
= \frac{1}{2} \ln \left( \frac{2\mu+r}{2\mu-r} \right) +  \frac{\mu}{r}+ c
= \ln \left( \left(\frac{2\mu+r}{2\mu-r}\right)^{\frac{1}{2}} e^{\frac{\mu}{r}} K \right),
\]
with a constant $c \coloneqq \ln K$. Introducing for $r \in \intervalleff{0}{2\mu}$:
\begin{equation}
f(r) \coloneqq \left(\frac{2\mu-r}{2\mu+r}\right)^{\frac{1}{2}} e^{-\frac{\mu}{r}},
\label{eq:f_Reeb_Vortex}
\end{equation}
with $f(0) = 0$, then the separating geodesics inside the punctured disk of radius $2\mu$ are exactly given
by the level sets of the function $F(r, \theta) \coloneqq f(r)\, e^\theta$.
Besides, the function $f$ satisfies the condition (7) and thus, it defines a generalized
$2$-dimensional compact Reeb component on the annulus $\intervalleff{0}{2\mu} \times \Sphere^1$, see 
Figure~\ref{fig:Reeb_vortex_1}.

\begin{figure}[ht!]
\def\sizefigreeb{0.3}
\centering
\includegraphics[width=\sizefigreeb\textwidth]{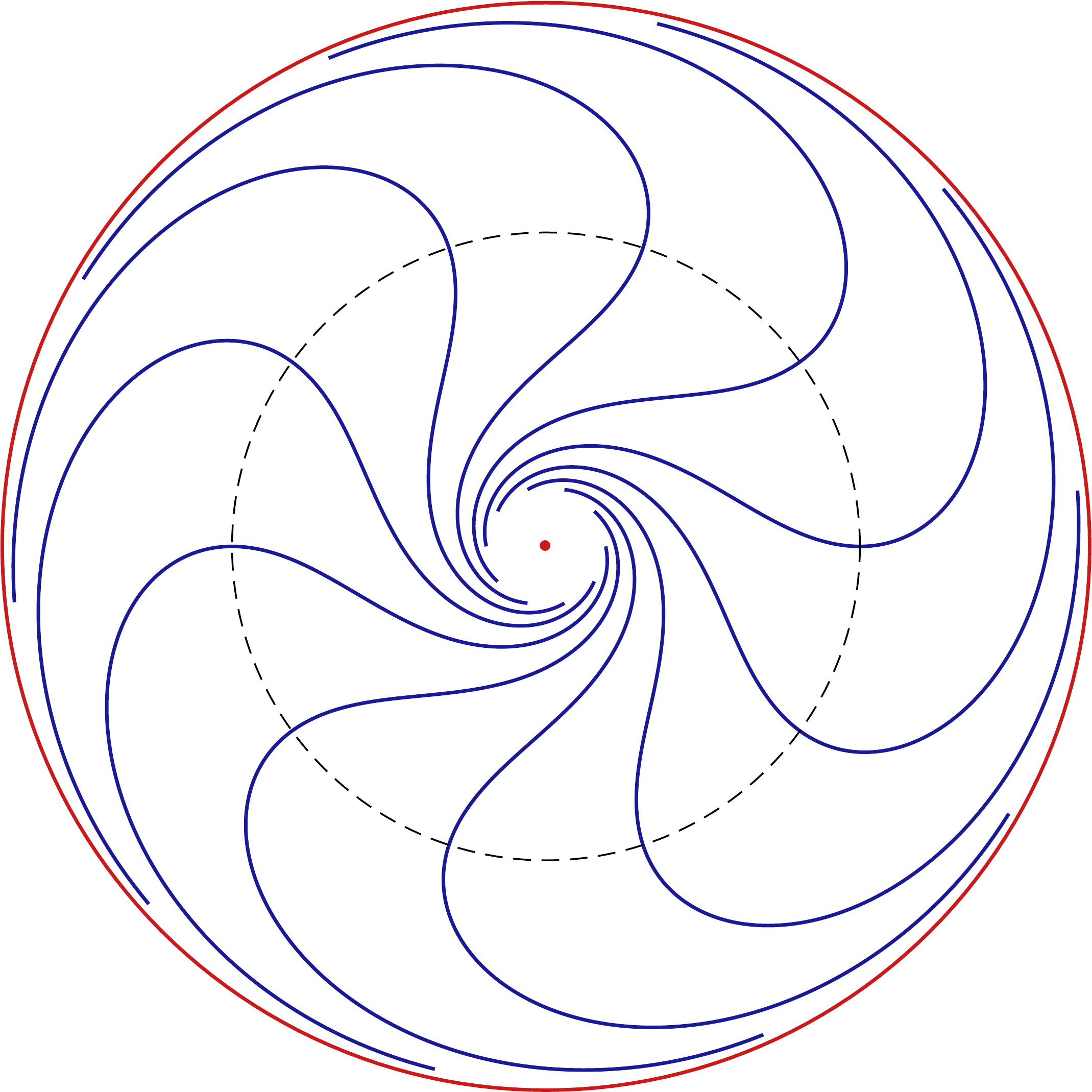}
\caption{Separating geodesics in the punctured disk of radius $2\mu$ in the $(x_1,x_2)$-plane.
The vortex is placed at the origin and represented by a red dot. The red circle is the Reeb circle of radius $2\mu$
while the black and dashed circle is the circle of radius $2\mu/\sqrt{3}$ where $\dot{\theta}=0$ along the separating
geodesics.}
\label{fig:Reeb_vortex_2}
\end{figure}

The separating geodesics inside the punctured disk of radius $2\mu$ correspond to the non compact leaves of the foliation
while the unique circle geodesic, that is the separating geodesic with $r_0 = 2\mu$, corresponds to one of the two compact
leaves. The vortex corresponds to the second compact leaf. Figure~\ref{fig:Reeb_vortex_2} shows the foliation on the
punctured disk of radius $2\mu$ in the $(x_1,x_2)$-plane of the Vortex application. 

\begin{dfntn}
In the Vortex application, the circle of radius $2\abs{\mu}$ (assuming $\mu\ne0$) is called the \emph{Reeb circle}.
\label{def:Reeb_circle}
\end{dfntn}

\subsection{Symmetries of the extremal curves}
\label{sec:Symmetries}

There exists three continuous symmetries and one discrete symmetry of main interest. The first continuous symmetry was considered at the
beginning of the paper when we have fixed $\umax = 1$. The second symmetry is the rotational symmetry that was mentioned
at the end of Section \ref{sec:classification}. As a consequence of this symmetry of revolution, one can easily prove
that:
\[
\forall \Delta\theta \in \R ~:~ V(r_0, \theta_0, r_f, \theta_f, \mu) = 
V(r_0, \theta_0 + \Delta \theta, r_f, \theta_f + \Delta \theta, \mu),
\]
where $V$ is the value function \eqref{eq:value_function} expressed in polar coordinates. Hence, one can fix for instance
$\theta_0 = 0$ and reduce the set of variables of the value function. Note that the Reeb circle is the unique geodesic invariant
by the rotational symmetry. The third continuous symmetry comes from the invariance of the set of geodesics
by homothety combined with time/circulation dilatation. Let us detail this third symmetry. We start by introducing the
change of variables
$
(\rt, \thetat) \coloneqq (\lambda r, \theta)
$
for a given $\lambda \in \Rsp$. The associated Mathieu transformation leads to
\[
\pt \coloneqq (p_{\rt}, p_{\thetat}) = \left( \frac{p_r}{\lambda}, p_\theta \right)
\]
and the Hamiltonian system \eqref{eq:hamiltonian_system_polar} in polar coordinates becomes
\begin{equation*}
\begin{aligned}
\dot{\rt}           &= \lambda \dot{r}              = \frac{\lambda\, p_r}{\norme{p}_r} 
= \frac{\lambda\, p_r}{\lambda \norme{\pt}_{\rt}}
= \lambda \frac{p_{\rt}}{\norme{\pt}_{\rt}}, \\
\dot{\thetat}       &= \dot{\theta}                 = \frac{1}{r^2} \left( \mu + \frac{p_\theta}{\norme{p}_r} \right)
= \frac{\lambda^2}{\rt^2} \left( \mu + \frac{p_{\thetat}}{\lambda\norme{\pt}_{\rt}} \right)
= \frac{\lambda}{\rt^2} \left( \lambda \mu + \frac{p_{\thetat}}{\norme{\pt}_{\rt}} \right)
\end{aligned}
\end{equation*}
for the state and
\begin{equation*}
\dot{p}_{\rt}         = \frac{\dot{p}_r}{\lambda}    
= \frac{p_\theta}{\lambda r^3} \left( 2 \mu + \frac{p_\theta}{\norme{p}_r} \right)
= \frac{\lambda^2 p_{\thetat}}{\rt^3} \left( 2 \mu + \frac{p_{\thetat}}{\lambda \norme{\pt}_{\rt}} \right)
= \frac{\lambda p_{\thetat}}{\rt^3} \left( 2 \lambda \mu + \frac{p_{\thetat}}{\norme{\pt}_{\rt}} \right), \quad
\dot{p}_{\thetat}    = 0,
\end{equation*}
for the covector, since
\[
\norme{\pt}_{\rt}^2 = \normeStyle{\left(\frac{p_r}{\lambda}, p_\theta\right)}_{\rt}^2 = \frac{p_r^2}{\lambda^2}+\frac{p_\theta^2}{\lambda^2r^2}
= \frac{1}{\lambda^2}\norme{p}_r^2.
\]
Introducing now $\mub \coloneqq \lambda \mu$, the time reparameterization $t = \psi(\tau) \eqqcolon {\tau}/{\lambda}$,
and the notation
$
(\rb, \thetab, \pb_r, \pb_\theta)(\tau) \coloneqq (\rt, \thetat, p_{\rt}, p_{\thetat})(\psi(\tau)),
$
then the system in the time $\tau$ writes
\begin{equation*}
\begin{aligned}
\frac{\diff}{\diff \tau}{\rb(\tau)} &= \frac{\diff}{\diff \tau}{\rt(\psi(\tau))} = \dot{\rt}(\psi(\tau)) \psi'(\tau) = \frac{\dot{\rt}(t)}{\lambda}
=  \frac{\pb_r(\tau)}{\norme{\pb(\tau)}_{\rb(\tau)}}, \\
\frac{\diff}{\diff \tau}{\thetab} \hspace{0.15cm}  &= \frac{1}{\rb^2} \left( \mub + \frac{\pb_\theta}{\norme{\pb}_{\rb}} \right), \\
\frac{\diff}{\diff \tau}{\pb}_r       &= \frac{\pb_\theta}{\rb^3} \left( 2 \mub + \frac{\pb_\theta}{\norme{\pb}_{\rb}} \right), \\
\frac{\diff}{\diff \tau}{\pb}_\theta  &= 0, \\
\end{aligned}
\end{equation*}
which is equivalent to the original system \eqref{eq:hamiltonian_system_polar}.
As a conclusion, the set of solutions of the problem \eqref{eq:minimum_time} is invariant by the dilatation
$
(r(\cdot), \theta(\cdot), t, \mu) \mapsto \sigma_\lambda(r(\cdot), \theta(\cdot), t, \mu) \coloneqq 
(\lambda r(\cdot), \theta(\cdot), \lambda t, \lambda \mu),
$
for any $\lambda \in \Rsp$, and we have
\[
\frac{1}{\lambda} V(\lambda r_0, \theta_0, \lambda r_f, \theta_f, \lambda \mu) = V(r_0, \theta_0, r_f, \theta_f, \mu).
\]
It is thus possible to fix for instance $r_0$ and again to reduce the set of variables of the value function.

\begin{rmrk}
Setting $\lambda \coloneqq 1/r_0$ and considering also the invariance by rotation, we can write:
\[
V(r_0, \theta_0, r_f, \theta_f, \mu) = r_0 V(1, 0, \Delta r, \Delta \theta, \frac{\mu}{r_0}),
\]
with $\Delta r \coloneqq r_f/r_0$ and $\Delta \theta \coloneqq \theta_f - \theta_0$. This shows that the value function depends
in fact only on the three variables $\Delta r$, $\Delta \theta$ and $\mu/r_0$.
\end{rmrk}

\begin{figure}[ht!]
\centering
\includegraphics[width=0.3\textwidth]{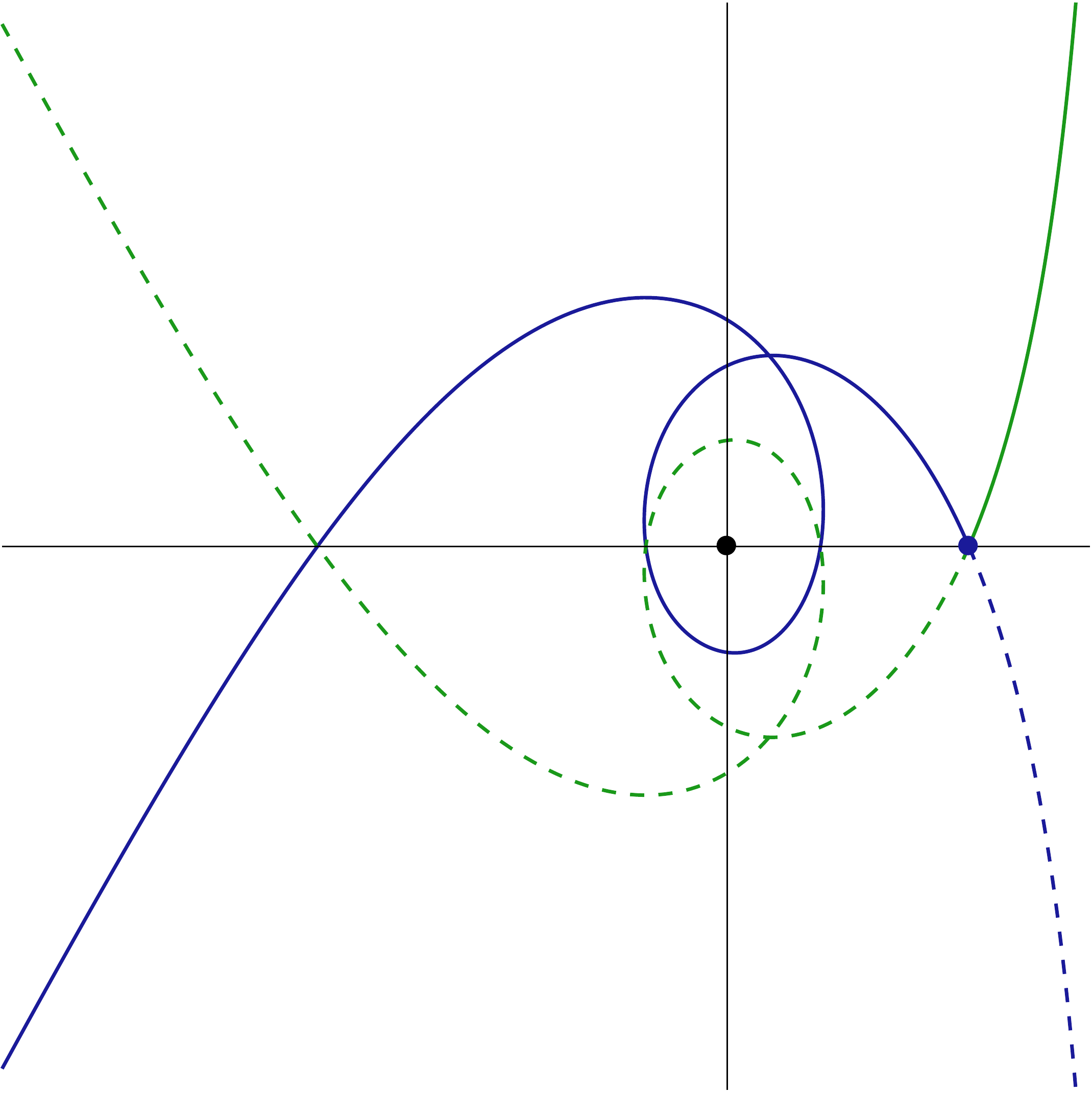}
\hspace{3.0em}
\includegraphics[width=0.3\textwidth]{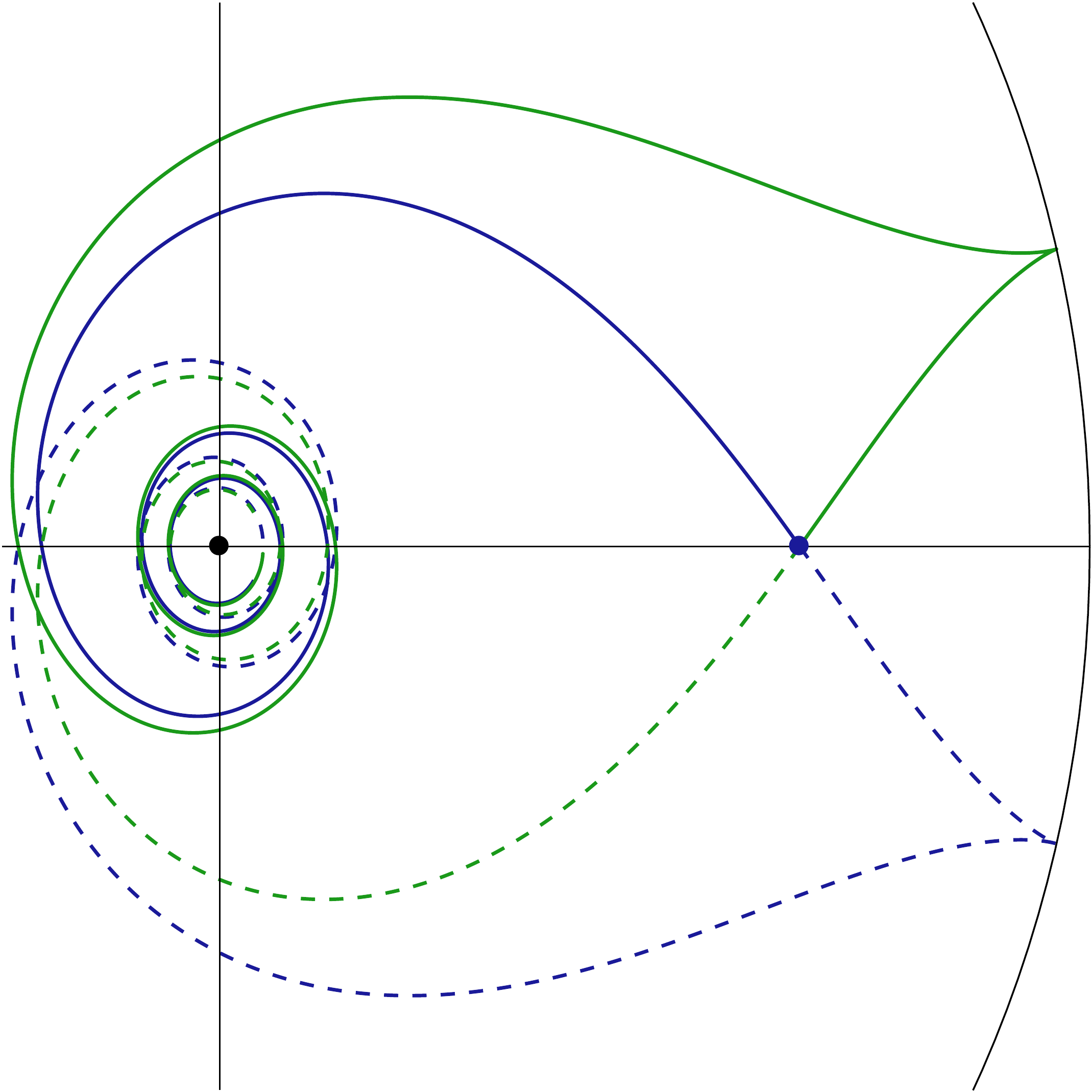}
\caption{\textbf{Geodesics in the state space in the $(x_1,x_2)$-plane.}
The vortex is placed at the origin and represented by a black dot. The initial point, denoted $x_0$, is the black dot at
an intersection of both geodesics. The geodesic in green is obtained by a reflection of axis $\R x_0$ from the blue 
geodesic and vice-versa. Assuming that the initial point is given at time $t=0$, then the plain lines are obtained 
for times $t>0$ while the dashed parts are for negative times. On the left subgraph is represented a normal geodesic 
with both extremities going to infinity while on the right subgraph is represented the two abnormals (in the case of 
a strong drift) with a cusp at $r=\abs{\mu}$. Both extremities of the abnormal geodesics go the vortex.}
\label{fig:discrete_symmetry}
\end{figure}

We end this section presenting the discrete symmetry on the set of extremals.
Let $(r(\cdot), \theta(\cdot), p_r(\cdot), p_\theta)$ be a reference extremal on $\intervalleff{0}{T}$
and introduce $\theta_0 \coloneqq \theta(0)$.
Let consider the following discrete symmetry:
$\zt(\cdot) \coloneqq (\rt(\cdot), \thetat(\cdot), p_{\rt}(\cdot), p_{\thetat})
\coloneqq (r(\cdot), 2\theta_0-\theta(\cdot), -p_r(\cdot), p_\theta)$.
Then, $\zt(\cdot)$ is solution of 
\[
\dot{z} = -\vv{\Htrue}(z), \quad z(0) = (r(0), \theta_0, -p_r(0), p_\theta).
\]
This means that $\zt(\cdot)$ is obtained by backward integration from the same initial condition in the state space 
than $z(\cdot)$ but with the initial covector $(-p_r(0), p_\theta)$. Let us extend the reference extremal to obtain 
a maximal solution of $\dot{z} = \vv{\Htrue}(z)$, $z(0) = (r(0), \theta_0, p_r(0), p_\theta)$.
In this case, its orbit in the state space is given by $\gamma \coloneqq \im (r(\cdot), \theta(\cdot))$
and $\zt$ is also a maximal solution of its associated differential equation, and we denote by 
$\gammat \coloneqq \im (\rt(\cdot), \thetat(\cdot))$ its orbit in the state space.
Then, we have:
\[
\begin{pmatrix}
1 & 0 \\ 0 & -1
\end{pmatrix}
\gamma
+
\begin{pmatrix}
0 \\ 2\theta_0
\end{pmatrix}
=
\gammat,
\]
that is $\gammat$ is obtained by an affine reflection of axis $\theta = \theta_0$. In cartesian coordinates, one has
a linear reflection of axis $\R x_0$, see Figure~\ref{fig:discrete_symmetry}.



\subsection{Properties of the value function and its level sets in the weak case}

\subsubsection{Continuity of the value function and characterization of the cut points}

Let us introduce the mapping $V_\mu(x_0, x_f) \coloneqq V(x_0, x_f, \mu)$, where $V$ is the value function defined by \eqref{eq:value_function}.
We are first interested in the continuity of $x_f \mapsto V_\mu(x_0, x_f)$ at points $x_f$ where the drift is weak, that is for target 
$x_f$ such that $\norme{x_f} > \abs{\mu}$. This situation is more general than the Finslerian case with Randers metric since the drift
is not necessarily weak all along the geodesics. We have:

\begin{prpstn}
\label{prop:continuity_weak_drift}
Let $\mu \ne 0$ be the vortex circulation and $x_0 \in \xdom$ the initial condition.
Then, $V_{\mu}(x_0, \cdot)$ is continuous at any point $x_f$ such that $\norme{x_f} > \abs{\mu}$.
\end{prpstn} 

\begin{proof}
Since we consider that the drift is weak at the target $x_f$ then local controllability around $x_f$ holds, that is we have the following:
\begin{equation}
\label{eq:local_controllability}
\forall\, \veps > 0,~ \exists\, \eta_\veps > 0 \text{~s.t.~} \forall\, (x_1, x_2) \in 
\Ball(x_f, \eta_\veps)\times\Ball(x_f, \eta_\veps) :
V_\mu(x_1, x_2) \le \veps.
\end{equation}
Let us fix $\veps > 0$ and consider $x \in \Ball(x_f, \eta_\veps)$. By definition of the value function and by
\eqref{eq:local_controllability}, we have
\[
V_\mu(x_0, x) - V_\mu(x_0, x_f) \le V_\mu(x_f, x) \le \veps \quad \text{and} \quad
V_\mu(x_0, x_f) - V_\mu(x_0, x) \le V_\mu(x, x_f) \le \veps,
\]
which prove the result.
\end{proof}

We are now interested in the construction of the optimal synthesis associated to problem \eqref{eq:minimum_time} for given
values of $\mu$ and $x_0$, that is we want to get for any $x_f \in \xdom$ the value of $V_\mu(x_0, x_f)$ together with the optimal
geodesic  joining the points $x_0$ and $x_f$. To construct such an optimal synthesis, a classical approach is to compute the
level sets of $V_\mu(x_0, \cdot)$, which corresponds in Riemmanian or Finslerian geometry to compute the spheres centered in $x_0$.
To compute the spheres, one preliminary work is to compute the \emph{cut locus}, that is the set of 
\emph{cut points} where the geodesics cease to be optimal. Let us introduce some notations before characterizing the cut points.
Following Section \ref{sec:abnormal_extremals} and fixing $x_0 \in \xdom$, we introduce the notation in polar coordinates
$p_0(\alpha) \coloneqq (\cos \alpha, r_0 \sin \alpha)$, with $\alpha \in \intervallefo{0}{2\pi}$ and where 
$r_0 \coloneqq \norme{x_0}$.
%
%
%
Let $t \in \Rp$, we introduce the \emph{wavefront} from $x_0$ at time $t$ by
\[
\WF(x_0, t) \coloneqq \enstq{\exp_{x_0}(t, p_0(\alpha))}{\alpha \in \intervallefo{0}{2\pi}
~\text{s.t.}~ t < t_\alpha},
\]
where the exponential mapping is given by \eqref{eq:exponential_mapping} and depends on the circulation $\mu$ supposed to be given,
and where we recall that $t_\alpha \in \Rsp \cup \{+\infty\}$ is the maximal positive time such that the associated geodesic is
defined over $\intervallefo{0}{t_\alpha}$.
Then, we define what corresponds to the balls and spheres in Riemannian or Finslerian geometry:
\[
\begin{aligned}
\MBall(x_0, t)       \coloneqq \enstq{x \in \xdom}{V_\mu(x_0, x) < t}, \\
\MBallClosed(x_0, t) \coloneqq \enstq{x \in \xdom}{V_\mu(x_0, x) \le t}, \\
\MSphere(x_0, t)     \coloneqq \enstq{x \in \xdom}{V_\mu(x_0, x) = t}. \\
\end{aligned}
\]
We use the terminology ball and sphere by analogy with Geometry. Note that the spheres $\MSphere(x_0, \cdot)$ are the level
sets of $V_\mu(x_0, \cdot)$. To construct the optimal synthesis, one of the main important loci to compute is the \emph{cut locus}.
We define first the \emph{cut times} as:
\[
\tcut(x_0, \alpha) \coloneqq \sup \enstq{t \in \Rsp}{\exp_{x_0}(\cdot, p_0(\alpha)) 
~\text{is optimal over}~ \intervalleff{0}{t}}.
\]
Then, the \emph{cut locus} from $x_0$ is simply defined by:
\[
\Cut(x_0) \coloneqq \enstq{x \in \xdom}{\exists\, (t, \alpha) \in \Rsp \times \intervallefo{0}{2\pi} ~\text{s.t.}~
x = \exp_{x_0}(t, p_0(\alpha)) ~\text{and}~ t = \tcut(x_0, \alpha)}.
\]
Any point of the cut locus is called a \emph{cut point}. From the cut times, we define also the \emph{injectivity radius}
from $x_0$ as:
\[
\tinj(x_0) \coloneqq \inf_{\alpha \in \intervallefo{0}{2\pi}} \tcut(x_0, \alpha).
\]
\begin{rmrk}
This time is of particular interest since if the drift is weak at the initial point $x_0$, then $\exp_{x_0}$ is a 
diffeomorphism from $\intervalleoo{0}{\tinj(x_0)} \times \Sphere^1$ to 
$\MBall(x_0, \tinj(x_0)) \setminus \{x_0\}$. This is similar to the Riemannian and Finslerian cases
and for times $0 \le t \le \tinj(x_0)$, we have $\MSphere(x_0, t) = \WF(x_0, t)$.
\end{rmrk}
%
In Finslerian geometry \cite{BaoChernShen} (as in Riemannian geometry \cite[Proposition 2.2, Chapter 13]{Carmo:1988}), 
the cut locus is a part of the union of the conjugate locus and the \emph{splitting locus}, which is defined by:
\[
\Split(x_0) \coloneqq \enstq{x \in \xdom}{\exists (t, \alpha_1, \alpha_2) \in \Rsp \times 
\intervallefo{0}{2\pi}^2 ~\text{s.t.}~
\alpha_1 \ne \alpha_2 ~\text{and}~ \exp_{x_0}(t, p_0(\alpha_1)) = \exp_{x_0}(t, p_0(\alpha_2)) = x}.
\]
We have the following similar result:
\begin{prpstn}
\label{prop:caracterisation_cut_points}
Let $x(\cdot)$ be a geodesic. If $x_c \coloneqq x(t_c)$ is a cut point along $x(\cdot)$ such that
$V_\mu(x_0, \cdot)$ is continuous at $x_c$, then:
\begin{enumerate}[label=(\alph*)]
\item either $x_c$ is a conjugate point along $x(\cdot)$,
\item or $x_c$ is a splitting point along $x(\cdot)$.
\end{enumerate}
Conversely if a) or b) holds, then there exists a time $0 \le t_1 \leq t_c $ such that $x(t_1)$ is a cut point along $x(\cdot)$.
\end{prpstn}

\begin{rmrk}
The proposition above gives a characterization of the cut points where the mapping $V_\mu(x_0, \cdot)$ is continuous.
In Riemannian and Finslerian geometry, the mapping $V_\mu(x_0, \cdot)$ is replaced by the distance function which is continuous
with respect to its second argument and so all the cut points are either conjugate or splitting points. By Proposition
\ref{prop:continuity_weak_drift}, if the drift is weak at the cut points, then $V_\mu(x_0, \cdot)$ is continuous
at this point and so the characterization holds. In the strong case, we can have abnormal minimizers and so we can have
cut points which are neither conjugate nor splitting points. However, if the drift is weak at the initial point $x_0$,
then there are no abnormals and we can expect that the mapping $V_\mu(x_0, \cdot)$ is continuous at any point, even where the
drift is strong. This result would be useful to conclude in the following part since we are interested in the construction 
of the optimal synthesis for a given initial condition $x_0$ where the drift is weak. Note that the numerical experiments
of the following section suggest that $V_\mu(x_0, \cdot)$ is continuous at any point of $\xdom$.
\end{rmrk}

\subsubsection{Level sets of the value function and optimal synthesis}

In this section, we fix for the numerical experiments $x_0 \coloneqq (3,0)$ and $\mu \coloneqq 0.6\, \norme{x_0}$.
The drift is thus weak at the initial position $x_0$. We decompose the construction of the optimal synthesis in three steps.
In the first step, we compute the splitting locus. In a second time we give the cut locus and we finish by the construction
of the spheres and balls.

\medskip\paragraph{\textbf{Step 1: Computation of the splitting locus.}} Let us introduce the mapping
\[
F_\mathrm{split}(t, \alpha_1, x, \alpha_2) \coloneqq (x - \exp_{x_0}(t, p_0(\alpha_1)), x - \exp_{x_0}(t, p_0(\alpha_2)))
\in \R^4.
\]
The splitting line is then given by solving $F_\mathrm{split} = 0$ since we have
\[
\Split(x_0) = \enstq{x \in \xdom}{ \exists\, (t, \alpha_1, \alpha_2) ~\text{s.t.}~
\alpha_1 \ne \alpha_2 ~\text{and}~ F_\mathrm{split}(t, \alpha_1, x, \alpha_2) = 0}.
\]
Introducing $y \coloneqq (t, \alpha_1, x) \in \R^4$ and $\lambda \coloneqq \alpha_2 \in \R$, we have to solve 
$F_\mathrm{split}(y, \lambda) = 0 \in \R^4$. Numerically, we compute the splitting line with the differential path following
method (or homotopy method) of the \hampath\ software.
Under some regularity assumptions, the set $F_\mathrm{split}^{-1}(0)$ is a disjoint union
of differential curves, each curve being called a \emph{path of zeros}. 
To compute a path of zeros, we look for a first point on the curve by fixing the homotopic parameter $\lambda$ to a certain 
value $\lb$ and calling a Newton method to solve  $F_\mathrm{split}(\cdot, \lb)=0$. Then, we use a Prediction-Correction (PC)
method to compute the differential curve. The \hampath\ code implements a PC method with a high-order Runge-Kutta scheme
with variable step-size for the prediction, hence the discretization grid of the homotopic parameter is computed by the numerical
integrator. Besides, the Jacobian of $F_\mathrm{split}$ is computed by automatic differentiation combined with variational equations
of the exponential mapping.

\begin{rmrk}
If there exists several paths of zeros, then each path has to be found manually and compared.
\end{rmrk}

\begin{figure}[ht!]
\centering
\def\x{1425}
\def\y{1058}
\includegraphics[width=0.12\textwidth]{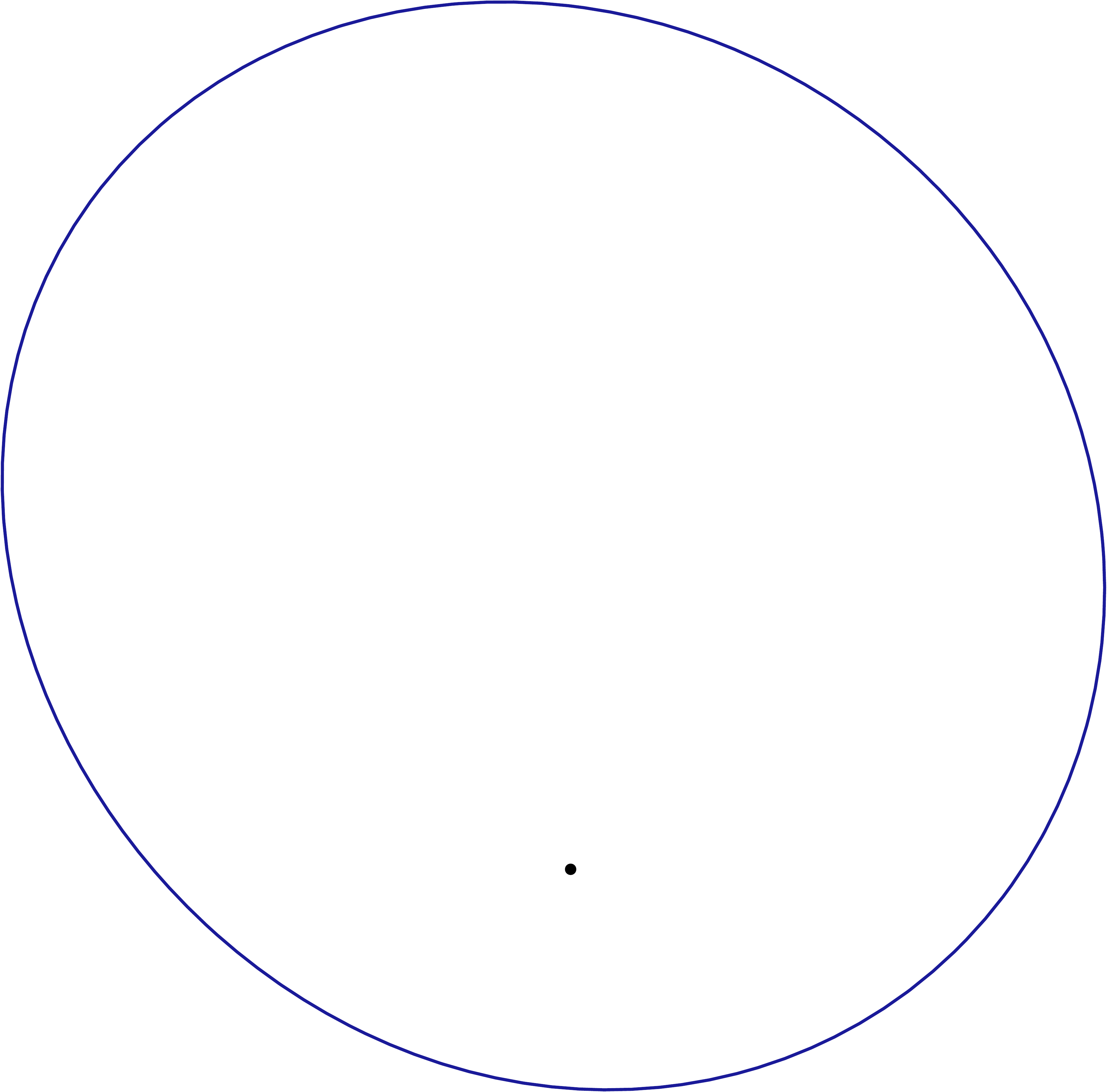}
\hspace{0.0em}
\includegraphics[width=0.20\textwidth]{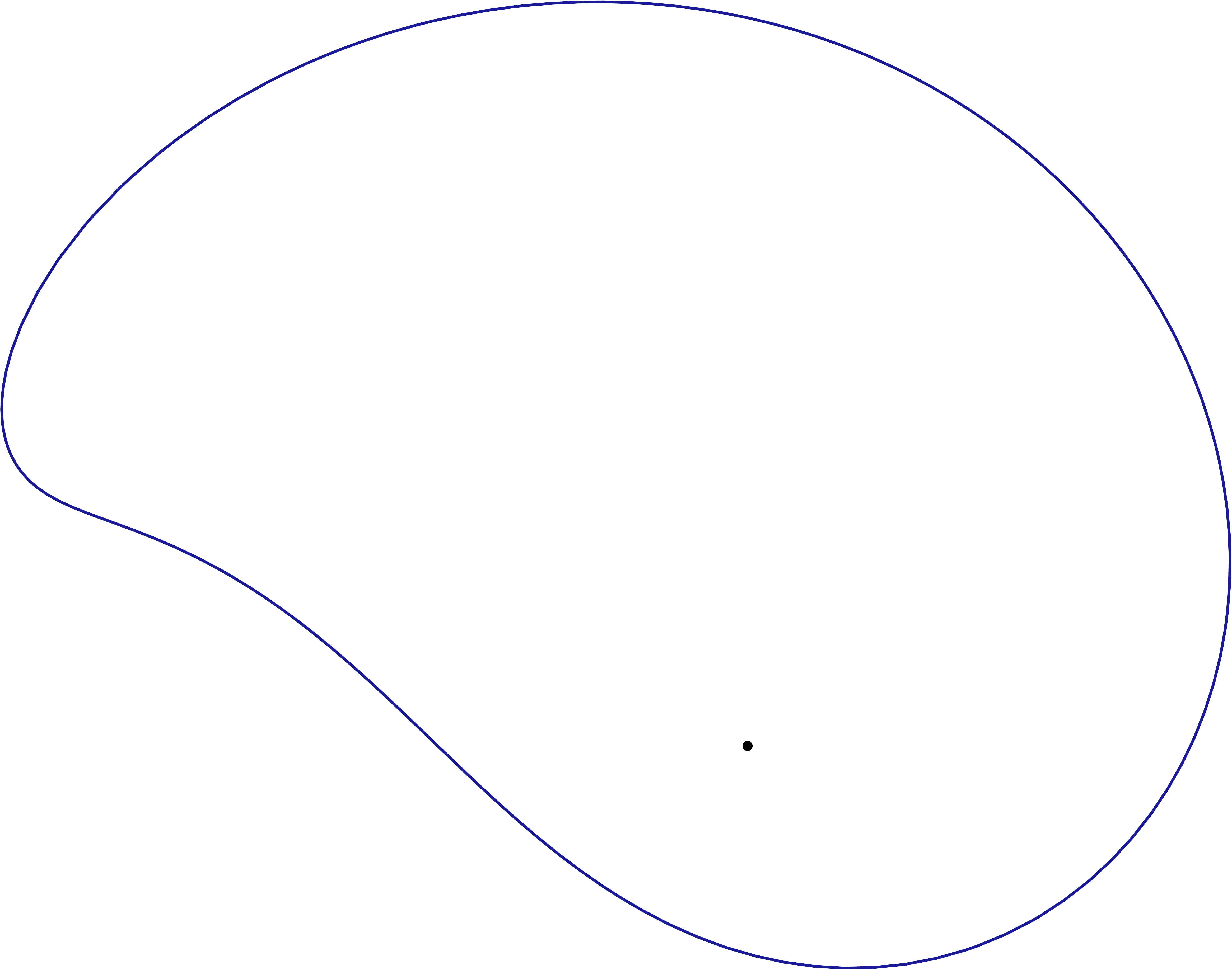}
\hspace{0.0em}
\includegraphics[width=0.28\textwidth]{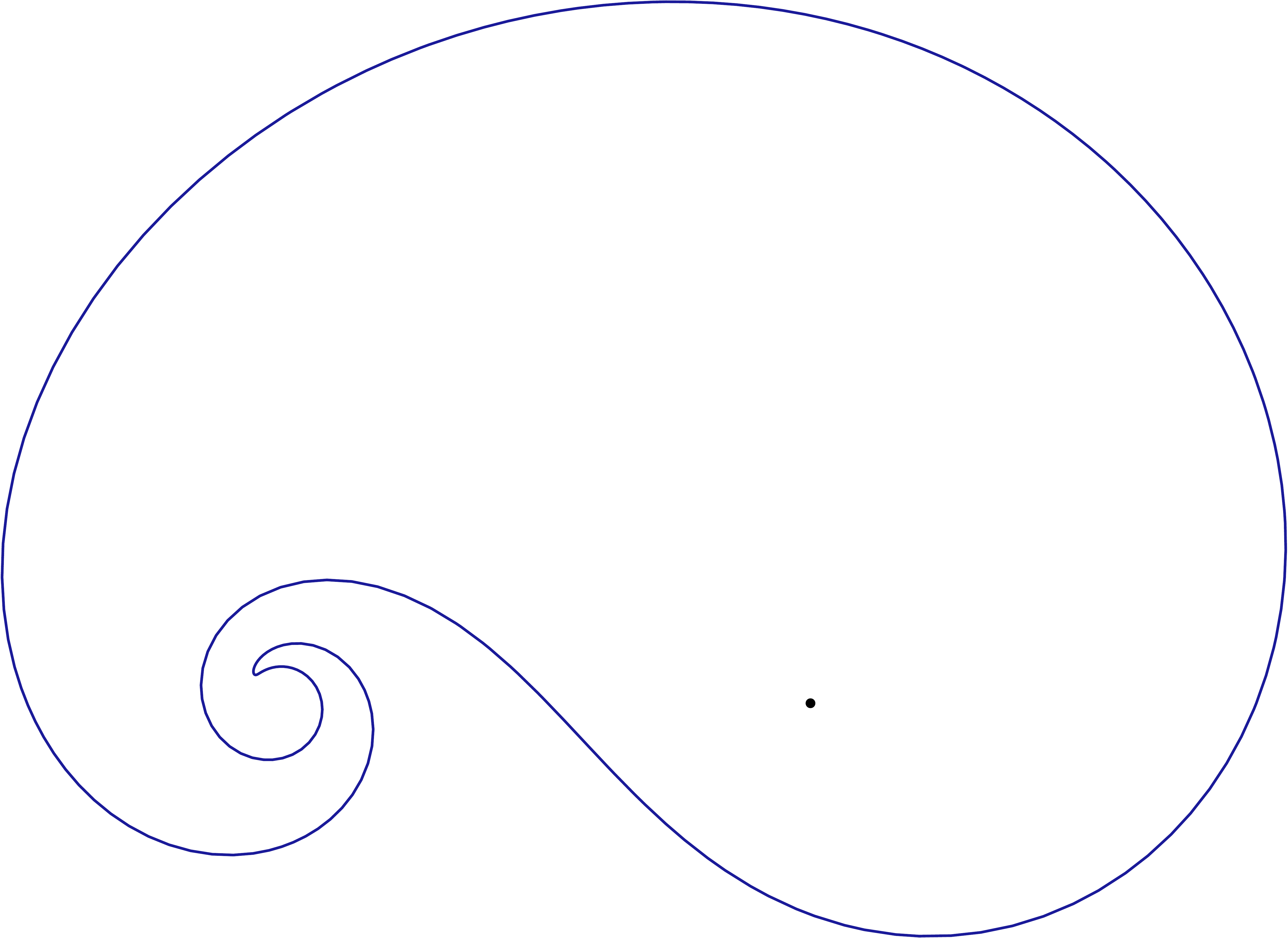}
\hspace{0.0em}
\begin{tikzgraphics}{0.34\textwidth}{\x}{\y}{WF_3p500}
\pxcoordinate{0.46*\x}{0.73*\y}{A}; \draw (A) node {\tiny $1$};
\pxcoordinate{0.39*\x}{0.64*\y}{A}; \draw (A) node {\tiny $2$};
\pxcoordinate{0.335*\x}{0.665*\y}{A}; \draw (A) node {\tiny $3$};
\end{tikzgraphics}
\caption{Four wavefronts with different scaling at times $t \in \{0.5, 2, 2.8, 3.5\}$. 
The black dot is $x_0$. For $t = 3.5$, the wavefront has at least three self-intersections labelled $1$, $2$ and $3$.}
\label{fig:wavefronts}
\end{figure}

To compute a splitting curve (that is a path of zeros), we need to get a first point on the curve.
This step is easy since the splitting points are self-intersections of the wavefronts.
On the right subgraph of Figure~\ref{fig:wavefronts}, we can see that the wavefront\footnote{The wavefronts are
computed by homotopy. We compute only rough approximations of the extremities of the wavefronts which are not closed curves, 
like $\WF(x_0, 3.5)$. Hence, there may exists others self-intersections but which are not relevant for rest of the analysis.}
$\WF(x_0, 3.5)$ has at least three self-intersections labelled $1$, $2$ and $3$. 

\begin{figure}[ht!]
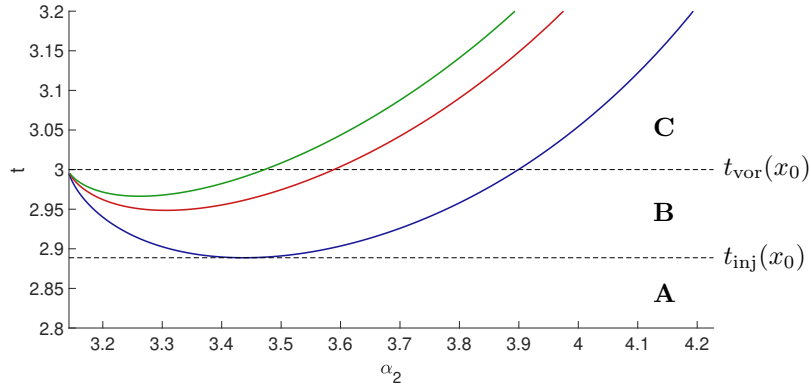

\centering
\def\x{1389}
\def\y{729}
\begin{tikzgraphics}{0.55\textwidth}{\x}{\y}{splitting_curve_time_alpha_zoom_all}
\pxcoordinate{1.0*\x}{0.43*\y}{A}; \draw (A) node[right] {$\tvor(x_0)$};
\pxcoordinate{1.0*\x}{0.675*\y}{A}; \draw (A) node[right] {$\tinj(x_0)$};
\pxcoordinate{0.9*\x}{0.32*\y}{A}; \draw (A) node[right] {\textbf{C}};
\pxcoordinate{0.9*\x}{0.55*\y}{A}; \draw (A) node[right] {\textbf{B}};
\pxcoordinate{0.9*\x}{0.77*\y}{A}; \draw (A) node[right] {\textbf{A}};
\end{tikzgraphics}
\vspace{-1em}
\caption{Graph of the function $t(\alpha_2)$ along $\Sigma_1$ (blue), $\Sigma_2$ (red) and $\Sigma_3$ (green) for times $t \le 3.2$.
The balls are of three types: A, B and C. See the description detailed in the step 3.}
\label{fig:splitting_curve}
\end{figure}

Each self-intersection of the wavefront $\WF(x_0, 3.5)$ is a point of a splitting curve. We denote by $\Sigma_1$, $\Sigma_2$
and $\Sigma_3$ the three splitting curves associated respectively to the self-intersections $1$, $2$ and $3$, and
we have
$\Sigma_1 \cup \Sigma_2 \cup \Sigma_3 \subset \Split(x_0)$.
On Figure~\ref{fig:splitting_curve} is represented the graph of the function $\alpha_2 \mapsto t(\alpha_2)$ along the three
splitting curves for times $t \le 3.2$, the homotopic parameter $\alpha_2$ being strictly monotone along the splitting curves.
Since $\Sigma_1$ is below $\Sigma_2$ and $\Sigma_3$ on the figure,\footnote{We mean that the graph of $\alpha_2 \mapsto t(\alpha_2)$
associated to $\Sigma_1$ is below the ones associated to $\Sigma_2$ and $\Sigma_3$.} it is clear that only $\Sigma_1$ has
a chance to be part of the cut locus, compared to $\Sigma_2$ and $\Sigma_3$.
Note that there may exist others splitting curves but the numerical experiments suggest that $\Sigma_1$ is below any of them.

\begin{rmrk}
Note that the extremities of a splitting curve are the vortex and infinity according to the numerical experiments.
Hence, we cannot compute all the curve. Numerically, we use an option of the \hampath\ code to stop the homotopy if the
curve goes out the annulus of smaller circle of radius $0.005$ and of larger one of radius $100$.
\end{rmrk}

\begin{figure}[ht!]
\centering
\def\sizefig{0.27}
\includegraphics[width=\sizefig\textwidth]{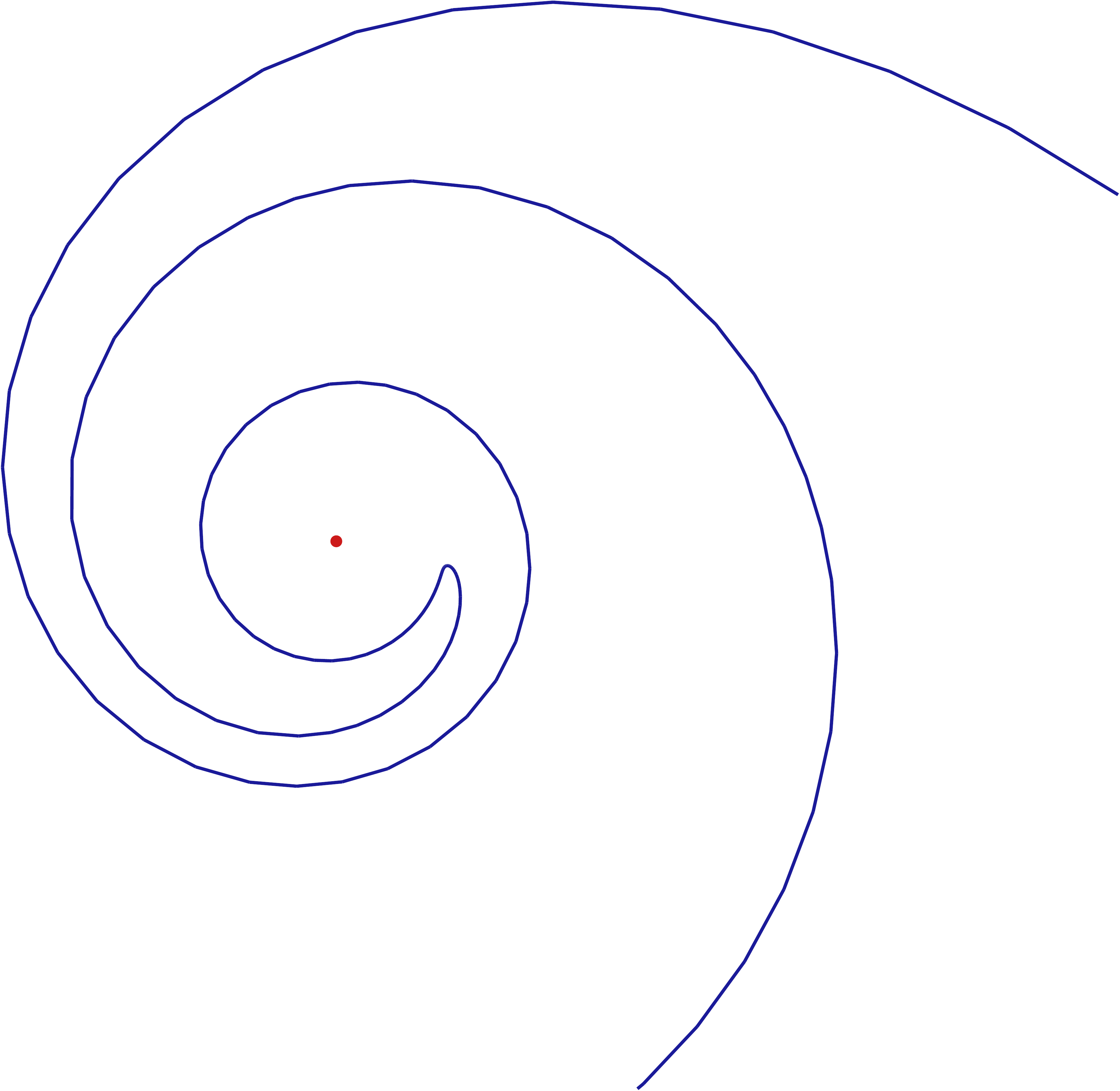}
\hspace{0.0em}
\includegraphics[width=\sizefig\textwidth]{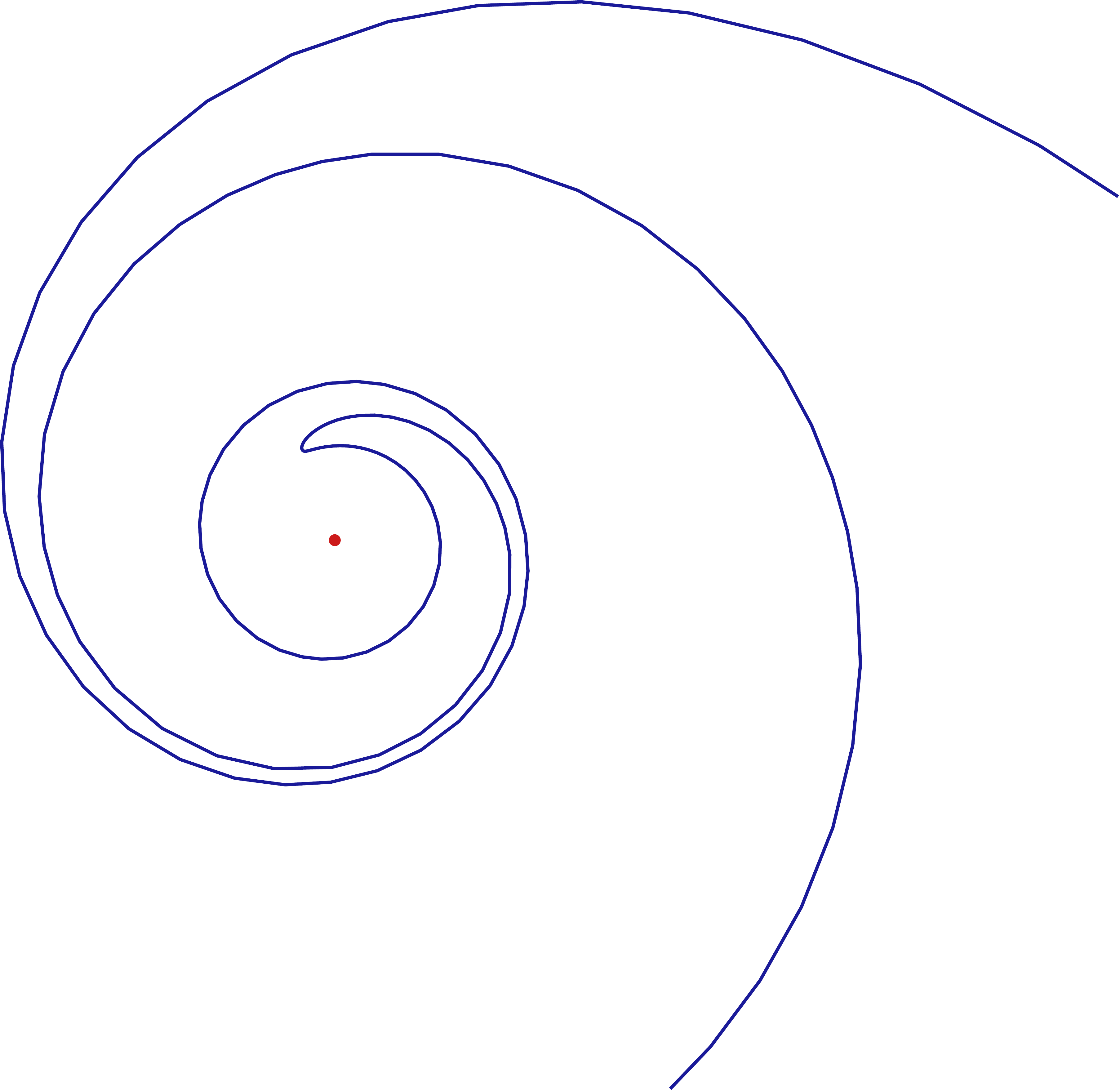}
\hspace{0.0em}
\includegraphics[width=\sizefig\textwidth]{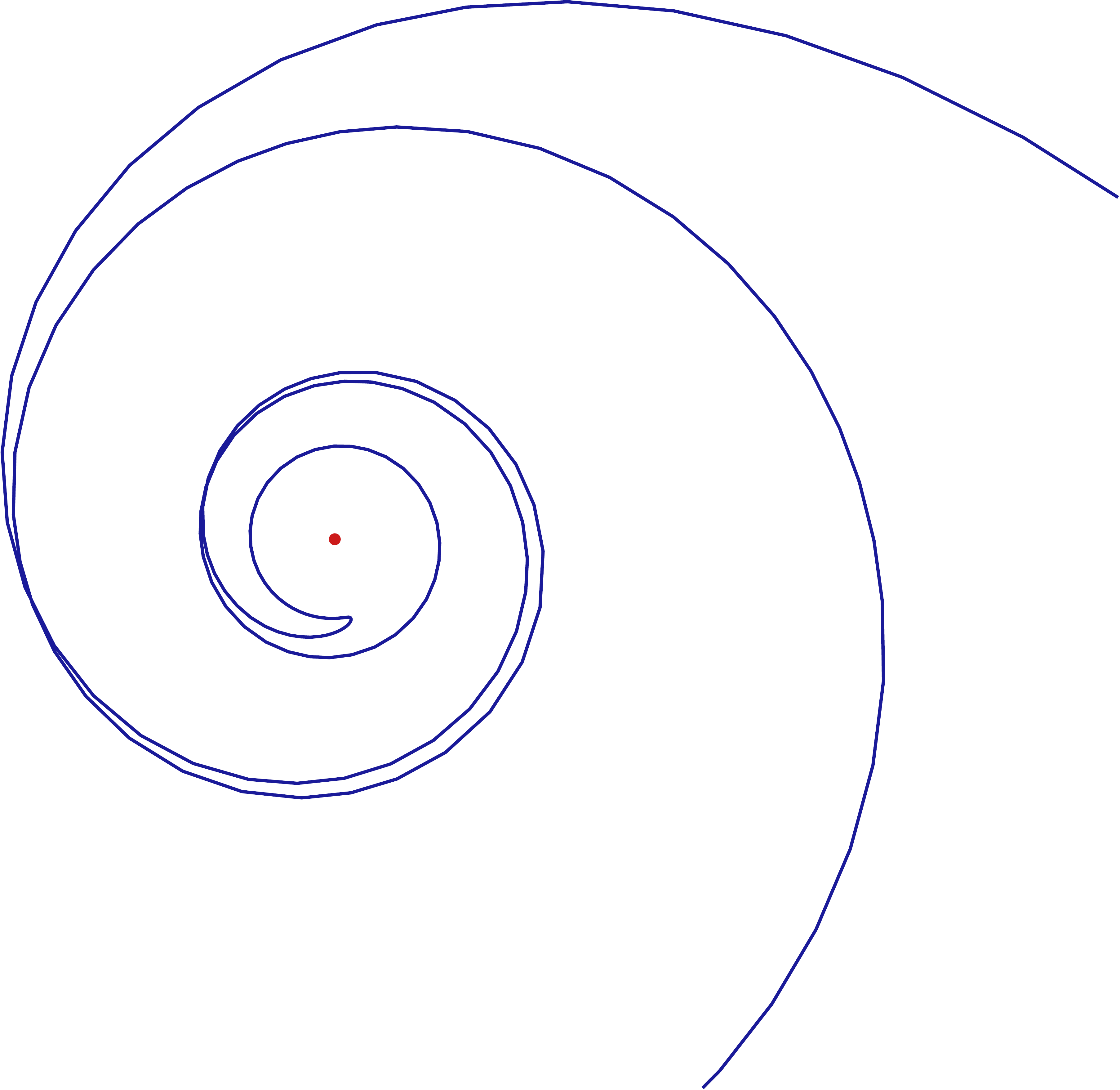}
\caption{Three wavefronts with the same scaling and a zoom around the vortex (red dot) for $t \in \{2.86, 2.88, 2.9\}$. 
For $t = 2.9 > \tinj(x_0) \approx 2.889$, the wavefront has two self-intersections.}
\label{fig:wavefronts_zoom}
\end{figure}

Let us describe the evolution of the wavefronts according to Figure~\ref{fig:splitting_curve}. For $t=0$, the wavefront is reduced
to the singleton $\{x_0\}$. For times $t \in \intervalleoo{0}{\tinj(x_0)}$ with $\tinj(x_0) = \inf t(\cdot)$ along $\Sigma_1$,
the wavefronts have no self-intersections. For $t = \tinj(x_0)$, the wavefront has one self-intersection, while for 
$t \in \intervalleoo{\tinj(x_0)}{\tvor(x_0)}$, where $\tvor(x_0) \coloneqq \norme{x_0}$ is the minimal time to reach the vortex,
the wavefronts have two self-intersections contained in $\Sigma_1$, see Figure~\ref{fig:wavefronts_zoom}, and up to
six self-intersections in $\Sigma_1 \cup \Sigma_2 \cup \Sigma_3$. Finally, for times $t \ge \tvor(x_0)$, the wavefronts have one
self-intersection in $\Sigma_1$ and three in $\Sigma_1 \cup \Sigma_2 \cup \Sigma_3$, see the right subgraph of Figure~\ref{fig:wavefronts}.



\medskip
\paragraph{\textbf{Step 2: Computation of the cut locus.}} 
According to Proposition \ref{prop:continuity_weak_drift}, at a final point $x_f$ where the drift is weak the mapping
$V_\mu(x_0, \cdot)$ is continuous. In this case, from Proposition \ref{prop:caracterisation_cut_points}, if $x_f$ is a cut point,
then either it is a conjugate point or a splitting point. Since, there are no conjugate points from Conjecture \ref{conj:no_conj},
$x_f$ is a splitting point. Finally, from step 1, we can conclude that if $x_f$ is a cut point then it belongs to $\Sigma_1$.
In this section, we consider moreover that the drift is weak at the initial condition $x_0$. In this case, we conjecture the following:

\begin{conjecture}
If the drift is weak at $x_0$, then $V_\mu(x_0, \cdot)$ is continuous at any point of $\xdom$.
\end{conjecture}

This conjecture means that the characterization of cut points from Proposition \ref{prop:caracterisation_cut_points} holds in all the
state space $\xdom$. Hence, we have $\Cut(x_0) \subset \Sigma_1$. Now, by the converse part of Proposition
\ref{prop:caracterisation_cut_points} and according to the previous conjecture, we can conclude that:
\[
\Cut(x_0) = \Sigma_1.
\]

\medskip
\paragraph{\textbf{Step 3. Computations of the spheres and balls.}}

\begin{figure}[ht!]
\centering
\def\sizefig{0.3}
\includegraphics[width=\sizefig\textwidth]{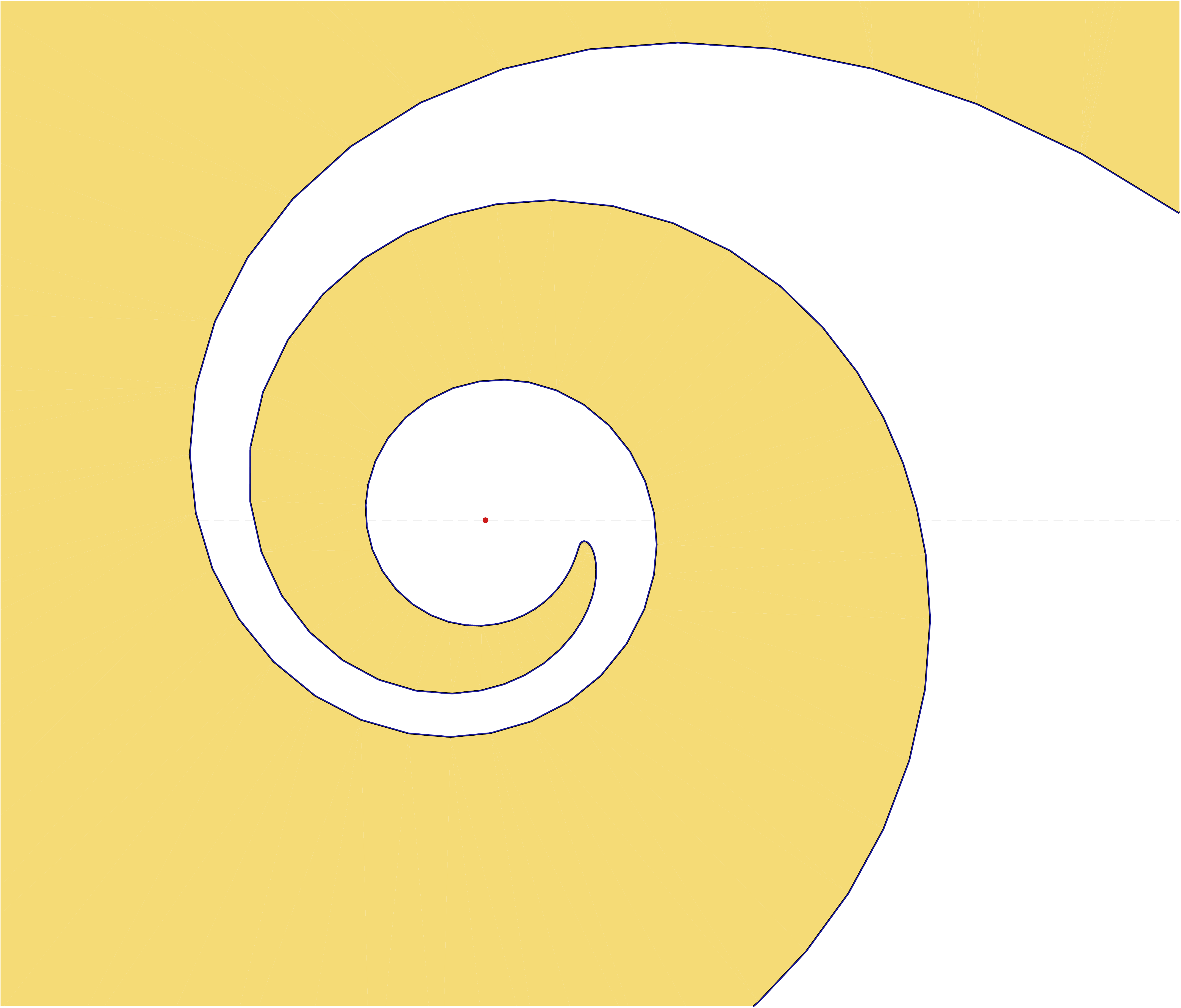}
\hspace{0.0em}
\includegraphics[width=\sizefig\textwidth]{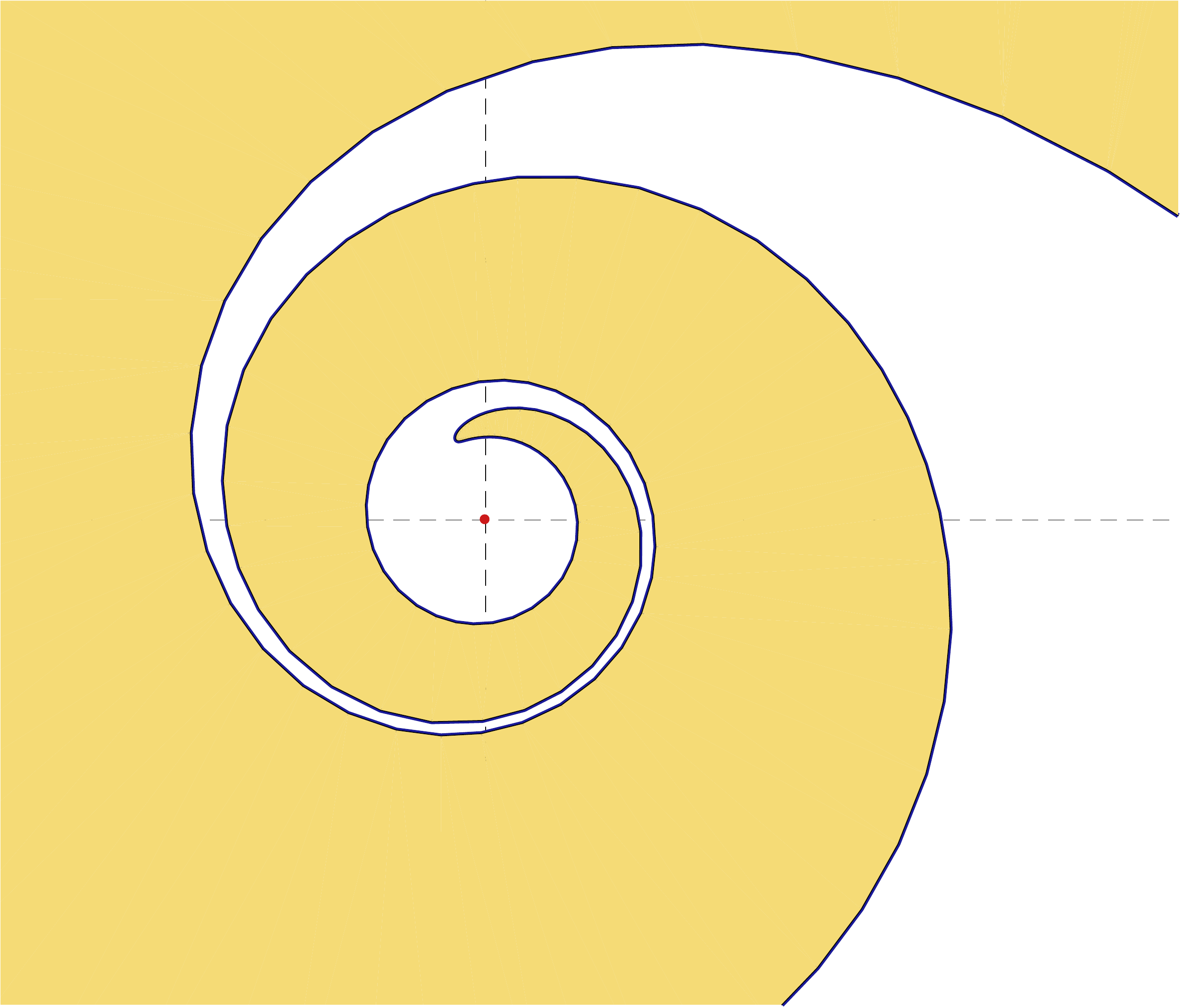}
\hspace{0.0em}
\includegraphics[width=\sizefig\textwidth]{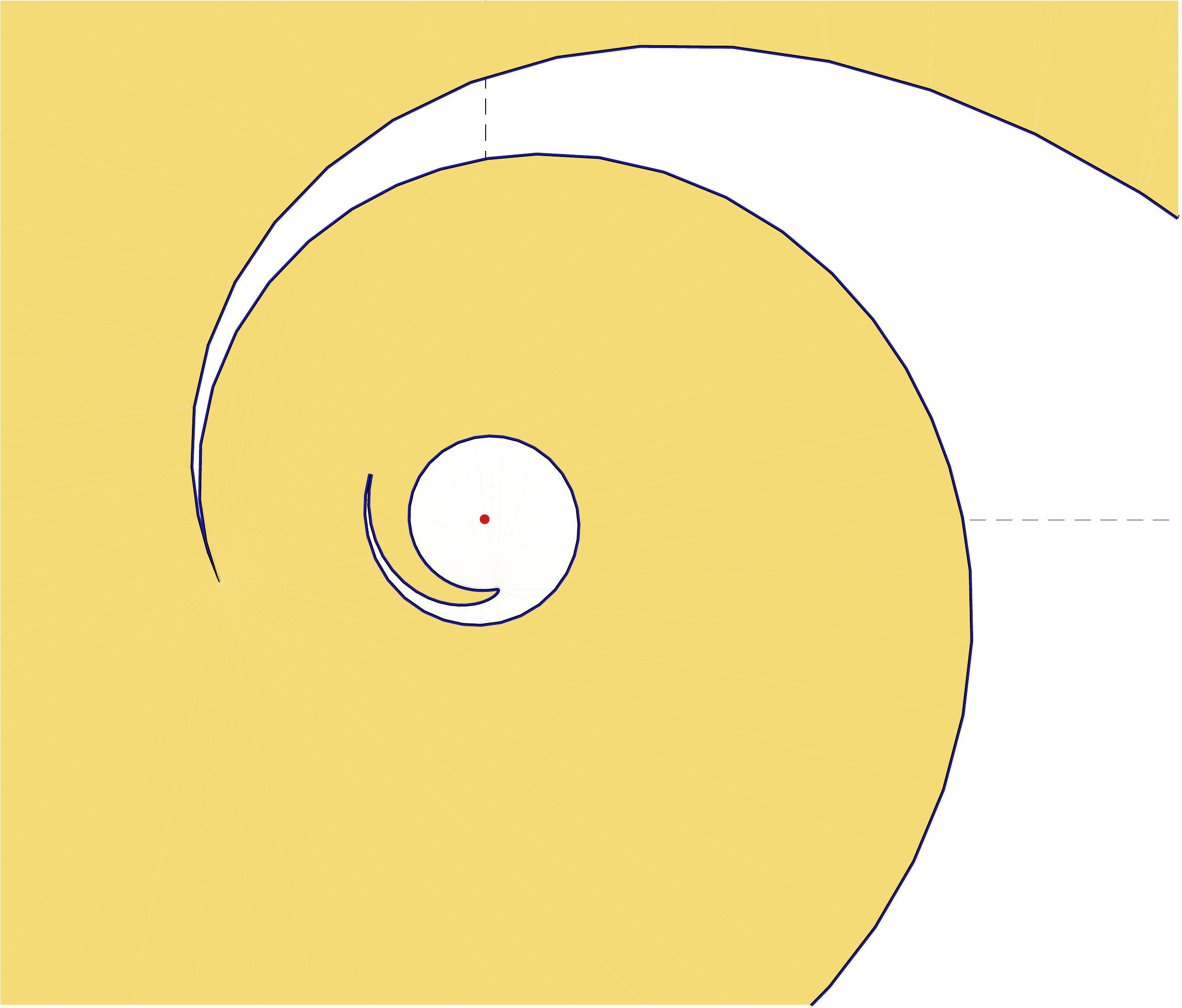}
\caption{Three balls with a zoom around the vortex (red dot) at times $t \in \{2.86, 2.88, 2.9\}$. 
Compare to Figure~\ref{fig:wavefronts_zoom}. 
For $t = 2.9$ the ball $\MBallClosed(x_0, t)$ is of type C and has a small hole with a fetus shape.}
\label{fig:balls_zoom}
\end{figure}

Figure~\ref{fig:spheres_balls} represents four balls that we can group into three categories (see also Figure~\ref{fig:splitting_curve}):
\begin{itemize}
\item \textbf{Type A.} The ball is simply connected in $\R^2$ with a smooth boundary;
\item \textbf{Type B.} The  ball is not simply connected in $\R^2$;
\item \textbf{Type C.} The ball is simply connected in $\R^2$ with a singularity on its boundary.
\end{itemize}
For times $t \in \intervalleof{0}{\tinj(x_0)}$, we have $\WF(x_0, t) = \MSphere(x_0, t)$ and the ball $\MBallClosed(x_0, t)$
is of type A. The two balls on the top of Figure~\ref{fig:spheres_balls} are of type A. 
For times $t \in \intervalleoo{\tinj(x_0)}{\tvor(x_0)}$, the balls have a hole around the vortex and so they are of type B.
The ball at the bottom-left subgraph of Figure~\ref{fig:spheres_balls} is of type B. One can see on Figure~\ref{fig:balls_zoom},
the evolution of the balls around the vortex for times close to $\tinj(x_0)$. The hole appears when the ball self-intersects,
like a snake biting its tail. For times $t \ge \tvor(x_0)$, the balls are simply connected in $\R^2$ since the hole has
vanished (the vortex is reached at $t=\tvor(x_0)$) and the spheres have a singularity at the cut point. In this case the balls have
a shape of an apple and are of type C.
The bottom-right subgraph of Figure~\ref{fig:spheres_balls} represents a ball of type C.

\begin{figure}[ht!]
\centering
\def\sizefig{0.44}
\def\x{1604}
\def\y{1170}
\begin{tikzgraphics}{\sizefig\textwidth}{\x}{\y}{Balls_1p500}
\pxcoordinate{0.1*\x}{0.1*\y}{A}; \draw (A) node[right] {\textbf{A}};
\end{tikzgraphics}
\hspace{0.0em}
\begin{tikzgraphics}{\sizefig\textwidth}{\x}{\y}{Balls_2p800}
\pxcoordinate{0.1*\x}{0.1*\y}{A}; \draw (A) node[right] {\textbf{A}};
\end{tikzgraphics}

\vspace{1.0em}
\begin{tikzgraphics}{\sizefig\textwidth}{\x}{\y}{Balls_2p900}
\pxcoordinate{0.1*\x}{0.1*\y}{A}; \draw (A) node[right] {\textbf{B}};
\end{tikzgraphics}    
\hspace{0.0em}
\begin{tikzgraphics}{\sizefig\textwidth}{\x}{\y}{Balls_3p500}
\pxcoordinate{0.1*\x}{0.1*\y}{A}; \draw (A) node[right] {\textbf{C}};
\end{tikzgraphics}
\caption{Four balls in the $(x_1,x_2)$-plane at times $t \in \{ 1.5, 2.8, 2.9, 3.5\}$ with the initial condition $x_0$
represented by a black dot and the vortex by a red dot. The balls $\MBallClosed(x_0, 1.5)$ and $\MBallClosed(x_0, 2.8)$
on top are of type A. The ball $\MBallClosed(x_0, 2.9)$ (Bottom-Left) is of type B while $\MBallClosed(x_0, 3.5)$
(Bottom-Right) is of type C.}
\label{fig:spheres_balls}
\end{figure}

\begin{figure}[ht!]
\centering
\includegraphics[width=0.45\textwidth]{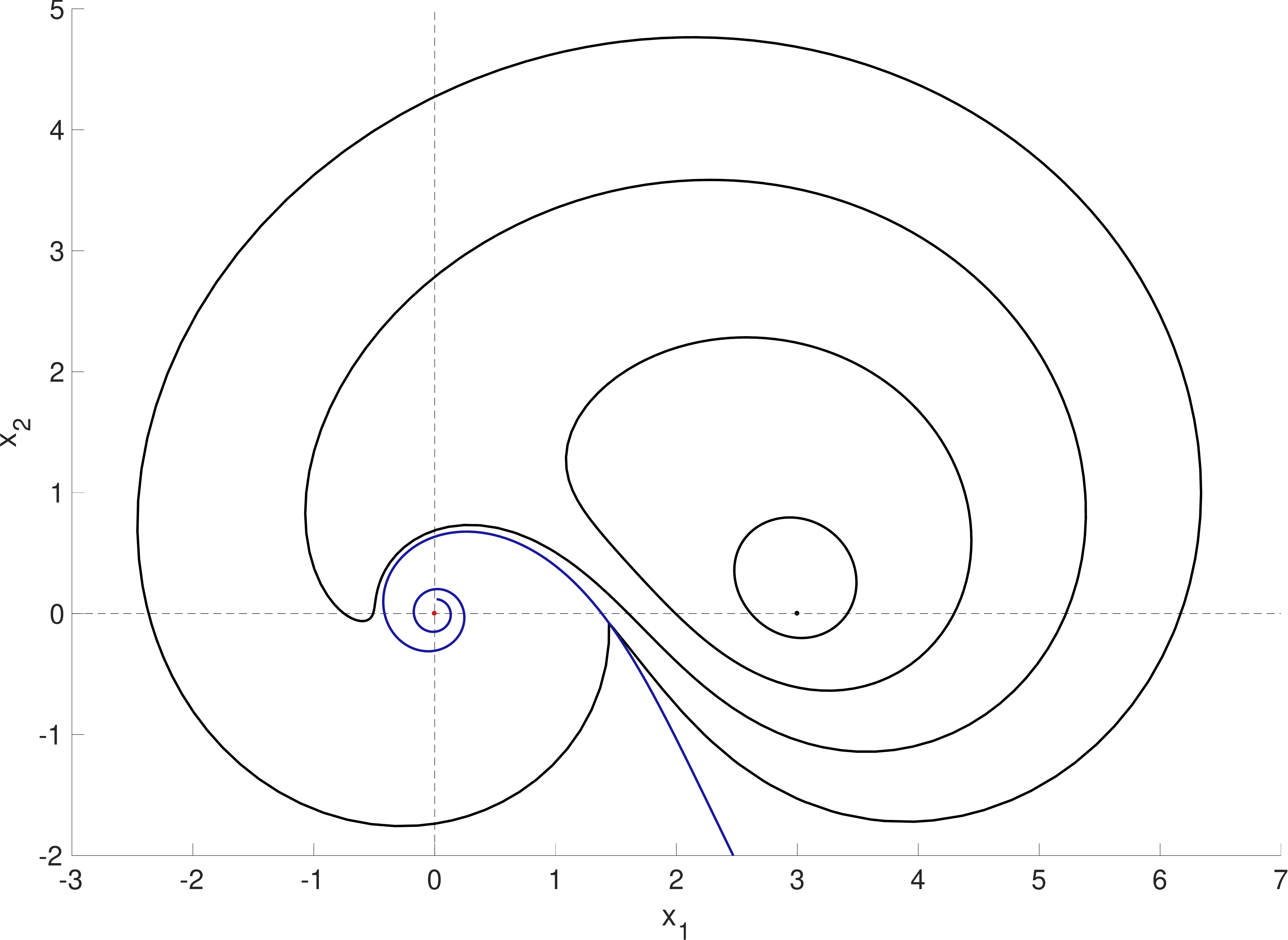}

\caption{Spheres $\MSphere(x_0, t)$ at times $t \in \{0.5, 1.5, 2.5, 3.5\}$ in black together with the cut locus $\Sigma_1$
in blue. The cut locus should go to the vortex but to gain in clarity, it is represented up to a distance of 0.12 to the
vortex. The vortex being represented by a red dot while the initial condition is represented by a black dot.}
\label{fig:spheres_and_cut}
\end{figure}

\begin{figure}[ht!]
\centering
\includegraphics[width=0.45\textwidth]{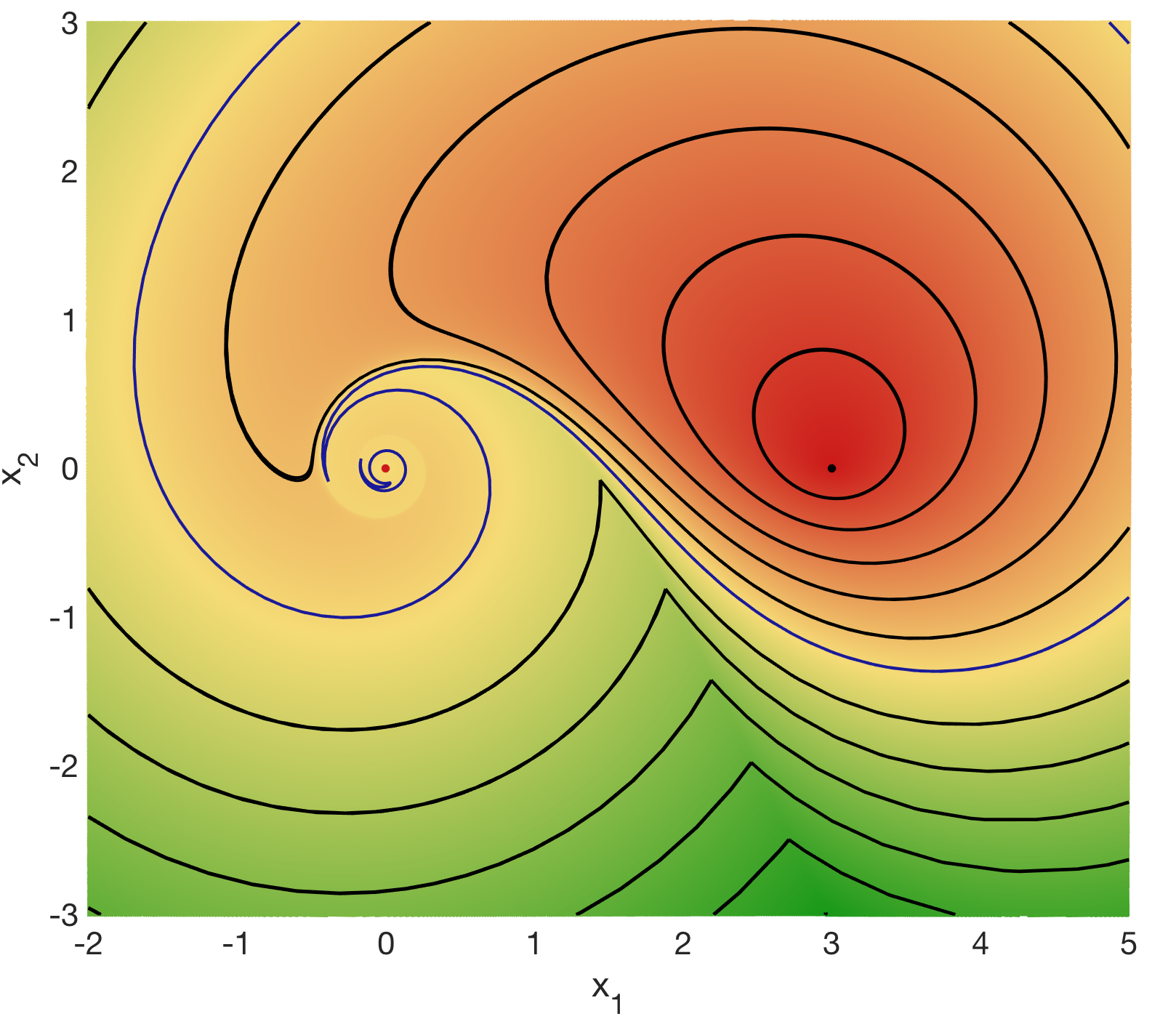}

\caption{Optimal synthesis. The black curves represent the spheres of types A and C. The two blue curves is a single sphere of type
B. Since the balls of type B are not simply connected and have one hole, the corresponding spheres have two connected components.}
\label{fig:synthesis}
\end{figure}

One can see on Figure~\ref{fig:spheres_and_cut} the evolution of the spheres $\MSphere(x_0, t)$ at times 
$t \in \{0.5, 1.5, 2.5, 3.5\}$ together with the cut locus $\Sigma_1$. One can notice that the singularity on $\MSphere(x_0, 3.5)$
in contained in the cut locus. Putting all together, we have constructed the optimal synthesis which is given on 
Figure~\ref{fig:synthesis}.

\section{Conclusion}

In this article, we have analyzed the Zermelo navigation problem with a vortex singularity. We have shown the existence of an optimal
solution. Thanks to the integrability properties we made a micro-local classification of the extremal solutions, showing remarkably
the existence of a Reeb foliation. The spheres and balls are described in the case where at the initial condition the current is weak.
This gives us the optimal synthesis for any initial condition such that $\norme{x_0} > \abs{\mu}$.

Note that the current is not necessarily weak all along the geodesics and so the situation we have analyzed is more general than the
Randers case from Finsler geometry. Still, this analysis has to be completed to the case of a strong current at the initial condition,
in particular in relation with the abnormal extremals.
Note that the historical work \cite{Caratheodory} is a case study of the strong current situation which can be applied to our study.
A possible further extension concerns the case of several vortices by analogy with the N-body problem.

Finally in a forthcoming article, the 2D-Zermelo navigation problem will be generalized to the case of arbitrary current and metric, 
with the symmetry of revolution. It will cover the historical problem and the vortex case. The analysis will be based on tools from
Hamiltonian dynamics from \cite{BolsinovFomenko:2004} and in particular the concept of Reeb graph. 

\begin{acknowledgement}
We are grateful to Carlos Balsa for helpful discussions about vortex theory and for drawing authors' attention to this topic,
and to Ulysse Serres for helpful discussions about Zermelo-like navigation problems.
\end{acknowledgement}

\bibliographystyle{model1-num-names}


\end{document}